\documentclass[10pt]{article}
\usepackage[utf8]{inputenc}
\usepackage[a4paper,top=2.5cm,bottom=2cm,left=1.75cm,right=1.75cm,marginparwidth=2cm]{geometry}
\usepackage{amsmath}
\usepackage{amssymb}
\numberwithin{equation}{section}
\usepackage[english]{babel}
\usepackage{amsthm}
\usepackage{bbm}
\usepackage{amsfonts}
\usepackage{mathrsfs}
\usepackage{array}
\usepackage{hyperref}
\usepackage{titlesec}
\usepackage{subcaption}
\usepackage{booktabs}

\hypersetup{
    colorlinks=true,
    linkcolor=blue,
    citecolor=red,
    filecolor=magenta,      
    urlcolor=black,
    pdftitle={ALR},
    pdfpagemode=FullScreen,
    }
\usepackage{graphicx, tikz}
\usepackage{float}

\usepackage{xcolor}
\usepackage{bbm}
\usepackage[format=plain,
            font=it]{caption}

\newtheorem{theorem}{Theorem}[section]
\newtheorem{lemma}[theorem]{Lemma}
\newtheorem{proposition}[theorem]{Proposition}
\newtheorem{remark}[theorem]{Remark}
\newtheorem{definition}[theorem]{Definition}
\theoremstyle{definition}
\newtheorem*{notation}{Notation}
\theoremstyle{plain}
\numberwithin{equation}{section}
\newcommand{\HH}{\ell^2(\mathbb N)\oplus\ell^2(\mathbb N)}
\newcommand{\TT}{\mathcal T}
\newcommand{\OO}{\mathcal O}

               \def\D{\Delta }

                 \def\n{\nu}         
                     
     \def\t{\tau }

\usepackage{todonotes}

\title{Phase transitions and sharp localisation in anomalous localised resonance}

\author{Konstantinos Alexopoulos$^1$}
\date{}
\begin{document}

\maketitle

\vspace{-0.5cm}

{\footnotesize
$^1$ MAP5, Université Paris Cité, CNRS, Paris, France.
}

\begin{abstract}
We study how detuning and absorption determine the strength and spatial extent of anomalous localised resonance in a two-dimensional quasistatic coated structure with shell permittivity $-1+c\ell^p+i\ell$, where $\ell>0$ is the loss and $c\ne0$, $p>0$ are fixed. Detuning can carry the material path through narrow resonances, producing intermittent rather than sustained power growth. For an explicit family of noncircular shells and a smooth exterior dipole source, both fixed, we classify the dissipated power under a quantitative smallness condition on the deformation. Depending on the detuning exponent and conformal source distance, the power diverges along the entire path, exhibits divergent peaks with vanishing gaps, or decays, including the threshold cases. We also identify the sharp conformal boundary beyond which the field, normalized by the square root of the dissipated power, converges locally to zero along every divergent-power sequence. For resonant sources, tangential reference crossings yield an explicit nonzero limit for the exterior gradient energy divided by the dissipated power at the critical radius. 
\end{abstract}

\section{Introduction}
Plasmonic structures can concentrate electromagnetic fields near their interfaces through the excitation of surface modes. In the quasistatic regime, anomalous localised resonance describes a singular form of this response: as material loss vanishes, increasingly fine spatial oscillations produce large fields in a bounded region, while the response farther away remains controlled. This separation between local amplification and exterior behavior underlies its connection with cloaking and superlensing \cite{milton,review}. Related near-field enhancement has been observed experimentally through plasmon-assisted transmission across silver films \cite{fangtransmission} and sub-diffraction-limited imaging with a silver superlens \cite{fangsuperlens}. These experiments demonstrate the physical relevance of surface-mode amplification and its sensitivity to absorption, rather than the zero-loss limit studied here.

In this paper, we consider a two-dimensional coated structure whose normalized core and exterior permittivities are equal to one, while the shell permittivity is $s_\ell=-1+c\ell^p+i\ell$, with fixed $c\ne0$ and $p>0$. The positive imaginary part $\ell$ represents material loss and enters the dissipated power through the absorption density $\ell|\nabla u_\ell|^2$, where $u_\ell$ denotes the complex electrostatic potential and $-\nabla u_\ell$ is the corresponding electric field. The real increment $c\ell^p$ represents detuning from the reference value $-1$. Thus, $p$ controls whether detuning is smaller than, comparable to, or larger than loss. We call the curve $\ell\mapsto s_\ell$ in the complex permittivity plane the \emph{material path}. As $\ell$ decreases to zero, this curve approaches $-1$ from the upper half-plane, with $c$ fixing the sign of the real detuning and $p$ determining the manner of approach. In particular, for $0<p<1$ the approach is tangential to the real axis. Physically, this path describes a family of shell materials with progressively weaker absorption, whose real permittivity is simultaneously tuned towards $-1$, while the geometry and source remain fixed. It specifies how closely the shell is tuned as its damping is reduced: the real part controls proximity to the surface-mode resonances, whereas the imaginary part limits their amplification. Thus, the exponent $p$ specifies the balance between these two changes, rather than a change in the shape of the coating. This is a prescribed relation between material parameters, not necessarily a frequency sweep in a single material. Our question is how this relative rate determines both the growth of the dissipated power and the spatial extent of the resonant response for a prescribed source.

The dependence on the material path is particularly important when detuning dominates loss. As the shell approaches its limiting permittivity, the path can pass through narrow resonant windows separated by much weaker responses. Reducing loss need not, therefore, imply sustained power growth. The source position also matters, since it determines how strongly the increasingly fine modes are excited. Sensitivity to relative dissipation rates has been established for slab configurations \cite{meklachi}, while related competition between material dispersion and damping appears in \cite{alexopoulosdavies}. Here we vary the real detuning within the shell and retain a fixed source and geometry throughout the limit.

The mathematical analysis of this phenomenon has linked dissipated-power blow-up to the spectrum of a two-interface Neumann--Poincar\'e-type operator \cite{ammari}, with further developments concerning spectral accumulation, nonradial configurations and complementary media \cite{ando,kohn,nguyen}. We work with the noncircular family generated by $\Psi(w)=w+\varepsilon w^{-2}$ from the annulus $\rho<|w|<1$, and a smooth exterior dipole source. Under an explicit condition relating the fixed deformation to $\rho$ and $p$, we obtain a complete classification in terms of the detuning exponent and conformal source distance. The dissipated power either diverges along the entire loss path, exhibits divergent peaks separated by vanishing responses, or tends to zero. The classification includes the threshold cases and distinguishes persistent absorption from intermittent excitation of narrow resonances.

We also determine how far this amplification remains visible after the field is divided by the square root of the dissipated power. The sharp conformal observation radius is $r_c=\rho^{-1/(2\min(p,1))}$: beyond its image under $\Psi$, the normalized field converges locally to zero along every divergent-power sequence. In the resonant source range, this boundary cannot be moved closer to the shell. For tangential reference crossings, we further compute an explicit nonzero limit of the exterior gradient energy divided by the dissipated power at $r_c$. Local decay away from the critical curve, therefore, coexists with residual integrated energy near it. These conclusions concern normalized fields for a prescribed source, rather than invisibility of arbitrary objects.

Our coefficient formulation is motivated by the use of Faber polynomials and Grunsky coefficients in closely related transmission problems. In \cite{junglim} construct geometry-adapted harmonic bases are constructed and layer-potential operators are represented by matrices with Grunsky coefficients, while in \cite{choikimlim}, corresponding geometric multipole expansions are developped. Particularly close to our setting is the work carried in \cite{choilimshipman}, which treats the Neumann--Poincar\'e spectrum between two conformal level curves. We use this framework to encode the noncircular core response and the source in compatible coordinates. A related modal analysis of field concentration appears in \cite{alexopoulosdaviesmillien}, where spherical harmonics and radial Bessel functions are used to describe non-singular hotspots between high-index particles. Here, the concentration instead arises from vanishing loss in a fixed geometry. The difficulty here is to control a forced response as loss vanishes on one fixed shell, rather than to determine spectral locations alone. Geometric coupling must remain controlled even near increasingly narrow resonances. Uniform comparisons with a circular auxiliary problem resolve the power growth, and weighted comparisons retain the outward attenuation needed for sharp localisation. Our analysis is quasistatic. Finite-frequency resonance and limiting absorption involve separate issues addressed in \cite{nguyenfinite,onofrei,liliu,kettunen,nguyenlap}.

The paper is organised as follows. In Section~\ref{framework}, we formulate the transmission problem and state the main results. Sections~\ref{boundaryformulation} and~\ref{conformal} establish its boundary-integral and coefficient representations. In Sections~\ref{sec:uniform} and~\ref{sec:power}, we derive the uniform comparison and prove the power classification. Sections~\ref{sec:weighted} and~\ref{sec:localisation} establish the weighted comparison, sharp localisation and endpoint limit. Section~\ref{numerics} illustrates these results numerically. Detailed supporting proofs are collected in the appendices.

\begin{notation}
Let $A(\ell),B(\ell)\ge0$ for $\ell>0$. As $\ell\downarrow0$, we use the conventions
\begin{align*}
 A\lesssim B
 &\quad\Longleftrightarrow\quad
 \exists\,C,\ell_0>0:\quad A(\ell)\le C B(\ell)
 \quad(0<\ell<\ell_0),\\
 A\asymp B
 &\quad\Longleftrightarrow\quad A\lesssim B\ \text{and}\ B\lesssim A.
\end{align*}
Thus, $A\asymp B$ means that $c_0B(\ell)\le A(\ell)\le CB(\ell)$ for some $0<c_0\le C<+\infty$ and all sufficiently small $\ell>0$. The constants are independent of $\ell$ and, for estimates involving mode indices, uniform over the indices in the stated range. They may depend on fixed parameters and fixed observation radii. The same conventions apply to other limiting variables, including sequences as $n\to+\infty$. If $B$ is positive in the limiting regime, we write
\begin{align*}
 A\sim B\quad\Longleftrightarrow\quad
 \lim_{\ell\downarrow0}\frac{A(\ell)}{B(\ell)}=1.
\end{align*}
\end{notation}

\section{Problem and main results}\label{framework}

Let us start by formulating the transmission problem, specify the fixed geometry and source, and state the power and localisation results.

\subsection{Problem setting}
Let $\Omega_i$ and $\Omega_e$ be bounded simply connected domains in $\mathbb R^2$ with real-analytic Jordan boundaries, with $\overline{\Omega_i}\subset\Omega_e$. We write
\begin{align*}
 \Gamma_i=\partial\Omega_i,\qquad \Gamma_e=\partial\Omega_e \qquad \text{ and } \qquad
 \Omega_s=\Omega_e\setminus\overline{\Omega_i}.
\end{align*}
The unit normals $\nu_i$ and $\nu_e$ point out of the bounded domains $\Omega_i$ and $\Omega_e$, respectively. In particular, $\nu_i$ points \emph{into} the shell, and the outward shell normal on $\Gamma_i$ is $-\nu_i$. For each interface, the symbols $+$ and $-$ denote limits from the exterior and interior of its bounded domain, taken along this same normal.

Let $\ell>0$, $c\in\mathbb R\setminus\{0\}$, and $p>0$. Set
\begin{align*}
 s_\ell=-1+c\ell^p+i\ell \qquad \text{ and } \qquad
 a_\ell(x)=\begin{cases}
 1,&x\in\Omega_i\cup(\mathbb R^2\setminus\overline{\Omega_e}),\\
 s_\ell,&x\in\Omega_s.
 \end{cases}
\end{align*}
\begin{figure}
\centering
\includegraphics[width=.6\linewidth]{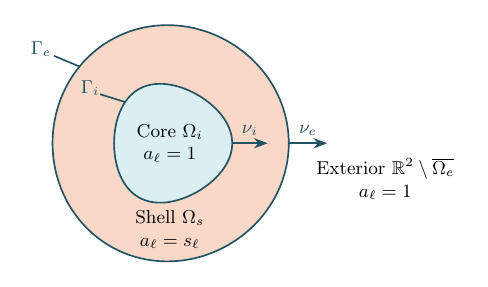}
\caption{The core $\Omega_i$, shell $\Omega_s$ and exterior, with their normalized relative permittivities. The interfaces are $\Gamma_i$ and $\Gamma_e$, and both normals point out of their bounded interiors.}\label{fig:domains}
\end{figure}
The function $a_\ell$ is the spatially varying relative permittivity, normalized by the positive exterior permittivity, and $s_\ell$ is its value in the shell. The core and exterior are lossless and have the same normalized permittivity. With the time convention $e^{-i\omega t}$, the imaginary part $\operatorname{Im}s_\ell=\ell>0$ describes absorption in the shell. We call $\ell$ the loss parameter and $d_\ell=\operatorname{Re}s_\ell+1=c\ell^p$ the real detuning, as it measures the displacement of the real shell permittivity from $-1$, the reference contrast in the quasistatic resonance model \cite{ammari}.

\begin{remark}
The exponent $p$ compares these two scales, since
\begin{align*}
 \lim_{\ell\downarrow0}\frac{|d_\ell|}{\ell}
 =\lim_{\ell\downarrow0}|c|\ell^{p-1}
 =
 \begin{cases}
 0,&p>1,\\
 |c|,&p=1,\\
 +\infty,&0<p<1.
 \end{cases}
\end{align*}
Thus, detuning is smaller than loss for $p>1$, comparable to loss for $p=1$, and larger than loss for $p<1$. In the last case the permittivity approaches the real axis tangentially, allowing the material path to resolve narrow resonant peaks.
\end{remark}

We consider a fixed real source $f\in C_c^\infty(\mathbb R^2\setminus\overline{\Omega_e})$, separated from $\Gamma_e$, satisfying $\int_{\mathbb R^2}f=0$. The complex electrostatic potential $u_\ell$, whose negative gradient is the electric field, solves
\begin{equation}
\left\{
\begin{aligned}
& \nabla\cdot(a_\ell\nabla u_\ell)=f\quad&&\text{in }\mathbb R^2,\\
&u_\ell(x)=O(|x|^{-1})\quad &&\text{ as }|x|\to+\infty.
 \label{pde}
\end{aligned}
\right.
\end{equation}
Equivalently, the governing transmission problem is
\begin{equation}
 \left\{\begin{aligned}
 &\Delta u_\ell=0&&\text{in }\Omega_i\cup\Omega_s,\\
 &\Delta u_\ell=f&&\text{in }\mathbb R^2\setminus\overline{\Omega_e},\\
 &u_\ell|_+-u_\ell|_-=0&&\text{on }\Gamma_i\cup\Gamma_e,\\
 &s_\ell\,\frac{\partial}{\partial \nu_i}u_\ell|_+=\frac{\partial}{\partial \nu_i}u_\ell|_-&&\text{on }\Gamma_i,\\
 &\frac{\partial}{\partial \nu_e}u_\ell|_+=s_\ell\,\frac{\partial}{\partial \nu_e}u_\ell|_-&&\text{on }\Gamma_e,\\
 &u_\ell(x)=O(|x|^{-1})&&\text{as }|x|\to+\infty.
 \end{aligned}\right.
 \label{transmission}
\end{equation}
The solution is required to belong to $H^1_{\rm loc}(\mathbb R^2)$ with $\nabla u_\ell\in L^2(\mathbb R^2)$. Since $s_\ell\ne0$, the shell equation in \eqref{transmission} is equivalent to the divergence-form equation there.

Let
\begin{equation}
 \Phi(x)=\frac1{2\pi}\log|x| \qquad \text{ and } \qquad
 F(x)=\int_{\mathbb R^2}\Phi(x-y)f(y)\,dy.
 \label{freepotential}
\end{equation}
Then, since $\int_{\mathbb{R}^2}f=0$, we have 
\begin{align*}
\Delta\Phi=\delta_0, \qquad \Delta F=f \qquad \text{and} \qquad F(x)=O(|x|^{-1}).
\end{align*}
The free field is harmonic in a neighborhood of $\overline{\Omega_e}$. We call $u_\ell-F$ the scattered field. We define the dissipated power $E(\ell)$ to be
\begin{equation}
 E(\ell)=\ell\int_{\Omega_s}|\nabla u_\ell|^2\,dx.
 \label{physicalpower}
\end{equation}
It is a measure of the absorption inside the material. Finally, for a bounded operator $A$, let us denote its spectrum by $\operatorname{Sp}(A)$.

\subsection{Fixed geometry, source and observation regions}\label{problemgeometry}
We specialize the problem to a fixed noncircular family. Fix $0<\rho<1$ and real $\varepsilon$ with $0<|\varepsilon|<\rho^3/16$, and set
\begin{align*}
 \begin{gathered}
 \Psi(w)=w+\varepsilon w^{-2},\qquad
 \Omega_s=\Psi(\{\rho<|w|<1\}),\qquad
 \Gamma_i=\Psi(\{|w|=\rho\}) \qquad \text{ and } \qquad
 \Gamma_e=\Psi(\{|w|=1\}).
 \end{gathered}
\end{align*}
The bounded interiors define $\Omega_i$ and $\Omega_e$. 

\begin{figure}
\centering
\includegraphics[width=.87\linewidth]{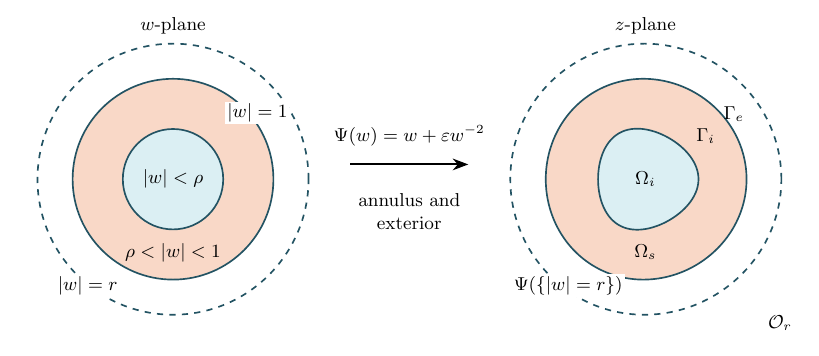}
\caption{The conformal annulus and its physical image under $\Psi$. Only the annulus and exterior are mapped as $\Psi$ is not used to parametrize the core interior.}\label{fig:conformal}
\end{figure}

For the remainder of the paper, fix $R>1$ and let $z_0=\Psi(R)$. Choose
$0<b<\operatorname{dist}(z_0,\overline{\Omega_e})$,
and let $\chi\in C_c^\infty(B_b(0))$ be a real-valued radial function satisfying $\int_{\mathbb R^2}\chi(x)\,dx=1$. We use the source
\begin{align*}
 f=\frac{\partial}{\partial x_1}\chi(\cdot-z_0)\qquad \text{ and } \qquad
 \alpha=\frac{2\log R}{\log(1/\rho)}>0.
\end{align*}
\begin{remark}
The parameter $\alpha$ measures the conformal source distance. Larger $\alpha$ means weaker excitation of high angular modes. 
\end{remark}
The source is held fixed as $\ell\to0$. For $r>1$, let us define
\begin{equation}
 \mathcal O_r=\Psi(\{|w|>r\}) \qquad \text{ and } \qquad
 S_r(\ell)=\int_{\mathcal O_r}|\nabla(u_\ell-F)|^2.
 \label{observables}
\end{equation}
The quantity $S_r$ measures the scattered field remaining outside the observation curve $\Psi(\{|w|=r\})$. 

\subsection{Main results}\label{results}
There are two distinct thresholds to identify. The first concerns whether the source transfers an unbounded amount of power into the shell as absorption vanishes. The second concerns how far the resulting field remains observable after unit-power normalization. The first theorem classifies the response and the second identifies its spatial extent.

We state the results in terms of the preceding transmission problem, fixed source and physical energies. In addition to the geometric bound $0<|\varepsilon|<\rho^3/16$, assume
\begin{equation}
 \tau=\min(p,1) \qquad \text{ and } \qquad
 \eta=\sqrt{2|\varepsilon|/\rho^3}<\rho^{1/\tau}.
 \label{standingcoupling}
\end{equation}
For a tangential approach, individual reference resonances become narrow compared with their spacing. Condition~\eqref{standingcoupling} compares the geometric coupling scale $\eta^n$ with the crossing loss scale $\rho^{n/\tau}$. It is a sufficient condition for the uniform estimates, not a necessary geometric threshold.

We refer to these fixed-geometry, fixed-source and material-path conditions as the standing hypotheses in the following sections. The constants implicit in $\asymp$ and $\lesssim$ may depend on $\rho,\varepsilon,R,c,p$, but not on loss. A nonzero deformation is chosen once. It does not shrink along a loss sequence.

\begin{theorem}[Complete fixed-source power classification]\label{mainpower}
Under the standing hypotheses:
\begin{enumerate}
\item If $p\ge1$, then, 
\begin{align}
    \left\{
    \begin{aligned}
        &\lim_{\ell\downarrow0}E(\ell)=+\infty, \qquad &&\text{ for } 0<\alpha\le1,\\
        &\lim_{\ell\downarrow0}E(\ell)=0, &&\text{ for } \alpha>1.
    \end{aligned}
    \right.
\end{align}
\item If $0<p<1$, the classification is
\begin{align}
    \left\{
    \begin{aligned}
        &\lim_{\ell\downarrow0}E(\ell) = +\infty, && \text{ for }0<\alpha\le2-1/p,\\
        &\limsup_{\ell\downarrow0}E(\ell) = +\infty, \quad\liminf_{\ell\downarrow0}E(\ell)=0, \qquad \qquad &&\text{ for } \max(0,2-1/p)<\alpha\le1/p,\\
        &\lim_{\ell\downarrow0}E(\ell)=0,&&\text{ for }\alpha>1/p.
    \end{aligned}
    \right.
\end{align}
The first range is empty when $p\le1/2$. Both positive threshold equalities are on the divergent side.
\end{enumerate}
\end{theorem}

\begin{theorem}[Sharp normalized conformal boundary]\label{mainlocal}
Under the standing hypotheses, set $r_c=\rho^{-1/(2\tau)}$. For every positive sequence $(\ell_j)$ with $\lim_{j\to+\infty}\ell_j=0$ and $\lim_{j\to+\infty}E(\ell_j)=+\infty$,
\begin{equation}
 \lim_{j\to+\infty}\frac{S_r(\ell_j)}{E(\ell_j)}=0
 \quad\text{for every fixed }r>r_c \qquad \text{ and } \qquad
 \lim_{j\to+\infty}\frac{u_{\ell_j}}{\sqrt{E(\ell_j)}}=0
 \quad\text{in }H^1_{\mathrm{loc}}(\OO_{r_c}).
 \label{mainupper}
\end{equation}
Suppose additionally that the source is in the resonant range: $0<\alpha\le1/p$ if $p<1$, or $0<\alpha\le1$ if $p\ge1$. Then, there is one such divergent-power sequence $(\ell_j)_{j\geq0}$ for which
\begin{align*}
 \lim_{j\to+\infty}\frac{S_r(\ell_j)}{E(\ell_j)}=+\infty, \quad \text{ for } \quad 1<r<r_c \qquad \text{ and } \qquad 
 \liminf_{j\to+\infty}\frac{S_{r_c}(\ell_j)}{E(\ell_j)}>0.
\end{align*}
Thus, in these source ranges, convergence in the full exterior gradient seminorm along every divergent-power sequence holds exactly for $r>r_c$.
\end{theorem}

The sharpness statement distinguishes a decay estimate from an optimal observation boundary. The following theorem gives the quantitative response at that boundary on tangential paths.

\begin{theorem}[Tangential endpoint energy]\label{mainendpoint}
Under the standing hypotheses, let $0<p<1$ and $0<\alpha\le1/p$. For all sufficiently large $n$, let $\ell_n$ be the unique sufficiently small positive loss satisfying the exact crossing equation
\begin{equation}
 |\operatorname{Re}\zeta(\ell_n)|=\rho^n \qquad \text{ and } \qquad
 \zeta(\ell)=\frac{c\ell^p+i\ell}{2-c\ell^p-i\ell}
 \label{crossings}
\end{equation}
Then, $(\ell_n)$ is decreasing with $\lim_{n\to+\infty}\ell_n=0$, and these crossings satisfy
\begin{equation}
 \lim_{n\to+\infty}\frac{S_{r_c}(\ell_n)}{E(\ell_n)}
 =\frac12\left(\frac{|c|}{2}\right)^{1/p}.
 \label{endpoint}
\end{equation}
\end{theorem}

\begin{remark}[Intermittent resonance]
In the intermittent regime, decreasing the loss does not produce a monotone increase in power: the material path passes through increasingly narrow resonant windows with vanishing responses between them. The localisation theorem concerns the divergent subsequences, or the entire path when the power diverges throughout.
\end{remark}

\begin{remark}[Physical interpretation and scope]
These are statements about a prescribed smooth source normalized by its dissipated power, not a bound on the unnormalized exterior field or invisibility of arbitrary inserted objects. The hypotheses permit a nonzero fixed deformation for every fixed $p>0$, but not one fixed deformation uniformly for all $p\downarrow0$.
\end{remark}

\section{Boundary-integral formulation}\label{boundaryformulation}
The layer representation separates the prescribed free field from the interface response and establishes well-posedness before we pass to coefficient coordinates.

\subsection{Layer potentials and trace spaces}
We represent the response of each permittivity interface by a layer potential whose density is the induced normal derivative jump. 
For a smooth closed curve $\Gamma$, let $H^{1/2}(\Gamma)$ be the Dirichlet trace space and $H^{-1/2}(\Gamma)$ its dual. Write
\begin{align*}
 H^{-1/2}_0(\Gamma)=\{\varphi\in H^{-1/2}(\Gamma):
                  \langle\varphi,1\rangle_\Gamma=0\}
 \qquad \text{ and } \qquad
 \mathcal X=H^{-1/2}_0(\Gamma_i)\oplus H^{-1/2}_0(\Gamma_e).
\end{align*}
\begin{remark}
The pairing with $1$ represents the total charge.
\end{remark}

\begin{definition}[Layer potentials and boundary operators]\label{layers}
For $j\in\{i,e\}$, define the single- and double-layer potentials by
\begin{align}
 \mathcal S_j[\varphi](x)&=\int_{\Gamma_j}\Phi(x-y)\varphi(y)\,ds_y,
       &&x\notin\Gamma_j,\label{single}\\
 \mathcal D_j[\psi](x)&=\int_{\Gamma_j}\frac{\partial}{\partial \nu_j(y)}\Phi(x-y)\psi(y)\,ds_y,
       &&x\notin\Gamma_j.\label{double}
\end{align}
Their boundary operators are
\begin{align}
 \mathcal K_j[\psi](x)
  &=\frac1{2\pi}\operatorname{p.v.}\int_{\Gamma_j}
     \frac{(y-x)\cdot\nu_j(y)}{|x-y|^2}\psi(y)\,ds_y,\label{Kdef}\\
 \mathcal K_j^*[\varphi](x)
  &=\frac1{2\pi}\operatorname{p.v.}\int_{\Gamma_j}
     \frac{(x-y)\cdot\nu_j(x)}{|x-y|^2}\varphi(y)\,ds_y,
       &&x\in\Gamma_j.\label{Kstardef}
\end{align}
Here $\mathcal K_j:H^{1/2}(\Gamma_j)\to H^{1/2}(\Gamma_j)$ and
$\mathcal K_j^*:H^{-1/2}(\Gamma_j)\to H^{-1/2}(\Gamma_j)$ are dual operators. For $k\ne j$, define the cross-interface operators $\mathcal S_{kj}$ and $\mathcal N_{kj}$ by
\begin{align}
\begin{aligned}
 \mathcal S_{kj}:H^{-1/2}(\Gamma_j)&\to H^{1/2}(\Gamma_k)\\
 \varphi &\mapsto \mathcal S_j[\varphi]\Big |_{\Gamma_k}
\end{aligned}
\end{align}
and
\begin{align}
\begin{aligned}
 \mathcal N_{kj}:H^{-1/2}(\Gamma_j)&\to H^{-1/2}(\Gamma_k)\\
 \varphi &\mapsto \frac{\partial}{\partial \nu_k}\mathcal S_j[\varphi]\Big|_{\Gamma_k},
\end{aligned}
\end{align}
respectively. The first index is the observation interface and the second is the source interface. We denote by $\mathcal S_{jj}[\varphi]$ the common interior and exterior boundary trace of $\mathcal S_j[\varphi]$ on $\Gamma_j$.
\end{definition}

\begin{theorem}[Layer-potential jump relations]\label{jumpthm}
For $\varphi\in H^{-1/2}(\Gamma_j)$ and $\psi\in H^{1/2}(\Gamma_j)$, the one-sided traces satisfy
\begin{align}
 \mathcal S_j[\varphi]|_+&=\mathcal S_j[\varphi]|_-=\mathcal S_{jj}[\varphi],\label{Sjump}\\
 \frac{\partial}{\partial \nu_j}\mathcal S_j[\varphi]|_\pm
   &=\left(\pm\tfrac12 I+\mathcal K_j^*\right)\varphi,\label{normaljump}\\
 \mathcal D_j[\psi]|_\pm
   &=\left(\mp\tfrac12 I+\mathcal K_j\right)\psi.\label{Djump}
\end{align}
The single-layer traces and double-layer traces are taken in $H^{1/2}(\Gamma_j)$, and the normal traces in $H^{-1/2}(\Gamma_j)$.
\end{theorem}
\begin{proof}
The jump identities with the positive logarithmic fundamental solution and outward normals are given in \cite[Section 2]{ammari}. They apply first to smooth densities. Continuity of the layer trace maps and density of smooth functions in the stated Sobolev spaces yield the asserted identities.
\end{proof}

\begin{remark}
	\begin{enumerate}
		\item The diagonal $\mathcal K_j^*$ are compact on $H^{-1/2}(\Gamma_j)$ for the analytic curves considered here.
		\item The cross-interface maps are \emph{smoothing}: their kernels have no singularity because the two curves are disjoint, so for every $m\ge0$, $\mathcal S_{kj},\mathcal N_{kj}:H^{-1/2}(\Gamma_j)\longrightarrow H^m(\Gamma_k)$, with $k\ne j$, are bounded. Thus, even a distributional boundary density produces a smooth field on the other interface.
		\item The single layers with zero charge have $O(|x|^{-1})$ decay and finite global Dirichlet energy.
	\end{enumerate}
\end{remark}

\subsection{Integral formulation of the transmission problem}
The two flux conditions become a coupled boundary system. Its diagonal blocks describe each interface separately and its off-diagonal blocks describe their interaction. Positive loss gives uniqueness through an energy identity even though the real part of the shell coefficient is negative. We isolate the ingredients of this argument in the following lemmas.

For $\zeta_\ell=\zeta(\ell)$ defined in~\eqref{crossings}, set
\begin{equation}
 \lambda_\ell=\frac{s_\ell+1}{2(s_\ell-1)}
             =-\frac{\zeta_\ell}{2}.
 \label{lambdazeta}
\end{equation}
In particular $\operatorname{Im}\zeta_\ell>0$ and $\lambda_\ell\ne0$.
Define the two-interface operator $\mathbb K^*$ and the forcing vector $\mathbf g_F$ by
\begin{equation}
 \mathbb K^*=
 \begin{pmatrix}
 -\mathcal K_i^*&-\mathcal N_{ie}\\
 \mathcal N_{ei}&\mathcal K_e^*
 \end{pmatrix} \qquad \text{ and } \qquad
 \mathbf g_F=\binom{-\frac{\partial}{\partial \nu_i}F}{\frac{\partial}{\partial \nu_e}F},
 \label{blockNP}
\end{equation}
respectively. The sign in the first row reflects that $\nu_i$ is opposite to the outward shell normal. We do not assume that $\mathbb K^*$ is self-adjoint in the usual product $L^2$ inner product.

\begin{lemma}[Interface equations]\label{layerflux}
Let $\ell>0$ and $\boldsymbol\varphi\in\mathcal X$. The field $v=F+\mathcal S_i[\varphi_i]+\mathcal S_e[\varphi_e]$ is continuous across both interfaces. Its flux conditions are equivalent to $(\lambda_\ell I-\mathbb K^*)\boldsymbol\varphi=\mathbf g_F$.
\end{lemma}
\begin{proof}
See Appendix~\ref{proof:layerflux}.
\end{proof}

\begin{lemma}[Charge conservation and the Fredholm property]\label{layercharges}
The vector $\mathbf g_F$ belongs to $\mathcal X$. For arbitrary $\psi_j\in H^{-1/2}(\Gamma_j)$,
\begin{align*}
 \langle\mathcal K_j^*\psi_j,1\rangle=\tfrac12\langle\psi_j,1\rangle,\quad
 \langle\mathcal N_{ie}\psi_e,1\rangle=0 \quad \text{ and } \quad
 \langle\mathcal N_{ei}\psi_i,1\rangle=\langle\psi_i,1\rangle.
\end{align*}
Consequently $\mathbb K^*$ restricts to a compact operator on $\mathcal X$, and $\lambda_\ell I-\mathbb K^*$ is Fredholm of index zero.
\end{lemma}
\begin{proof}
See Appendix~\ref{proof:layercharges}.
\end{proof}

\begin{lemma}[Decay of zero-charge layers]\label{layerdecay}
If $\psi_j\in H^{-1/2}(\Gamma_j)$ has zero charge, then,
\begin{align*}
 \mathcal S_j[\psi_j](x)=O(|x|^{-1})\qquad \text{ and } \qquad
 \nabla\mathcal S_j[\psi_j](x)=O(|x|^{-2})
\end{align*}
as $|x|\to+\infty$. Its gradient has finite exterior energy.
\end{lemma}
\begin{proof}
See Appendix~\ref{proof:layerdecay}.
\end{proof}

\begin{lemma}[Triviality of the transmission kernel]\label{layerkernel}
For every $\ell>0$, the kernel of $\lambda_\ell I-\mathbb K^*$ on $\mathcal X$ is $\{0\}$.
\end{lemma}
\begin{proof}
See Appendix~\ref{proof:layerkernel}.
\end{proof}

\begin{lemma}[Construction of a physical solution]\label{layerconstruction}
If $\boldsymbol\varphi\in\mathcal X$ solves $(\lambda_\ell I-\mathbb K^*)\boldsymbol\varphi=\mathbf g_F$, then, $F+\mathcal S_i[\varphi_i]+\mathcal S_e[\varphi_e]$ solves~\eqref{pde} in the stated energy and decay class.
\end{lemma}
\begin{proof}
See Appendix~\ref{proof:layerconstruction}.
\end{proof}

\begin{lemma}[Recovery of the layer densities]\label{layerconverse}
Every solution $v$ of~\eqref{transmission} in the prescribed class has densities $\varphi_j=\frac{\partial}{\partial \nu_j}v|_+-\frac{\partial}{\partial \nu_j}v|_-$, with $(\varphi_i,\varphi_e)\in\mathcal X$, and satisfies $v=F+\mathcal S_i[\varphi_i]+\mathcal S_e[\varphi_e]$. These densities solve the system in Lemma~\ref{layerflux}.
\end{lemma}
\begin{proof}
See Appendix~\ref{proof:layerconverse}.
\end{proof}

\begin{theorem}[Boundary-integral representation]\label{BIE}
For every $\ell>0$, the operator $\lambda_\ell I-\mathbb K^*$ is boundedly invertible on $\mathcal X$. The unique solution of \eqref{transmission} is
\begin{equation}
 u_\ell=F+\mathcal S_i[\varphi_i]+\mathcal S_e[\varphi_e] \qquad \text{ and } 
 \qquad
 \binom{\varphi_i}{\varphi_e}
  =(\lambda_\ell I-\mathbb K^*)^{-1}\mathbf g_F.
 \label{BIrepresentation}
\end{equation}
Conversely, every solution of \eqref{transmission} in the stated energy and decay class has this representation, with unique densities in $\mathcal X$.
\end{theorem}
\begin{proof}
Lemmas~\ref{layercharges} and~\ref{layerkernel} give bounded invertibility by the Fredholm alternative. The forcing lies in $\mathcal X$, so its inverse image defines a density pair. Lemma~\ref{layerconstruction} gives the solution in~\eqref{BIrepresentation}. Conversely, Lemma~\ref{layerconverse} recovers densities in $\mathcal X$ from any physical solution. Lemma~\ref{layerflux} and invertibility identify them uniquely, proving the converse and uniqueness.
\end{proof}

\begin{remark}[Exact spectral parameter]
The parameter in the boundary-integral resolvent is $\lambda_\ell$, not the material increment $c\ell^p+i\ell$. The coefficient realization below uses $\zeta_\ell=-2\lambda_\ell$. Approximating this coordinate before locating a narrow resonance can change the location by more than the linewidth. Also, the general integral formulation above does not yet give the power classification for arbitrary shapes or sources.
\end{remark}

\section{Conformal coefficient formulation}\label{conformal}
We return to the geometry and source fixed in Section~\ref{problemgeometry}. The following coefficient representation retains the effect of the noncircular core while making modal decay and the physical energy explicit.

\subsection{Grunsky coefficients and core impedance}
Our next task is to express the response of the core in coordinates compatible with the Fourier expansion in the shell. Although conformal coordinates reduce the shell to an annulus, the map $\Psi$ does not parametrize the core interior, whose geometry must still enter the inner transmission condition. To retain this information, we will use the classical Grunsky coefficients associated with $\Psi$, following their use in series representations of transmission problems in \cite{junglim,choikimlim} and of two-interface operators in \cite{choilimshipman}. These coefficients arise from a logarithmic expansion of the conformal map and describe how the boundary traces of harmonic polynomial fields couple different angular modes. We will assemble their normalized values into a matrix $G$ and establish decay estimates that quantify this coupling at high mode indices.

\begin{lemma}\label{geometry}
The map $\Psi$ is injective on $|w|\ge\rho$ and has nonzero derivative there. Define symmetric coefficients by
\begin{align*}
 \log\frac{\Psi(w)-\Psi(v)}{w-v}
 =-\sum_{n,m\ge1}\beta_{nm}w^{-n}v^{-m}.
\end{align*}
Their values and the energy-normalized matrix are
\begin{equation}
 \beta_{nm}=\begin{cases}
 \displaystyle\frac{\varepsilon^k}{k}\binom{k}{n-k},
 &n+m=3k,\ k\le n\le2k,\\0,&\text{otherwise},
 \end{cases}
 \qquad \text{ and } \qquad G_{nm}=\sqrt{nm}\,\beta_{nm}\rho^{-n-m}.
 \label{grunsky}
\end{equation}
Put $b_*=(2|\varepsilon|/\rho^3)^{1/3}<1/2$. With $\eta$ as in~\eqref{standingcoupling}, one has $\eta=b_*^{3/2}$. Then, $G$ is real, compact, self-adjoint, and
\begin{equation}
 \gamma:=\|G\|\le\frac{3b_*^2}{2(1-b_*)}<1,
 \qquad \|Ge_n\|\le C_G\eta^n
 \qquad \text{ and } \qquad C_G=\frac3{2\sqrt{1-b_*^2}}.
 \label{gdecay}
\end{equation}
\end{lemma}
\begin{proof}
See Appendix~\ref{proof:geometry}.
\end{proof}

To use these geometric coefficients in the transmission conditions, we next introduce compatible coordinates for the boundary potential and its normal derivative. Their different normalizations account for the different Sobolev regularities of Dirichlet and normal traces.

\begin{lemma}[Trace coordinates for the core]\label{coretraces}
Let $\mathbb T=\mathbb R/(2\pi\mathbb Z)$, and define
\begin{equation}
 \gamma_i(\theta)=\Psi(\rho e^{i\theta}) \qquad \text{ and } \qquad
 h_i(\theta)=|\gamma_i'(\theta)|,\qquad ds=h_i(\theta)\,d\theta.
 \label{coremetric}
\end{equation}
For $v\in H^{1/2}(\Gamma_i)$ set $\widehat v=v\circ\gamma_i$. For $q\in H^{-1/2}(\Gamma_i)$ define its metric pullback by
\begin{equation}
 \langle\widetilde q,\psi\rangle_{\mathbb T}
 =\langle q,\psi\circ\gamma_i^{-1}\rangle_{\Gamma_i},
 \qquad \psi\in H^{1/2}(\mathbb T).
 \label{fluxpullback}
\end{equation}
For $s=\pm1/2$, let $H^s_{\mathrm{ev},0}(\mathbb T)$ denote the subspace of distributions $w$ satisfying $w(-\theta)=w(\theta)$ and $\langle w,1\rangle_{\mathbb T}=0$, and set
\begin{equation}
 \mathcal V_i=\{v\in H^{1/2}(\Gamma_i):
                  \widehat v\in H^{1/2}_{\mathrm{ev},0}(\mathbb T)\} \qquad \text{ and } \qquad
 \mathcal Q_i=\{q\in H^{-1/2}(\Gamma_i):
                  \widetilde q\in H^{-1/2}_{\mathrm{ev},0}(\mathbb T)\}.
 \label{coretracespaces}
\end{equation}
The maps defined by
\begin{equation}
 \begin{aligned}
 (\mathscr D_i v)_n
 =\frac{\sqrt n}{\pi}\int_0^{2\pi}\widehat v(\theta)\cos(n\theta)\,d\theta \qquad \text{ and } \qquad
 (\mathscr N_i q)_n
 =\frac1{\pi\sqrt n}\langle\widetilde q,\cos(n\theta)\rangle_{\mathbb T},
 \qquad n\ge1,
 \end{aligned}
 \label{corenormalization}
\end{equation}
are bounded isomorphisms from $\mathcal V_i$ and $\mathcal Q_i$, respectively, onto $\ell^2(\mathbb N)$. For $\mathbf v=\mathscr D_i v$ and $\mathbf k=\mathscr N_i q$, their inverse series are
\begin{equation}
 \widehat v=\sum_{n\ge1}\frac{\mathbf v_n}{\sqrt n}\cos(n\theta) \qquad \text{ and } 
 \qquad
 \widetilde q=\sum_{n\ge1}\sqrt n\,\mathbf k_n\cos(n\theta),
 \label{coreinverses}
\end{equation}
with convergence in $H^{1/2}(\mathbb T)$ and $H^{-1/2}(\mathbb T)$, respectively, and
\begin{align*}
 \|v\|_{H^{1/2}(\Gamma_i)}^2\asymp\|\mathbf v\|_{\ell^2}^2 \qquad \text{ and } \qquad
 \|q\|_{H^{-1/2}(\Gamma_i)}^2\asymp\|\mathbf k\|_{\ell^2}^2.
\end{align*}
Moreover, the core Dirichlet-to-Neumann map
\begin{equation}
 \Lambda_i v=\frac{\partial}{\partial\nu_i}U_v|_{\Gamma_i},
 \qquad
 \Delta U_v=0\ \text{in }\Omega_i \qquad \text{ and } \qquad U_v|_{\Gamma_i}=v,
 \label{coreDN}
\end{equation}
is bounded from $\mathcal V_i$ to $\mathcal Q_i$.
\end{lemma}
\begin{proof}
See Appendix~\ref{proof:coretraces}.
\end{proof}

\begin{remark}
	\begin{enumerate}
		\item For smooth $q$ this definition reads $\widetilde q(\theta)=h_i(\theta)q(\gamma_i(\theta))$. Thus, it includes the arclength factor.
		\item Pairings in these coordinate definitions are bilinear, without complex conjugation.
	\end{enumerate}
\end{remark}

\paragraph*{Faber polynomials associated with the core.}
We now need harmonic fields in the physical core whose boundary potential and normal derivative can both be expressed in the coordinates of Lemma~\ref{coretraces}. For this purpose we use harmonic polynomials adapted to the conformal geometry, obtained from the Faber polynomials.
We write $\mathcal F_n$ for the $n$th Faber polynomial associated with the normalized exterior map $\Psi$, not for a Fourier coefficient of the incident potential $F$. These are the monic polynomials of degree $n$ defined, for fixed $z$ and sufficiently large $|w|$, by
\begin{equation}
 \log\frac{\Psi(w)-z}{w}
 =-\sum_{n\ge1}\frac{\mathcal F_n(z)}{n}w^{-n}.
 \label{faberdefinition}
\end{equation}
The branch of the logarithm tends to zero at infinity. This is the classical Faber construction \cite{faber}, and we refer to \cite{junglim,choikimlim} for its use in transmission problems. For our map,
\begin{align*}
 \mathcal F_1(z)=z,\qquad \mathcal F_2(z)=z^2 \qquad \text{ and } \qquad
 \mathcal F_3(z)=z^3-3\varepsilon.
\end{align*}
Substitution into the defining Grunsky expansion gives
\begin{align*}
 \mathcal F_n(\Psi(v))=v^n+\sum_{m\ge1}n\beta_{nm}v^{-m}.
\end{align*}
Thus, the real parts of these polynomials are harmonic fields in the physical core whose traces contain one prescribed growing mode together with the geometry-induced decaying modes. They allow us to compute the response of the actual core, rather than replacing it by a disk.

The preceding construction allows us to relate the two traces through the geometry. The following core impedance expresses the normal derivative supplied by the core in terms of its boundary potential, so that the inner transmission condition can be imposed without replacing the core by a disk.

\begin{lemma}[Core impedance]\label{core}
With the spaces and coordinate maps defined above, the coefficient representation of the core Dirichlet-to-Neumann map is
\begin{equation}
 L:=\mathscr N_i\Lambda_i\mathscr D_i^{-1}
   =(I-G)(I+G)^{-1}.
 \label{impedance}
\end{equation}
It is bounded, positive, invertible, and self-adjoint on $\ell^2$.
\end{lemma}
\begin{proof}
See Appendix~\ref{proof:core}.
\end{proof}

\subsection{The fixed source in conformal coordinates}
The source must remain fixed as the loss decreases, otherwise enhancement could be introduced through the forcing itself. Radial smoothing is the property that gives a regular source while preserving the dipole field near the shell exactly. Its conformal position determines the exponential decay of the incident coefficients, and hence, how efficiently it excites the high modes.

\begin{lemma}[Exterior field of the smooth dipole source]\label{smoothdipole}
For the fixed source of Section~\ref{problemgeometry}, suppose $\operatorname{supp}\chi\subset\overline{B_b(0)}$. The free potential in~\eqref{freepotential} satisfies
\begin{equation}
 F(z)=\frac1{2\pi}\operatorname{Re}\frac1{z-z_0},
 \qquad |z-z_0|>b.
 \label{dipole}
\end{equation}
Moreover, $F\in C^\infty(\mathbb R^2)$, $\nabla F\in L^2(\mathbb R^2)$, and
$F(z)=O(|z|^{-1})$, $\nabla F(z)=O(|z|^{-2})$ as $|z|\to+\infty$.
\end{lemma}
\begin{proof}
See Appendix~\ref{proof:smoothdipole}.
\end{proof}

Near the closed shell write
\begin{equation}
 F(\Psi(re^{i\theta}))=F_0+
 \sum_{n\ge1}\frac{A_n r^n+(D_{\rm inc})_n r^{-n}}{\sqrt n}\cos n\theta.
 \label{incident}
\end{equation}
\begin{lemma}\label{source}
The growing coefficients and their weights are
\begin{equation}
 A_n=-\frac{\sqrt n}{2\pi\Psi'(R)R^{n+1}},\qquad
 |A_n|^2=\kappa n t_n^\alpha,
 \quad t_n=\rho^n\qquad \text{and} \qquad \kappa=\frac1{4\pi^2\Psi'(R)^2R^2}>0.
 \label{A}
\end{equation}
Both $A$ and $D_{\rm inc}$ belong to $\ell^2$. With $T=\operatorname{diag}(t_n)$, they satisfy
\begin{equation}
 D_{\rm inc}=TGT A.
 \label{compatible}
\end{equation}
\end{lemma}
\begin{proof}
See Appendix~\ref{proof:source}.
\end{proof}

\subsection{Exact coefficient realization of the integral problem}
We can now combine the core response with the fixed-source expansion to express the solution through its shell coefficients. We will introduce an operator $J$ describing the geometric coupling and a vector $g$ describing the forcing, while retaining the material parameter from~\eqref{lambdazeta}. This separates the three ingredients of the response and gives an exact representation before any numerical truncation.

Write $\zeta=\zeta_\ell$ for the parameter in~\eqref{lambdazeta}, and set $x=\operatorname{Re}\zeta$, $y=\operatorname{Im}\zeta$. That relation yields
\begin{equation}
 s_\ell=\frac{\zeta-1}{\zeta+1},\qquad
 y=\frac{2\ell}{|1-s_\ell|^2}>0 \qquad \text{ and } \qquad
 \ell|1+\zeta|^2=2y.
 \label{coordinate}
\end{equation}
Set $W=(I-T^2)^{1/2}$ and $C=(I-G^2)^{1/2}$, both boundedly invertible. Define on $\HH$
\begin{equation}
 J=\begin{pmatrix}-TGT&-TC\\-CT&G\end{pmatrix},\quad
 M=\begin{pmatrix}W&0\\WGT&WC\end{pmatrix},\quad
 g=\binom A0,\qquad V=(\zeta-J)^{-1}g.
 \label{JM}
\end{equation}
Let $P_1:\HH\to\ell^2(\mathbb N)$ denote the first-component projection, $P_1(v_1,v_2)=v_1$.

The shell coordinates used below are
\begin{equation}
 u_\ell(\Psi(re^{i\theta}))=F_0+
 \sum_{n\ge1}\frac{b_nr^n+t_nq_nr^{-n}}{\sqrt n}\cos(n\theta).
 \label{shellexpansion}
\end{equation}
Here $b$ and $q$ record the growing and decaying parts of the shell field. The following lemmas separate the interface conditions, the operator metric, and the reconstruction of the physical solution.

\begin{lemma}[Interface equations in coefficient coordinates]\label{coeffinterfaces}
For square-summable shell vectors $b,q$ with constant mode $F_0$, transmission is equivalent to
\begin{align*}
 \zeta b+Tq=(1+\zeta)A \qquad \text{ and } \qquad (I-\zeta G)Tb+(\zeta I-G)q=0.
\end{align*}
Writing $q=GTb+Cz$, these equations are equivalent to $(\zeta-J)(b,z)=(1+\zeta)g$.
\end{lemma}
\begin{proof}
See Appendix~\ref{proof:coeffinterfaces}.
\end{proof}

\begin{lemma}[Self-adjoint realization and energy metric]\label{coeffmetric}
The operator $J$ in~\eqref{JM} is compact, self-adjoint and injective. Moreover $M$ is boundedly invertible,
\begin{align*}
 M^*M=I-J^2,\qquad \|J\|<1 \qquad \text{ and } \qquad
 \|(\zeta-J)^{-1}\|\le(\operatorname{Im}\zeta)^{-1}.
\end{align*}
\end{lemma}
\begin{proof}
See Appendix~\ref{proof:coeffmetric}.
\end{proof}

\begin{lemma}[Physical energy and exterior observation]\label{coeffenergy}
Let $b,z,q$ satisfy the equations in Lemma~\ref{coeffinterfaces} and set $V=(1+\zeta)^{-1}(b,z)$. The shell series and its continuous exterior extension have
\begin{align*}
 E=\pi\ell(\|Wb\|^2+\|Wq\|^2)
   =2\pi y\langle V,(I-J^2)V\rangle.
\end{align*}
The decaying scattered coefficients are $h=(1-\zeta^2)P_1V+\zeta A-D_{\rm inc}$ and satisfy $S_r=\pi\sum_{n\ge1}r^{-2n}|h_n|^2$.
\end{lemma}
\begin{proof}
See Appendix~\ref{proof:coeffenergy}.
\end{proof}

\begin{lemma}[Reconstruction of the energy solution]\label{coeffreconstruction}
The vectors $(b,z)=(1+\zeta)(\zeta-J)^{-1}g$ and $q=GTb+Cz$ reconstruct an $H^1_{\rm loc}$ solution of~\eqref{pde} with finite gradient energy and the prescribed decay. Its shell constant is $F_0$, and its logarithmic and sine modes vanish.
\end{lemma}
\begin{proof}
See Appendix~\ref{proof:coeffreconstruction}.
\end{proof}

Taken together, these lemmas show that the coefficient equations retain both the transmission conditions and the physical energy. We can, therefore, express the unique boundary-integral solution in these coordinates and read off its shell power and exterior response from the same coefficient vector.

\begin{theorem}[Exact forced realization]\label{physical}
The operator $J$ is compact, self-adjoint, injective, and $\|J\|<1$. The solution of Theorem~\ref{BIE} has the shell expansion~\eqref{shellexpansion}, with coefficients
\begin{equation}
 \binom bz=(1+\zeta)V \qquad \text{ and } \qquad q=GTb+Cz.
 \label{reconstruct}
\end{equation}
The exact energy and observation identities are
\begin{align}
 &M^*M=I-J^2,\label{metric}\\
 &E=\pi\ell(\|Wb\|^2+\|Wq\|^2) =2\pi y\langle V,(I-J^2)V\rangle,\label{power}\\
 &h=(1-\zeta^2)P_1V+\zeta A-D_{\rm inc},\qquad
 S_r=\pi\|D_r h\|^2,\quad D_r=\operatorname{diag}(r^{-n}),\label{observe}
\end{align}
In particular,
 \begin{align}
 2\pi(1-\|J\|^2)y\|V\|^2&\le E\le2\pi y\|V\|^2.\label{coercive}
\end{align}
\end{theorem}
\begin{proof}
Lemma~\ref{coeffmetric} establishes the stated operator properties and gives the unique vector $V$. Lemma~\ref{coeffinterfaces} identifies its coefficients with the transmission equations. Lemma~\ref{coeffreconstruction} identifies the reconstructed field with the unique solution of Theorem~\ref{BIE}. Lemma~\ref{coeffenergy} proves the energy and observation identities. Finally $0\le\|JV\|^2\le\|J\|^2\|V\|^2$ gives~\eqref{coercive}.
\end{proof}

\subsection{From layer densities to coefficient vectors}
Let us now make explicit the relation between the boundary-integral solution and the coefficient vectors. Alongside the inner metric in~\eqref{coremetric}, define the outer metric by
\begin{equation}
 h_e(\theta)=|e^{i\theta}\Psi'(e^{i\theta})| \qquad \text{ and } 
 \qquad ds_{\Gamma_e}=h_e(\theta)\,d\theta.
 \label{outermetric}
\end{equation}
Here $h_j$ is the boundary metric, not the scattered coefficient vector $h$ in \eqref{observe}. The radial derivative in the conformal plane satisfies
$\frac{\partial}{\partial \nu_j}=h_j^{-1}r\partial_r$ on either trace of $\Gamma_j$.

\begin{lemma}[Recovery of zero-charge densities]\label{densityrecovery}
For the coefficient solution in Theorem~\ref{physical}, the normal derivative jumps are the $H^{-1/2}$ series~\eqref{densityinner}--\eqref{densityouter}. Both have zero charge.
\end{lemma}
\begin{proof}
See Appendix~\ref{proof:densityrecovery}.
\end{proof}

\begin{lemma}[Invertibility of the trace coordinates]\label{tracecoordinates}
For the unique physical solution, let $v_i,v_e$ be the normalized nonconstant cosine trace vectors. Then, $b,q$ are given by~\eqref{traceinverse}, with
\begin{align*}
 \|b\|\le\frac{\|v_e\|+\rho\|v_i\|}{1-\rho^2} \qquad \text{ and } \qquad
 \|q\|\le\frac{\|v_i\|+\rho\|v_e\|}{1-\rho^2}.
\end{align*}
The recovered vector $V=(1+\zeta)^{-1}(b,C^{-1}(q-GTb))$ satisfies $(\zeta-J)V=g$.
\end{lemma}
\begin{proof}
See Appendix~\ref{proof:tracecoordinates}.
\end{proof}

These lemmas provide the maps in both directions: the shell coefficients determine the interface densities, and the two boundary traces recover the shell coefficients. The following theorem states the resulting equivalence with the boundary-integral formulation.

\begin{theorem}[Equivalence of the two realizations]\label{bridge}
For the fixed even source above, the densities in \eqref{BIrepresentation} are recovered from the shell coefficients in \eqref{reconstruct} by
\begin{align}
 \varphi_i(\Psi(\rho e^{i\theta}))
 &=\frac{1-s_\ell}{h_i(\theta)}
     \sum_{n\ge1}\sqrt n\,(t_nb_n-q_n)\cos n\theta,\label{densityinner}
\end{align}
and
\begin{align}
 \varphi_e(\Psi(e^{i\theta})) &=\frac{s_\ell-1}{h_e(\theta)}
     \sum_{n\ge1}\sqrt n\,(b_n-t_nq_n)\cos n\theta.\label{densityouter}
\end{align}
The series are interpreted in $H^{-1/2}$ and each density has zero charge.
Conversely, the layer solution determines $b,q$ by its two Dirichlet traces. If $v_i,v_e$ are the $\sqrt n$-normalized cosine coefficients of those traces after removing $F_0$, then,
\begin{equation}
 b=(I-T^2)^{-1}(v_e-Tv_i) \qquad \text{ and } \qquad
 q=(I-T^2)^{-1}(v_i-Tv_e).
 \label{traceinverse}
\end{equation}
The vector in the self-adjoint realization is
\begin{align*}
 z=C^{-1}(q-GTb) \qquad \text{ and } \qquad
 V=(1+\zeta)^{-1}(b,z).
\end{align*}
In particular \eqref{BIrepresentation} and \eqref{reconstruct} describe exactly the same solution for every positive loss.
\end{theorem}
\begin{proof}
Lemma~\ref{densityrecovery} gives the displayed densities and their zero charges. Lemma~\ref{tracecoordinates} gives the bounded inverse trace formulas and the coefficient equation. Theorem~\ref{BIE} identifies the recovered density pair uniquely, so both descriptions represent the same physical field.
\end{proof}

\section{Uniform comparison of forced resolvents}\label{sec:uniform}
The coefficient formulation gives an exact description of the physical field, but its coupled modes do not yet reveal the behavior as loss vanishes. We will compare this response with an auxiliary circular problem in which the modes can be computed separately, keeping the forcing vector fixed. The comparison must remain valid near increasingly narrow resonances, where a small change in the operator need not produce a small change in its inverse. To obtain this control, we will factor the geometric correction relative to the circular modal scales and use the decay established in~\eqref{gdecay}. The resulting power comparison will transfer the scalar estimates of Section~\ref{sec:power} to the noncircular shell. A more precise comparison along sequences with growing response will also be used in the sharpness and endpoint arguments of Section~\ref{sec:localisation}.
Define the circular comparator, with exactly the same forcing vector, by
\begin{equation}
 J_0=\begin{pmatrix}0&-T\\-T&0\end{pmatrix},\qquad
 \TT=|J_0|=\operatorname{diag}(T,T),\qquad
 V_0=(\zeta-J_0)^{-1}g,\qquad
 E_0=2\pi y\langle V_0,(I-J_0^2)V_0\rangle.
 \label{comparator}
\end{equation}
The comparator is an auxiliary circular problem, not the original physical source moved to a circle. Its incident amplitude includes $\Psi'(R)^{-1}$.

\begin{lemma}[Relative compact factorization]\label{relative}
If $\eta<\rho$, there is compact $K$ with $J=J_0(I+K)$ and boundedly invertible $I+K$. If $s>0$ and $\eta<\rho^{1+s}$, then, $K=\TT^s L_s$ with $L_s$ Hilbert--Schmidt.
\end{lemma}
\begin{proof}
See Appendix~\ref{proof:relative}.
\end{proof}

\begin{lemma}\label{multiplier}
For $a=|\zeta|$ sufficiently small, $y>0$, and fixed $s>0$,
\begin{align*}
 \|\zeta(\zeta-J_0)^{-1}\TT^s\|
 \le C_s\left(a^{\min(s,1)}+\frac{a^{1+s}}y\right).
\end{align*}
\end{lemma}
\begin{proof}
See Appendix~\ref{proof:multiplier}.
\end{proof}

\begin{lemma}[Uniform invertibility of the comparison operator]\label{uniformresolvent}
Under the standing hypotheses, put
\begin{align*}
 F_\zeta=I+K-\zeta(\zeta-J_0)^{-1}K.
\end{align*}
Then, $F_\zeta V=V_0$ and, along the material path~\eqref{crossings},
\begin{align*}
 \lim_{\ell\downarrow0}\|F_{\zeta(\ell)}-(I+K)\|&=0 \qquad \text{ and } \qquad
 \lim_{\ell\downarrow0}\|F_{\zeta(\ell)}^{-1}-(I+K)^{-1}\|&=0.
\end{align*} Consequently $\|V\|\asymp\|V_0\|$ and $E\asymp E_0$, with comparison constants independent of $g\in\HH$.
\end{lemma}
\begin{proof}
See Appendix~\ref{proof:uniformresolvent}.
\end{proof}

\begin{lemma}[Leading response for a fixed forcing]\label{leadingresolvent}
Under the standing hypotheses, fix $g\in\HH$ and a positive sequence $(\ell_j)$ satisfying $\lim_{j\to+\infty}\ell_j=0$ and $\lim_{j\to+\infty}\|V_0(\ell_j)\|=+\infty$. Here $V(\ell)$ and $V_0(\ell)$ denote the responses evaluated at $\zeta=\zeta(\ell)$. Then,
\begin{align*}
 \lim_{j\to+\infty}\frac{\|V(\ell_j)-V_0(\ell_j)\|}{\|V_0(\ell_j)\|}=0 \qquad \text{ and }
 \qquad
 \lim_{j\to+\infty}\frac{E(\ell_j)}{E_0(\ell_j)}=1.
\end{align*}
\end{lemma}
\begin{proof}
See Appendix~\ref{proof:leadingresolvent}.
\end{proof}

\begin{proposition}[Uniform and leading comparisons]\label{comparison}
Under the standing hypotheses, $\|V\|\asymp\|V_0\|$ and $E\asymp E_0$ uniformly for sufficiently small loss. The operator norm comparison is uniform over all $g\in\HH$, although only physically compatible $g$ is used for~\eqref{pde}. For fixed $g$ and every positive sequence $(\ell_j)$ with $\lim_{j\to+\infty}\ell_j=0$ and $\lim_{j\to+\infty}\|V_0(\ell_j)\|=+\infty$,
\begin{equation}
 \lim_{j\to+\infty}\frac{\|V(\ell_j)-V_0(\ell_j)\|}{\|V_0(\ell_j)\|}=0 \qquad \text{ and }
 \qquad
 \lim_{j\to+\infty}\frac{E(\ell_j)}{E_0(\ell_j)}=1.
 \label{leading}
\end{equation}
\end{proposition}
\begin{proof}
Lemma~\ref{uniformresolvent} proves both uniform comparisons. Lemma~\ref{leadingresolvent} gives~\eqref{leading} whenever the circular response norm diverges. The fixed-forcing condition in the latter lemma is precisely the condition used in the leading assertion.
\end{proof}

\section{Positive sums and proof of the power classification}\label{sec:power}
We now use Proposition~\ref{comparison} to determine how the shell power behaves as the loss tends to zero. It is enough to estimate the circular energy in~\eqref{comparator}, whose modes can be treated separately. We will first express this energy as a sum of nonnegative modal contributions. We can distinguish the response near an individual resonance from the combined contribution of the remaining modes, without having to control cancellations between them. This distinction will explain why some detuning paths produce divergence for every sufficiently small loss, whereas others produce large peaks separated by vanishing responses.

\begin{lemma}[Circular modal response and positive energy sum]\label{circularsum}
For the source coefficients in~\eqref{A} and the circular comparator in~\eqref{comparator}, the response at each index $n\ge1$ is
\begin{equation}
 (V_0)_n=\frac{A_n}{\zeta^2-t_n^2}\binom\zeta{-t_n} \qquad \text{ and } \qquad
 \|(V_0)_n\|^2=|A_n|^2\frac{a^2+t_n^2}{|\zeta^2-t_n^2|^2}.
 \label{V0}
\end{equation}
The corresponding energy is
\begin{equation}
 E_0=\pi y\sum_{n\ge1}|A_n|^2(1-t_n^2)
 \left(\frac1{|\zeta-t_n|^2}+\frac1{|\zeta+t_n|^2}\right).
 \label{positivesum}
\end{equation}
\end{lemma}
\begin{proof}
See Appendix~\ref{proof:circularsum}.
\end{proof}

The two denominators in~\eqref{positivesum} measure the distance of the material coordinate from the circular poles $t_n$ and $-t_n$. The imaginary part $y$ keeps each contribution finite at positive loss, while the source weights in~\eqref{A} determine how strongly the increasingly fine modes are excited. We will use $\mathscr{D}$ to measure the distance to the nearest positive pole and $H_\alpha$ to describe the size of the remaining response.
For small $q>0$, define
\begin{equation}
 \mathscr{D}(q)=\min_{n\ge1}|q-t_n| \qquad \text{ and } \qquad
 H_\alpha(q)=\begin{cases}
 q^{\alpha-2}\log(1/q),&0<\alpha<2,\\
 \log(1/q)^2,&\alpha=2,\\1,&\alpha>2.
 \end{cases}
 \label{scalarscales}
\end{equation}
We first estimate a sum in which the denominators stay away from resonance. Splitting the modes according to whether $t_n$ is larger or smaller than $q$ gives the scale $H_\alpha(q)$. The change at $\alpha=2$ reflects the balance between the decay of the source coefficients and the amplification of the small modal scales.

\begin{lemma}[Geometric head and tail sums]\label{geometricsums}
For $t_n=\rho^n$ and every fixed $\alpha>0$,
\begin{align*}
 \sum_{n\ge1}\frac{nt_n^\alpha}{(q+t_n)^2}\asymp H_\alpha(q)
 \qquad \text{ as }q\downarrow0.
\end{align*}
\end{lemma}
\begin{proof}
See Appendix~\ref{proof:geometricsums}.
\end{proof}

When the real part of the material coordinate is no larger than a fixed multiple of its imaginary part, no individual pole is approached on a scale smaller than $y$. The preceding sum, therefore, controls the whole response, as the next lemma makes precise.

\begin{lemma}[Sectorial response]\label{sectorialsums}
For each fixed $C>0$, the energy~\eqref{positivesum} satisfies $E_0\asymp yH_\alpha(y)$ uniformly when $|x|\le Cy$ and $y\downarrow0$.
\end{lemma}
\begin{proof}
See Appendix~\ref{proof:sectorialsums}.
\end{proof}

The situation changes when $y$ is much smaller than $|x|$. The material coordinate can pass close to an individual real pole, producing a narrow peak. We isolate this contribution and estimate the remaining modes together, which separates the peak height from the response between resonances.

\begin{lemma}[Nearest-pole decomposition]\label{tangentialsums}
Let $q=|x|>0$ tend to zero, with $y=y(q)>0$ satisfying $\lim_{q\downarrow0}y(q)/q=0$. With the quantities in~\eqref{A} and~\eqref{scalarscales}, choose any index
\begin{align*}
 m\in\operatorname*{arg\,min}_{n\ge1}|q-t_n| \qquad \text{ and }
 \qquad |q-t_m|=\mathscr{D}(q).
\end{align*}
Define the nearest-pole contribution and the remainder by
\begin{align}
 E_0^{\mathrm{near}}
 &:=\frac{\pi y|A_m|^2(1-t_m^2)}
              {(q-t_m)^2+y^2},\label{nearestpoleenergy}\\
 E_0^{\mathrm{rem}}
 &:=\pi y\sum_{\substack{n\ge1\\n\ne m}}
       \frac{|A_n|^2(1-t_n^2)}{(q-t_n)^2+y^2}
   +\pi y\sum_{n\ge1}
       \frac{|A_n|^2(1-t_n^2)}{(q+t_n)^2+y^2}.
       \label{remainderenergy}
\end{align}
Then,
\begin{align}
 E_0&=E_0^{\mathrm{near}}+E_0^{\mathrm{rem}},\label{nearestpoledecomposition}\\
 E_0^{\mathrm{near}}
 \asymp\frac{yq^\alpha\log(1/q)}{\mathscr{D}(q)^2+y^2}&
 \qquad \text{and} \qquad
 E_0^{\mathrm{rem}}\asymp yH_\alpha(q).
 \label{nearestpolebounds}
\end{align}
The comparison constants are independent of $q$, $y$ and the choice of $m$ in the event of a tie, for all sufficiently small $q$ with $0<y\le q$. In particular,~\eqref{envelopes} follows.
\end{lemma}
\begin{proof}
See Appendix~\ref{proof:tangentialsums}.
\end{proof}

\begin{remark}
    We will also need a lower bound that does not distinguish these two approaches to the origin. It follows by retaining modes whose scales are larger than $|\zeta|$, and will be used in the exterior estimates of Section~\ref{sec:localisation}.
\end{remark}

\begin{lemma}[An angle-independent lower bound]\label{anglelower}
For $0<\alpha<2$ and all sufficiently small $\zeta$ with $\operatorname{Im}\zeta=y>0$, one has $E_0\gtrsim ya^{\alpha-2}\log(1/a)$, where $a=|\zeta|$.
\end{lemma}
\begin{proof}
See Appendix~\ref{proof:anglelower}.
\end{proof}

\begin{lemma}[Scalar comparison]\label{sums}
With the quantities in~\eqref{scalarscales}, if $|x|\le C y$, then, $E_0\asymp y H_\alpha(y)$. If $q=|x|>0$ tends to zero with $y=y(q)>0$ and $\lim_{q\downarrow0}y(q)/q=0$, then,
\begin{equation}
 E_0\asymp yH_\alpha(q)+
 \frac{yq^\alpha\log(1/q)}{\mathscr{D}(q)^2+y^2}.
 \label{envelopes}
\end{equation}
All constants are independent of loss within the stated region. Also, for $0<\alpha<2$ and any small $\zeta$ in the upper half-plane,
\begin{equation}
 E_0\gtrsim y a^{\alpha-2}\log(1/a).
 \label{lowera}
\end{equation}
\end{lemma}
\begin{proof}
Lemma~\ref{sectorialsums} gives the sectorial estimate. Lemma~\ref{tangentialsums} gives the tangential decomposition, and Lemma~\ref{anglelower} gives~\eqref{lowera}.
\end{proof}

We can now return to the prescribed permittivity path. Proposition~\ref{comparison} transfers these scalar estimates to the physical shell power. For $p\ge1$, the path remains in the sectorial regime, so one estimate describes all sufficiently small losses. For $0<p<1$, the path approaches the real axis tangentially, and we must distinguish losses at which it crosses a pole from those at which it lies between neighboring poles. We retain the exact material coordinate from~\eqref{crossings}: approximating a crossing location can displace it by more than the width of a narrow resonance.

\begin{lemma}[Power on sectorial paths]\label{powersector}
Under the standing hypotheses and $p\ge1$,
\begin{align*}
 E\asymp
 \begin{cases}
 \ell^{\alpha-1}\log(1/\ell),&0<\alpha<2,\\
 \ell\log(1/\ell)^2,&\alpha=2,\\
 \ell,&\alpha>2.
 \end{cases}
\end{align*}
\end{lemma}
\begin{proof}
See Appendix~\ref{proof:powersector}.
\end{proof}

\begin{lemma}[Exact crossing and gap sequences]\label{powercrossings}
Under the standing hypotheses with $0<p<1$, $q(\ell)=|\operatorname{Re}\zeta(\ell)|$ is strictly increasing for sufficiently small $\ell>0$. For all sufficiently large $n$, there are unique small losses $\ell_n,\widehat\ell_n$ satisfying
\begin{align*}
 q(\ell_n)=t_n,\qquad q(\widehat\ell_n)=\frac{2\rho}{1+\rho}t_n \qquad \text{ and }
 \qquad \ell_n\asymp\widehat\ell_n\asymp t_n^{1/p}.
\end{align*}
\end{lemma}
\begin{proof}
See Appendix~\ref{proof:powercrossings}.
\end{proof}

\begin{lemma}[Power envelopes on tangential paths]\label{powerenvelopes}
Under the standing hypotheses with $0<p<1$, define
\begin{align*}
 B_\ell=\ell H_\alpha(\ell^p) \qquad \text{ and } \qquad
 P_\ell=\ell^{p\alpha-1}\log(1/\ell).
\end{align*}
For all sufficiently small $\ell$, $B_\ell\lesssim E\lesssim B_\ell+P_\ell$. On the sequences of Lemma~\ref{powercrossings},
\begin{align*}
 E(\ell_n)\asymp B_{\ell_n}+P_{\ell_n} \qquad \text{ and } \qquad
 E(\widehat\ell_n)\asymp B_{\widehat\ell_n}.
\end{align*}
\end{lemma}
\begin{proof}
See Appendix~\ref{proof:powerenvelopes}.
\end{proof}

These estimates identify the two mechanisms in the power classification. Divergence of the background forces the power to grow along the entire loss path. When the background vanishes but the peaks grow, the crossing and gap sequences instead give intermittent resonance. If both contributions vanish, the shell power tends to zero. We now use these alternatives to prove Theorem~\ref{mainpower}.

\begin{proof}[Proof of Theorem~\ref{mainpower}]
For $p\ge1$, Lemma~\ref{powersector} gives divergence for $\alpha\le1$ and decay for $\alpha>1$, including logarithmic divergence at $\alpha=1$.

For $0<p<1$, use the envelopes in Lemma~\ref{powerenvelopes}. Their background term satisfies
\begin{align*}
 B_\ell\asymp
 \begin{cases}
 \ell^{1+p(\alpha-2)}\log(1/\ell),&\alpha<2,\\
 \ell\log(1/\ell)^2,&\alpha=2,\\
 \ell,&\alpha>2.
 \end{cases}
\end{align*}
Thus, 
\begin{align*}
    \left\{
    \begin{aligned}
        & \lim_{\ell\downarrow0}B_\ell=+\infty, \qquad && \text{ for } \quad 0<\alpha\le2-1/p,\\
        & \lim_{\ell\downarrow0}B_\ell=0, && \text{ for } \alpha>\max(0,2-1/p).
    \end{aligned}
    \right.
\end{align*}
The peak term $P_\ell$ diverges for $\alpha\le1/p$ and vanishes for $\alpha>1/p$. Since $2-1/p<1/p$, the uniform lower bound gives full-path divergence in the first range. In the intermediate range, the exact crossing and gap sequences of Lemma~\ref{powercrossings} give divergent peaks and vanishing gaps. Above $1/p$, the uniform upper bound gives $\lim_{\ell\downarrow0}E(\ell)=0$. At each positive threshold the corresponding power of $\ell$ is zero and the logarithm diverges. The full-divergence range is empty when $p\le1/2$.
\end{proof}

\section{Weighted comparison and its kernel}\label{sec:weighted}
The power estimates describe how strongly the shell absorbs, but leave open how much of the resonant field reaches an exterior observer. This distinction matters because the fine angular modes excited at small loss also decay most rapidly with distance. To account for this attenuation, we will compare the physical and circular responses after applying the modal weights associated with~\eqref{comparator}. These weights retain the spatial information that is absent from a comparison of the total shell powers alone.

In a circular shell, each angular mode has its own decay rate. The boundary deformation couples these modes, so we must also account for how it redistributes the response among them. The next estimate expresses the decay of this coupling in both modal indices, which is what we mean here by two-sided smoothing. Write $P=J-J_0$ and $D=I-C$.

\begin{lemma}[Two-sided smoothing]\label{twosided}
For $\sigma\ge1/2$ and $\eta<\rho^{2\sigma}$, bounded compact extensions
$G_\sigma=T^{-\sigma}GT^{-\sigma}$ and $D_\sigma=T^{-\sigma}DT^{-\sigma}$ exist, and
\begin{equation}
 P=\TT^\sigma P_\sigma\TT^\sigma \qquad \text{ and } \qquad
 P_\sigma=\begin{pmatrix}-TG_\sigma T&TD_\sigma\\D_\sigma T&G_\sigma\end{pmatrix}
 \label{weightedP}
\end{equation}
is compact and self-adjoint.
\end{lemma}
\begin{proof}
See Appendix~\ref{proof:twosided}.
\end{proof}

This decay allows us to examine the geometric correction on the scales relevant to exterior observation. As the loss decreases, resonant amplification becomes stronger, while the observation weights increasingly suppress the modes involved. The condition $a^{2\sigma}/y=O(1)$ expresses the balance under which the weighted comparison has a controlled limit.

\begin{lemma}[Norm limit after compact multiplication]\label{weightedlimit}
Assume $\sigma\ge1/2$, $\eta<\rho^{2\sigma}$, $\zeta\to0$ in the upper half-plane, and $a^{2\sigma}/y=O(1)$. Define
\begin{align*}
 Z_{\sigma,\zeta}=\TT^\sigma(\zeta-J_0)^{-1}\TT^\sigma,\quad
 Z_{\sigma,0}=\begin{pmatrix}0&T^{2\sigma-1}\\T^{2\sigma-1}&0\end{pmatrix}
 \qquad \text{ and } \qquad F_{\sigma,\zeta}=I-Z_{\sigma,\zeta}P_\sigma.
\end{align*}
Set $F_{\sigma,0}:=I-Z_{\sigma,0}P_\sigma$. Along the approaches to zero specified above, for $v\in\HH$, we have
\begin{align*}
 \lim_{\zeta\to0}\|(Z_{\sigma,\zeta}-Z_{\sigma,0})v\|&=0
 \qquad \text{ and } \qquad
 \lim_{\zeta\to0}\|F_{\sigma,\zeta}-F_{\sigma,0}\|=0.
\end{align*}
\end{lemma}
\begin{proof}
See Appendix~\ref{proof:weightedlimit}.
\end{proof}

There is, however, a distinction between attenuating a field and discarding part of the solution space. Since the weights become small at high indices, we must check that the limiting description introduces no additional homogeneous response. The next lemma provides this consistency with the original, unweighted problem.

\begin{lemma}[Recovery of the unweighted kernel]\label{weightedkernel}
Under the hypotheses of Lemma~\ref{weightedlimit}, every $X\in\ker F_{\sigma,0}$ belongs to $\operatorname{Ran}\TT^\sigma$. Writing $X=\TT^\sigma u$, we have $Ju=0$, and hence, $X=0$.
\end{lemma}
\begin{proof}
See Appendix~\ref{proof:weightedkernel}.
\end{proof}

With this consistency established, the weighted circular response remains a reliable measure of the physical response as loss vanishes. The geometry modifies the relation between them, but the transformation and its inverse stay bounded. Thus, the comparison preserves the modal attenuation needed to study what remains visible outside the shell.

\begin{lemma}[Weighted inverse and response identity]\label{weightedinverse}
Under the hypotheses of Lemma~\ref{weightedlimit}, $F_{\sigma,\zeta}$ and $F_{\sigma,0}$ are boundedly invertible for sufficiently small loss,
\begin{align*}
 \TT^\sigma V=F_{\sigma,\zeta}^{-1}\TT^\sigma V_0 \qquad \text{ and } \qquad
 \lim_{\zeta\to0}\|F_{\sigma,\zeta}^{-1}-F_{\sigma,0}^{-1}\|=0
\end{align*}
where the limit is taken under the hypotheses of Lemma~\ref{weightedlimit}, and $F_{\sigma,0}^{-1}-I$ is compact.
\end{lemma}
\begin{proof}
See Appendix~\ref{proof:weightedinverse}.
\end{proof}

We summarize this conclusion below for use in the localisation argument. It allows us to compare exterior observations while retaining the noncircular geometry, rather than replacing the physical shell by its circular comparator.

\begin{lemma}[Weighted Fredholm comparison]\label{weightedcomparison}
Assume $\sigma\ge1/2$, $\eta<\rho^{2\sigma}$, $\zeta\to0$ in the upper half-plane, and $a^{2\sigma}/y=O(1)$. There are bounded operators $F_{\sigma,\zeta}$ such that
\begin{equation}
 \TT^\sigma V=F_{\sigma,\zeta}^{-1}\TT^\sigma V_0 \qquad \text{ and }
 \qquad \lim_{\zeta\to0}\|F_{\sigma,\zeta}^{-1}-F_{\sigma,0}^{-1}\|=0,
 \label{weightedcomparisonidentity}
\end{equation}
where the limiting inverse is bounded and differs from identity by a compact operator.
\end{lemma}
\begin{proof}
Lemma~\ref{weightedlimit} gives the operator-norm limit. Lemma~\ref{weightedkernel} rules out an additional kernel introduced by the weights. Lemma~\ref{weightedinverse} gives~\eqref{weightedcomparisonidentity} and the compactness assertion.
\end{proof}

\section{Proof of sharp normalized localisation}\label{sec:localisation}
We can now return to the spatial meaning of resonance. Along a sequence with divergent shell power, normalization by $\sqrt E$ measures the exterior field relative to the growing absorption. We seek the region where this normalized field disappears, and ask how close its boundary can lie to the shell. The balance between resonant excitation and outward attenuation will identify the critical curve and explain why a nonzero integrated energy can remain at the critical radius.

\subsection{Sufficiency and local convergence}
We begin outside the critical curve, where distance suppresses the fine angular structure of the resonant field. In conformal coordinates this suppression is measured by $r^{-n}$, so the modes most strongly attenuated are precisely those moving to higher indices as loss decreases. Beyond the critical radius, attenuation is sufficient to make the power-normalized field vanish. This describes a change in the observed asymptotic behavior, not a material interface or a discontinuity of the field.

\begin{lemma}[Exterior decay above the critical radius]\label{localsufficiency}
Under the standing hypotheses, let $(\ell_j)$ be positive with $\lim_{j\to+\infty}\ell_j=0$ and $\lim_{j\to+\infty}E(\ell_j)=+\infty$. Then,
\begin{align*}
 &\lim_{j\to+\infty}\frac{S_r(\ell_j)}{E(\ell_j)}=0
 \qquad\text{for every fixed }\ r>r_c
\end{align*}
and
\begin{align*}
 &\lim_{j\to+\infty}\frac{u_{\ell_j}}{\sqrt{E(\ell_j)}}=0
 \qquad\text{in } \ H^1_{\mathrm{loc}}(\mathcal O_{r_c}).
\end{align*}
\end{lemma}
\begin{proof}
See Appendix~\ref{proof:localsufficiency}.
\end{proof}

\subsection{Power concentration and lower bounds}
The preceding result does not yet tell us whether the same decay persists closer to the shell. At resonant losses, increasingly fine modes remain strong enough to survive power normalization in nearer observation regions. The lower bounds below make this distinction precise: they show that moving the critical curve inward would include regions where normalized localisation fails. In this sense, the critical boundary reflects the field itself, rather than merely a limitation of the upper estimate.

\begin{lemma}[Concentration and lower bounds]\label{localsharpness}
Under the standing hypotheses and in the resonant source range of Theorem~\ref{mainlocal}, use the exact crossings~\eqref{crossings} if $p<1$, and $\ell_n=t_n$ if $p\ge1$. In both cases $\lim_{n\to+\infty}E(\ell_n)=+\infty$. For $\beta=2\log r/\log(1/\rho)$,
\begin{align*}
 \frac{S_r(\ell_n)}{E(\ell_n)}\gtrsim t_n^{\beta-1/\tau}.
\end{align*}
Moreover, for every $1<r_0<r_c$, the same sequence prevents local $H^1$ convergence of $u_\ell/\sqrt E$ to zero on $\mathcal O_{r_0}$.
\end{lemma}
\begin{proof}
See Appendix~\ref{proof:localsharpness}.
\end{proof}

\subsection{Exact tangential endpoint}
At the critical radius itself, neither resonant amplification nor spatial attenuation dominates the normalized energy. Along tangential crossings, their balance leaves the finite, nonzero limit below. Its independence from the source amplitude reflects the fact that both the exterior energy and the shell power measure the same leading excitation. The source must nevertheless lie in the resonant range, so that the power used for normalization actually diverges.

\begin{lemma}[Critical energy along tangential crossings]\label{localendpoint}
Under the standing hypotheses with $0<p<1$ and $0<\alpha\le1/p$, the exact crossings~\eqref{crossings} satisfy
\begin{align*}
 \lim_{n\to+\infty}\frac{S_{r_c}(\ell_n)}{E(\ell_n)}
 =\frac12\left(\frac{|c|}{2}\right)^{1/p}.
\end{align*}
\end{lemma}
\begin{proof}
See Appendix~\ref{proof:localendpoint}.
\end{proof}

\begin{proof}[Proof of Theorem~\ref{mainlocal}]
Lemma~\ref{localsufficiency} proves~\eqref{mainupper} along every divergent-power sequence. In the resonant source range, Lemma~\ref{localsharpness} supplies one such sequence. If $r<r_c$, then, $\beta<1/\tau$ and its lower bound diverges, whereas if $r=r_c$, the exponent is zero and the lower limit is positive. Together with sufficiency, this proves the sharp full-exterior seminorm criterion. The local lower bound also rules out a smaller open exterior region. Local convergence in $\mathcal O_{r_c}$ is consistent with nonvanishing energy integrated over the whole of that open region, because its inner boundary is not contained in any compact subset.
\end{proof}

\begin{proof}[Proof of Theorem~\ref{mainendpoint}]
Lemma~\ref{localsharpness} ensures that the crossing sequence has divergent power in the stated source range. Lemma~\ref{localendpoint} gives the limit~\eqref{endpoint}.
\end{proof}

\section{Numerical illustrations}\label{numerics}
The analysis distinguishes the power absorbed by the shell from the field observed outside it. We now illustrate both quantities at finite loss, first comparing the power regimes of Theorem~\ref{mainpower}, then, examining the critical observation boundary in Theorems~\ref{mainlocal}--\ref{mainendpoint}. Spatial maps complete this comparison by showing where the dissipation and normalized exterior intensity are concentrated. Each experiment keeps the smooth source fixed as the loss decreases. These computations illustrate the proved asymptotic results rather than establish them. Implementation details and accuracy checks are collected in Appendix~\ref{numericalappendix}.

\subsection{Geometry and observation regions}
\begin{figure}
\centering
\includegraphics[width=.5\linewidth]{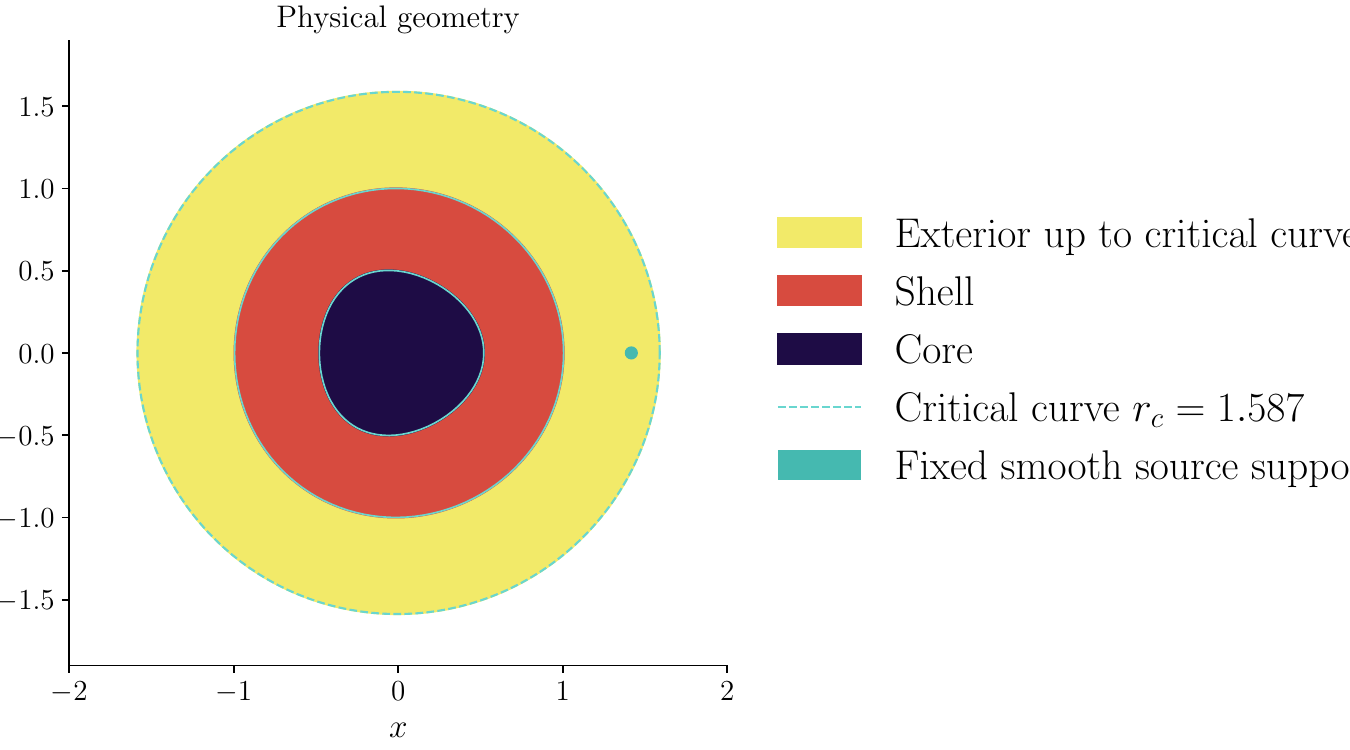}
\caption{The core, shell and exterior up to the critical curve for $p=0.75$, $\alpha=1$. The dashed curve is $\Psi(r_ce^{i\theta})$, with $r_c=2^{2/3}\simeq1.5874$. The small disk marks the smooth source support. The deformation is not magnified.}\label{fig:geometry}
\end{figure}

We first specify the geometry so that the observation curves in the subsequent figures can be distinguished from the material interfaces. Throughout this section,
$\rho=0.5, \varepsilon=0.005, c=1$ and $\Psi(w)=w+0.005w^{-2}$.
The standing geometric condition holds for both path exponents used below, $p=0.75$ and $p=1.2$. Figure~\ref{fig:geometry} shows the physical configuration for the intermittent case $p=0.75$, $\alpha=1$, with source center $\Psi(\sqrt2)$. The dashed critical curve in Figure~\ref{fig:geometry} separates observation radii for which the normalized exterior energy decays from those for which it need not decay.

\subsection{Power classification}
\begin{figure}
\centering
\includegraphics[width=.75\linewidth]{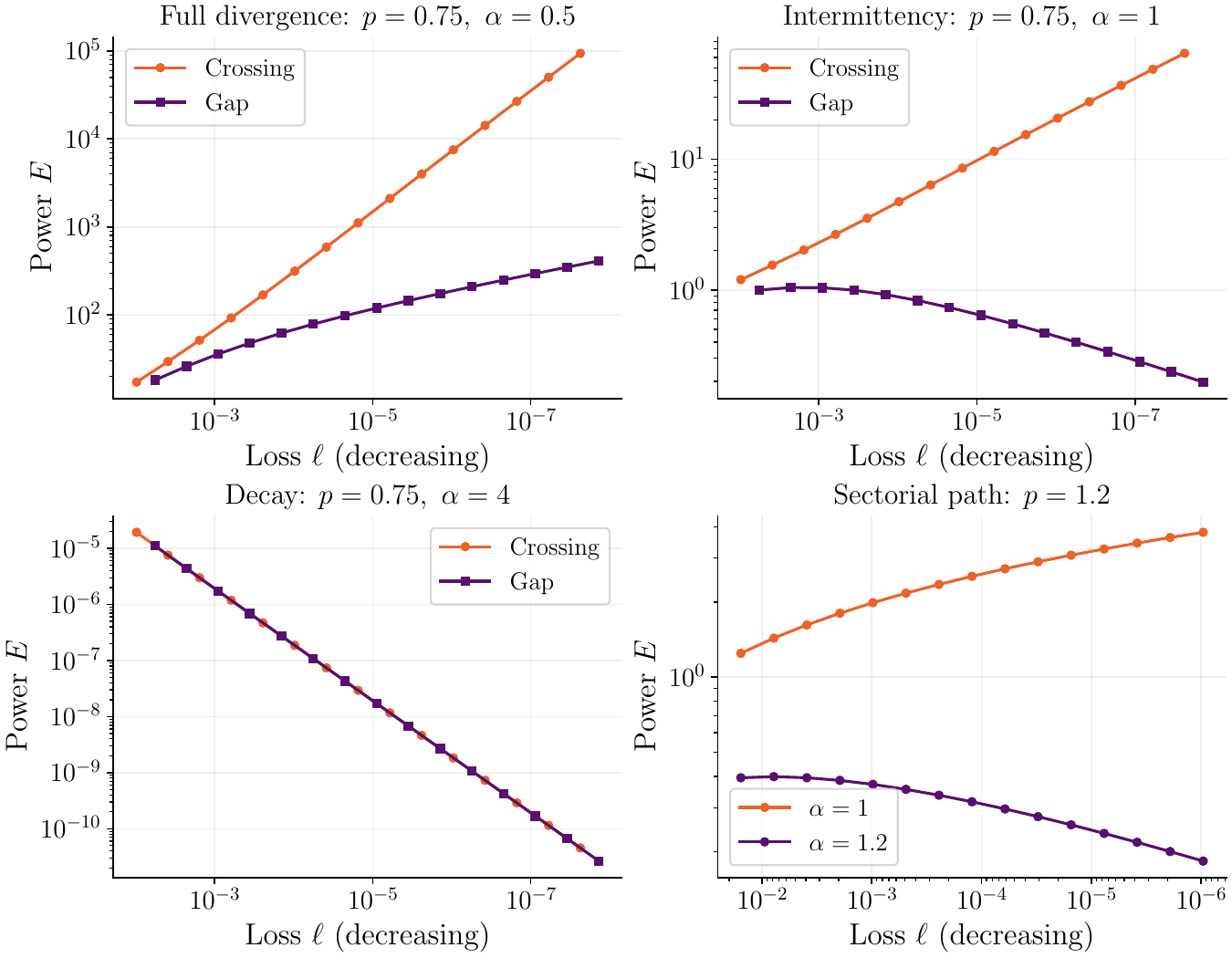}
\caption{Dissipated power plots.
The tangential panels compare reference crossings and gaps. The sectorial panel uses $\ell_n=\rho^n$. Loss decreases to the right.} \label{fig:regimes}
\end{figure}

We begin with the distinction between a growing background response and isolated growing peaks. For $p=0.75$, Theorem~\ref{mainpower} gives source thresholds $2/3$ and $4/3$. We choose one source exponent in each regime and compare losses at reference crossings with losses in the intervening gaps. The upper-left, upper-right and lower-left panels of Figure~\ref{fig:regimes} use $\alpha=0.5$, $1$ and $4$, respectively: crossing and gap powers both grow in the first case, separate into growing peaks and diminishing gaps in the second, and both decay in the third. Thus, the sampled curves make visible the background and peak mechanisms identified in Section~\ref{sec:power}. The lower-right panel treats the sectorial path $p=1.2$, for which the response does not split into narrow, separated peaks: slow growth at $\alpha=1$ is consistent with the logarithmic threshold law, whereas $\alpha=1.2$ gives decay. 

\subsection{The critical observation boundary}
\begin{figure}
\centering
\includegraphics[width=.7\linewidth]{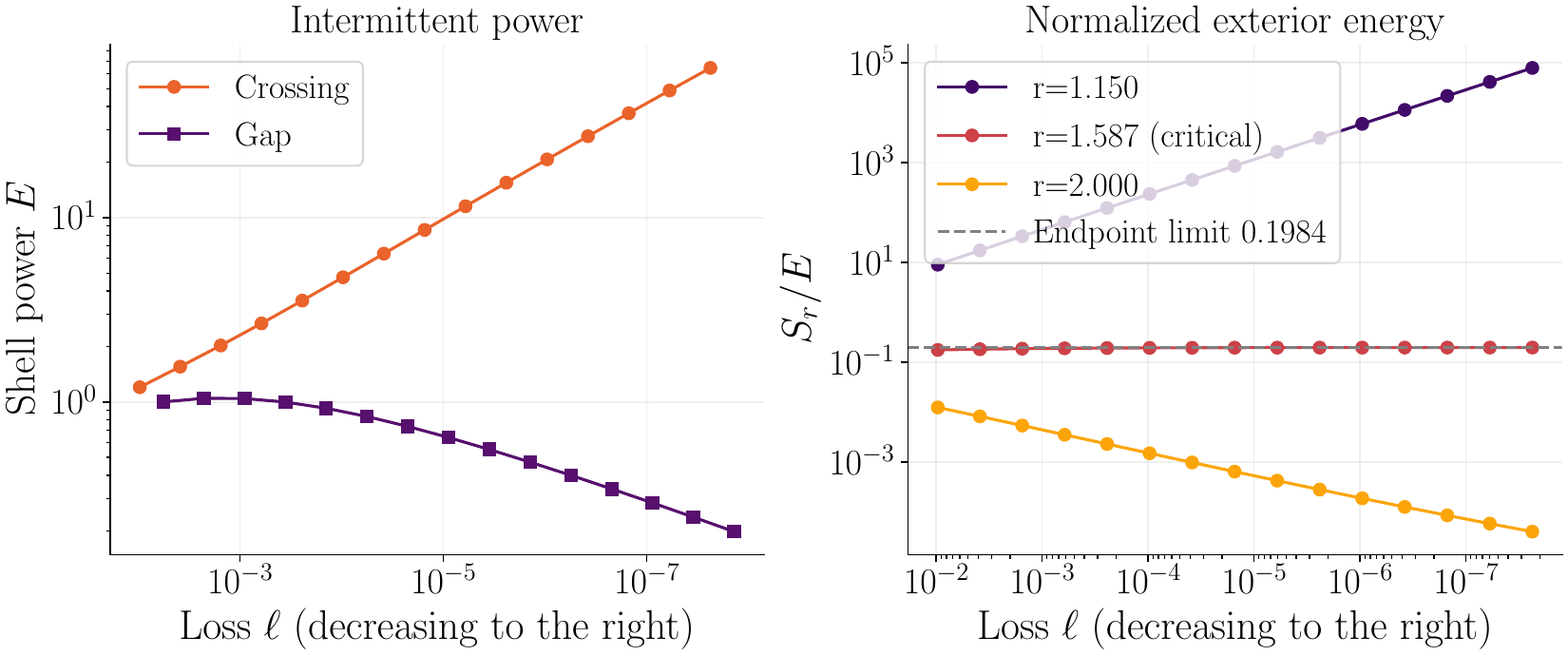}
\caption{Intermittent power and normalized exterior energy plots.
The right panel uses crossing samples at conformal radii $r=1.15,r_c,2$. The dashed horizontal line is the exact endpoint prediction. Physical observation boundaries are the corresponding images under $\Psi$.}\label{fig:ratios}
\end{figure}

Having identified an intermittent response, we next ask where its resonant field remains visible after normalization by the shell power. We retain the same source and compare observation regions below, at and above the critical radius. The left panel of Figure~\ref{fig:ratios} recalls the crossing and gap powers, while the right panel follows $S_r/E$ along the crossing sequence. The quotient grows at $r=1.15$, approaches a nonzero value at $r=r_c$, and decreases at $r=2$, illustrating the sharp distinction in Theorem~\ref{mainlocal}. The critical curve also permits a quantitative comparison, since Theorem~\ref{mainendpoint} predicts
\begin{align*}
 \lim_{n\to+\infty}\frac{S_{r_c}(\ell_n)}{E(\ell_n)}
 =\frac12\left(\frac12\right)^{4/3}\simeq0.198425.
\end{align*}
The last computed value on the critical-radius curve in Figure~\ref{fig:ratios} (right) is $0.198412$, close to the dashed horizontal line representing this limit. This is an integrated energy statement and not an abrupt pointwise disappearance of the field at the critical curve.

\subsection{Spatial concentration}
\begin{figure}
\centering
\includegraphics[width=0.8\linewidth]{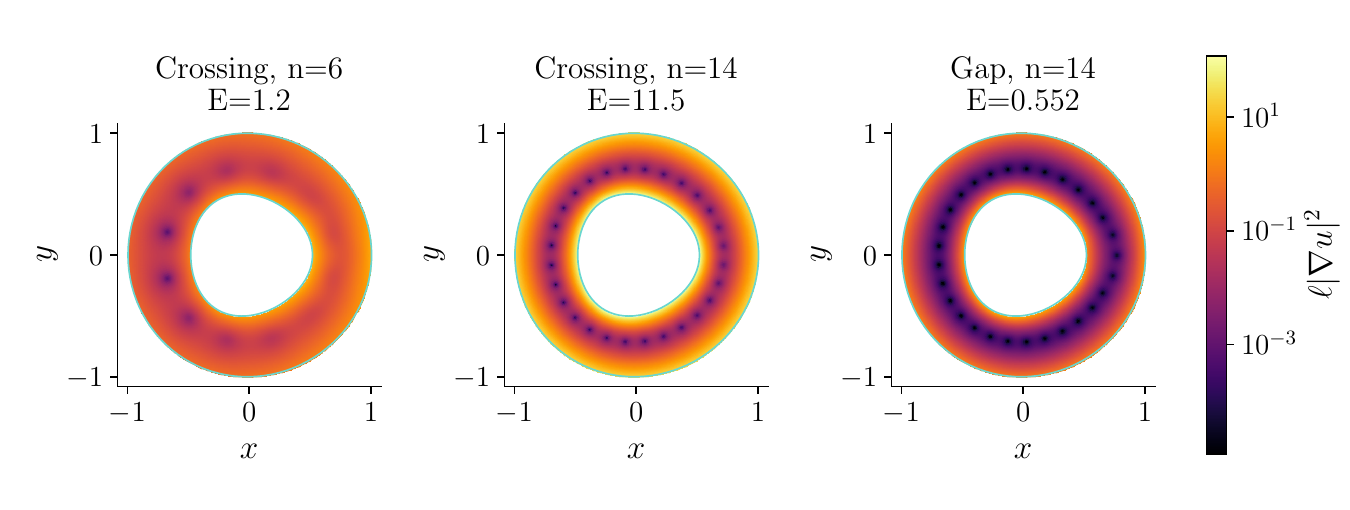}
\caption{Shell dissipation density $D_\ell$ at two reference crossings and one gap. The geometry is shown without magnifying its deformation.}\label{fig:shelldensity}
\end{figure}

\begin{figure}
\centering
\includegraphics[width=0.8\linewidth]{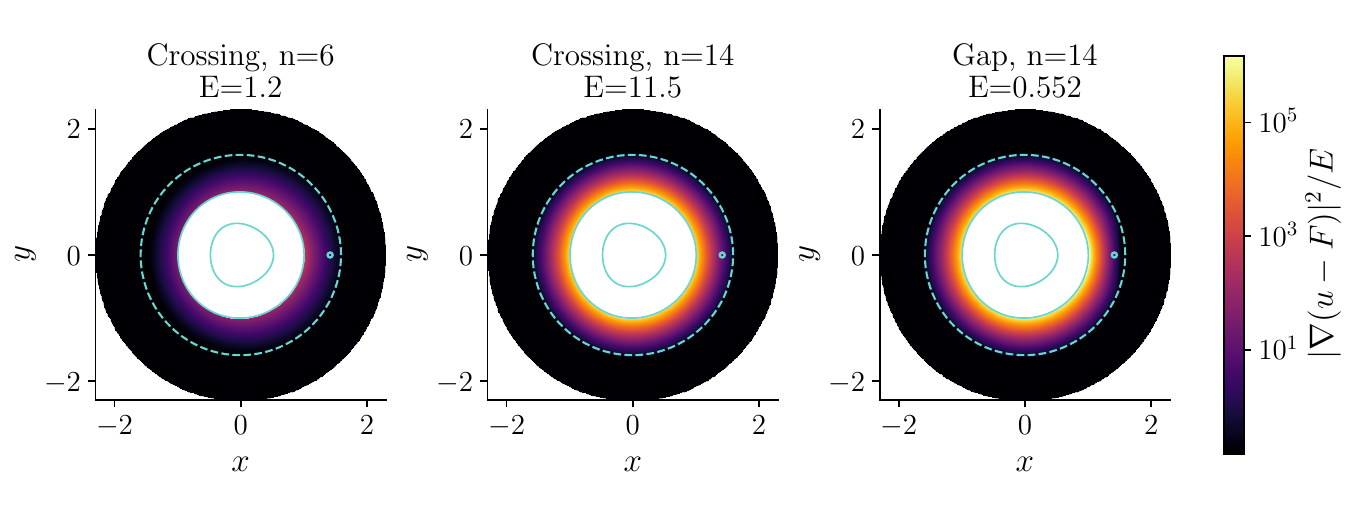}
\caption{Normalized exterior intensity $Q_\ell$ at the same samples as Figure~\ref{fig:shelldensity}. The dashed curve is the critical conformal level curve, and the small disk marks source support.} \label{fig:exteriordensity}
\end{figure}

The preceding plots compare integrated quantities. To see how those quantities are distributed in space, we now examine the same intermittent experiment at crossings $n=6$ and $14$ and the gap associated with $n=14$. Figures~\ref{fig:shelldensity} and~\ref{fig:exteriordensity} separate absorption within the shell from the normalized scattered-field intensity outside it. The plotted quantities are $\ell|\nabla u_\ell(z)|^2$ for $z\in\Omega_s$ in Figure~\ref{fig:shelldensity} and  $\frac{|\nabla(u_\ell-F)(z)|^2}{E(\ell)}$ for $z\notin\overline{\Omega_e}$ in Figure~\ref{fig:exteriordensity}. Their integrals are $E$ over the shell and $S_r/E$ over $\mathcal O_r$, respectively. Comparing the left and middle panels of Figure~\ref{fig:shelldensity}, the finer angular structure at the later crossing is consistent with excitation of higher-order modes. The common shell scale in Figure~\ref{fig:shelldensity} also displays the weaker absorption in the gap shown in its right panel.

The exterior maps in Figure~\ref{fig:exteriordensity} show where the normalized scattered-gradient intensity is concentrated. They complement the integrated comparison in Figure~\ref{fig:ratios} as the dashed critical curve is not a discontinuity of the field. In particular, a bright gap field need not indicate greater absorption, since normalization divides by a smaller $E$.

\section*{Declaration of AI assistance}
The author used OpenAI's ChatGPT/Codex to assist with literature searches and reference tracking, to assist with implementation and optimization of the numerical code, and to improve the language of the manuscript. The author takes full responsibility for the mathematical results, proofs, computations, references and final content.

\appendix
\section{Proofs of Section~\ref{boundaryformulation}}\label{app:layers}
Throughout this appendix, normal traces and boundary integrals are understood in the $H^{-1/2}$--$H^{1/2}$ duality. Green identities may first be established for smooth data and subsequently extended by continuity. Distributional pairings are bilinear and complex conjugation is written explicitly in energy identities.

\subsection{Proof of Lemma~\ref{layerflux}}\label{proof:layerflux}
\begin{proof}
Fix $\ell>0$ and write $s=s_\ell$, $\lambda=\lambda_\ell$. Then, $s\ne0$ and $s-1\ne0$.
Let $\boldsymbol\varphi=(\varphi_i,\varphi_e)\in\mathcal X$ and define
\begin{align*}
 v=F+\mathcal S_i[\varphi_i]+\mathcal S_e[\varphi_e],
 \qquad g_j=\frac{\partial}{\partial \nu_j}F|_{\Gamma_j},\quad j\in\{i,e\}.
\end{align*}
The cross-interface single layer is smooth across the other interface, and $F$ is harmonic near both interfaces. Consequently the four normal traces, with the signs fixed in \eqref{normaljump}, are
\begin{align}
 \frac{\partial}{\partial \nu_i}v|_+
 &=g_i+\left(\tfrac12I+\mathcal K_i^*\right)\varphi_i
          +\mathcal N_{ie}\varphi_e,\label{prooftraceiplus}\\
 \frac{\partial}{\partial \nu_i}v|_-
 &=g_i+\left(-\tfrac12I+\mathcal K_i^*\right)\varphi_i
          +\mathcal N_{ie}\varphi_e,\label{prooftraceiminus}\\
 \frac{\partial}{\partial \nu_e}v|_+
 &=g_e+\mathcal N_{ei}\varphi_i
          +\left(\tfrac12I+\mathcal K_e^*\right)\varphi_e,\label{prooftraceeplus}\\
 \frac{\partial}{\partial \nu_e}v|_-
 &=g_e+\mathcal N_{ei}\varphi_i
          +\left(-\tfrac12I+\mathcal K_e^*\right)\varphi_e.\label{prooftraceeminus}
\end{align}
On $\Gamma_i$, the flux condition in \eqref{transmission} is
$s\,\frac{\partial}{\partial \nu_i}v|_+-\frac{\partial}{\partial \nu_i}v|_-=0$. Substitution of
\eqref{prooftraceiplus}--\eqref{prooftraceiminus} gives
\begin{align*}
 (s-1)g_i+\tfrac{s+1}{2}\varphi_i
       +(s-1)\mathcal K_i^*\varphi_i
       +(s-1)\mathcal N_{ie}\varphi_e = 0.
\end{align*}
Dividing by $s-1$ and using $\lambda=(s+1)/(2(s-1))$ yields
\begin{equation}
 (-\lambda I-\mathcal K_i^*)\varphi_i-\mathcal N_{ie}\varphi_e=g_i.
 \label{innerBIE}
\end{equation}
On $\Gamma_e$, the flux condition is
$\frac{\partial}{\partial \nu_e}v|_+-s\,\frac{\partial}{\partial \nu_e}v|_-=0$. This time
\eqref{prooftraceeplus}--\eqref{prooftraceeminus} give
\begin{align*}
 (1-s)g_e+(1-s)\mathcal N_{ei}\varphi_i
       +\tfrac{1+s}{2}\varphi_e+(1-s)\mathcal K_e^*\varphi_e = 0.
\end{align*}
After division by $s-1$, this is
$-g_e-\mathcal N_{ei}\varphi_i+\lambda\varphi_e-\mathcal K_e^*\varphi_e=0$.
Thus,
\begin{equation}
 -\mathcal N_{ei}\varphi_i+(\lambda I-\mathcal K_e^*)\varphi_e=g_e.
 \label{outerBIE}
\end{equation}
Multiplying \eqref{innerBIE} by $-1$ and combining it with \eqref{outerBIE} gives
\begin{equation}
 \begin{pmatrix}
 \lambda I+\mathcal K_i^*&\mathcal N_{ie}\\
 -\mathcal N_{ei}&\lambda I-\mathcal K_e^*
 \end{pmatrix}
 \binom{\varphi_i}{\varphi_e}
 =
 \binom{-g_i}{g_e}
 \quad\Longleftrightarrow\quad
 (\lambda I-\mathbb K^*)\boldsymbol\varphi=\mathbf g_F.
 \label{proofblocksystem}
\end{equation}
Continuity of $v$ follows separately from \eqref{Sjump}: for example,
\begin{align*}
 v|_{\Gamma_i,+}-v|_{\Gamma_i,-}
 =\mathcal S_i[\varphi_i]|_+-\mathcal S_i[\varphi_i]|_-=0,
\end{align*}
since the other two terms have equal traces. The same identity holds on $\Gamma_e$.
\end{proof}

\subsection{Proof of Lemma~\ref{layercharges}}\label{proof:layercharges}
\begin{proof}
Write $g_j=\frac{\partial}{\partial \nu_j}F|_{\Gamma_j}$ and $\lambda=\lambda_\ell$.
First, for $j=i,e$, the divergence theorem and $\Delta F=0$ in $\Omega_j$ give
\begin{equation}
 \langle g_j,1\rangle_{\Gamma_j}
 =\int_{\Gamma_j}\frac{\partial}{\partial \nu_j}F\,ds
 =\int_{\Omega_j}\Delta F\,dx=0.
 \label{proofmeanforcing}
\end{equation}
hence, $\mathbf g_F\in\mathcal X$.

We next establish the mean identities for arbitrary densities, without yet assuming their charges vanish. Let $\psi_j\in H^{-1/2}(\Gamma_j)$. Since $\mathcal S_j[\psi_j]$ is harmonic inside $\Omega_j$, the interior jump formula gives
\begin{align}
\begin{aligned}
 0&=\int_{\Omega_j}\Delta\mathcal S_j[\psi_j]\,dx
   =\left\langle\frac{\partial}{\partial \nu_j}\mathcal S_j[\psi_j]|_-,1\right\rangle_{\Gamma_j}\\
 &=\left\langle(-\tfrac12I+\mathcal K_j^*)\psi_j,1\right\rangle_{\Gamma_j}
   =-\tfrac12\langle\psi_j,1\rangle_{\Gamma_j}
       +\langle\mathcal K_j^*\psi_j,1\rangle_{\Gamma_j}.
 \label{proofdiagonalmean}
\end{aligned}
\end{align}
Therefore,
\begin{equation}
 \langle\mathcal K_j^*\psi_j,1\rangle_{\Gamma_j}
 =\tfrac12\langle\psi_j,1\rangle_{\Gamma_j}.
 \label{proofKcharge}
\end{equation}
For the map from the outer interface to the inner one, the potential
$w_e=\mathcal S_e[\psi_e]$ is harmonic in a neighborhood of
$\overline{\Omega_i}$, so
\begin{equation}
 \langle\mathcal N_{ie}\psi_e,1\rangle_{\Gamma_i}
 =\int_{\Gamma_i}\frac{\partial}{\partial \nu_i}w_e\,ds
 =\int_{\Omega_i}\Delta w_e\,dx=0.
 \label{proofNiezero}
\end{equation}
For the opposite map let $w_i=\mathcal S_i[\psi_i]$. It is harmonic in
$\Omega_s$. Its shell normal is $\nu_e$ on $\Gamma_e$ and $-\nu_i$ on $\Gamma_i$.
Consequently,
\begin{align}
\begin{aligned}
 0
 &=\int_{\Omega_s}\Delta w_i\,dx
   =\int_{\Gamma_e}\frac{\partial}{\partial \nu_e}w_i\,ds
       -\left\langle\frac{\partial}{\partial \nu_i}w_i|_+,1\right\rangle_{\Gamma_i}\\
 &=\langle\mathcal N_{ei}\psi_i,1\rangle_{\Gamma_e}
       -\left\langle(\tfrac12I+\mathcal K_i^*)\psi_i,1\right\rangle_{\Gamma_i}.
 \label{proofshellmean}
\end{aligned}
\end{align}
Using \eqref{proofKcharge}, we obtain the exact identity
\begin{equation}
 \langle\mathcal N_{ei}\psi_i,1\rangle_{\Gamma_e}
 =\tfrac12\langle\psi_i,1\rangle_{\Gamma_i}
      +\langle\mathcal K_i^*\psi_i,1\rangle_{\Gamma_i}
 =\langle\psi_i,1\rangle_{\Gamma_i}.
 \label{proofNeicharge}
\end{equation}
These are weak Green identities. Equivalently, one can first prove them for smooth densities and extend by continuity of the operators and of the charge functional.

Now let $\boldsymbol\psi=(\psi_i,\psi_e)\in\mathcal X$, so both component charges vanish. Equations \eqref{proofKcharge}, \eqref{proofNiezero}, and \eqref{proofNeicharge} give
\begin{align}
 \left\langle(\mathbb K^*\boldsymbol\psi)_i,1\right\rangle_{\Gamma_i}
 &=-\langle\mathcal K_i^*\psi_i,1\rangle_{\Gamma_i}
      -\langle\mathcal N_{ie}\psi_e,1\rangle_{\Gamma_i}
   =-\tfrac12\langle\psi_i,1\rangle_{\Gamma_i}-0=0,\nonumber\\
 \left\langle(\mathbb K^*\boldsymbol\psi)_e,1\right\rangle_{\Gamma_e}
 &=\langle\mathcal N_{ei}\psi_i,1\rangle_{\Gamma_e}
      +\langle\mathcal K_e^*\psi_e,1\rangle_{\Gamma_e}
   =\langle\psi_i,1\rangle_{\Gamma_i}
       +\tfrac12\langle\psi_e,1\rangle_{\Gamma_e}=0.
 \label{proofinvariant}
\end{align}
This proves $\mathbb K^*(\mathcal X)\subset\mathcal X$. The space $\mathcal X$ is closed, since the charge functionals are continuous. The compact diagonal and cross-interface blocks, therefore, give a compact restriction $\mathbb K^*:\mathcal X\to\mathcal X$. As $\lambda\ne0$,
\begin{equation}
 \lambda I-\mathbb K^*=\lambda(I-\lambda^{-1}\mathbb K^*)
 \label{prooffredholm}
\end{equation}
is a Fredholm operator of index zero on $\mathcal X$.
\end{proof}

\subsection{Proof of Lemma~\ref{layerdecay}}\label{proof:layerdecay}
\begin{proof}
We record the estimates needed before taking an infinite-domain energy limit. Choose $R_0$ with $\overline{\Omega_e}\subset B_{R_0}$. For $|x|\ge2R_0$ and $y\in\Gamma_j$,
\begin{align*}
 \Phi(x-y)-\Phi(x)
 &=\frac1{4\pi}\log\left(1-\frac{2x\cdot y}{|x|^2}
                              +\frac{|y|^2}{|x|^2}\right)
   =O(|x|^{-1}),\\
 \nabla_x\Phi(x-y)-\nabla_x\Phi(x)
 &=\frac1{2\pi}\left(\frac{x-y}{|x-y|^2}-\frac{x}{|x|^2}\right)
   =O(|x|^{-2}).
\end{align*}
The bounds hold also after one tangential derivative in $y$, uniformly on the fixed curves. Hence, they hold in $H^{1/2}(\Gamma_j)$ as functions of $y$. If $\langle\psi_j,1\rangle_{\Gamma_j}=0$, then, the constant terms cancel:
\begin{align}
 \mathcal S_j[\psi_j](x)
 &=\langle\psi_j,\Phi(x-\cdot)-\Phi(x)\rangle_{\Gamma_j}
   =O(|x|^{-1}),\nonumber\\
 \nabla\mathcal S_j[\psi_j](x)
 &=\langle\psi_j,\nabla_x\Phi(x-\cdot)-\nabla_x\Phi(x)\rangle_{\Gamma_j}
   =O(|x|^{-2}).
 \label{prooflayerdecay}
\end{align}
In particular, their sum satisfies the decay condition in \eqref{transmission}, together with the derivative estimate needed for Green's identity. The gradient tail has finite energy because
$\int_{2R_0}^{+\infty} r^{-4}r\,dr<+\infty$.
\end{proof}

\subsection{Proof of Lemma~\ref{layerkernel}}\label{proof:layerkernel}
\begin{proof}
Write $s=s_\ell$, $\lambda=\lambda_\ell$, and choose $R_0$ with $\overline{\Omega_e}\subset B_{R_0}$. Boundary energy integrals denote weak normal-trace pairings when the data are not smooth.
Let $\boldsymbol\psi=(\psi_i,\psi_e)\in\mathcal X$ satisfy
$(\lambda I-\mathbb K^*)\boldsymbol\psi=0$, and set
\begin{align*}
 U=\mathcal S_i[\psi_i]+\mathcal S_e[\psi_e].
\end{align*}
Lemma~\ref{layerflux} with $F=0$ shows that $U$ is continuous across both interfaces, is harmonic in each region, and satisfies
\begin{equation}
 s\,\frac{\partial}{\partial \nu_i}U|_+=\frac{\partial}{\partial \nu_i}U|_- \qquad \text{ and } 
 \qquad
 \frac{\partial}{\partial \nu_e}U|_+=s\,\frac{\partial}{\partial \nu_e}U|_-.
 \label{proofhomogeneousflux}
\end{equation}
It has finite energy by the local mapping properties of the single layers and \eqref{prooflayerdecay}.

For $R>2R_0$, write $D_R=B_R\setminus\overline{\Omega_e}$ and let $\nu_R$ be the outward normal of $B_R$. Applying the first Green identity in each region gives
\begin{align}
 \int_{\Omega_i}|\nabla U|^2
 &=\int_{\Gamma_i}\overline U\,\frac{\partial}{\partial \nu_i}U|_-\,ds,\label{proofcoreenergy}\\
 s\int_{\Omega_s}|\nabla U|^2
 &=s\int_{\Gamma_e}\overline U\,\frac{\partial}{\partial \nu_e}U|_-\,ds
     -s\int_{\Gamma_i}\overline U\,\frac{\partial}{\partial \nu_i}U|_+\,ds,\label{proofshellenergy}\\
 \int_{D_R}|\nabla U|^2
 &=\int_{\partial B_R}\overline U\,\frac{\partial}{\partial \nu_R}U\,ds
     -\int_{\Gamma_e}\overline U\,\frac{\partial}{\partial \nu_e}U|_+\,ds.\label{proofexteriorenergy}
\end{align}
Continuity makes $\overline U$ the same trace in the two terms on each interface. Summing the three equations and grouping the interface terms, we obtain
\begin{align}
 \int_{\Omega_i}|\nabla U|^2+s\int_{\Omega_s}|\nabla U|^2
                         +\int_{D_R}|\nabla U|^2 =
 \int_{\partial B_R}\overline U\,\frac{\partial}{\partial \nu_R}U\,ds,
 \label{prooffiniteenergyidentity}
\end{align}
where we have used the two identities in \eqref{proofhomogeneousflux}.

By \eqref{prooflayerdecay}, there is a constant $C$ independent of $R$ such that
\begin{equation}
 \left|\int_{\partial B_R}\overline U\,\frac{\partial}{\partial \nu_R}U\,ds\right|
 \le 2\pi R\,\sup_{\partial B_R}|U|\,
                \sup_{\partial B_R}|\nabla U|
 \le 2\pi R(CR^{-1})(CR^{-2})
 =2\pi C^2R^{-2}.
 \label{proofboundaryvanish}
\end{equation}
Since $\lim_{R\to+\infty}2\pi C^2R^{-2}=0$, taking the limit as $R$ tends to infinity in \eqref{prooffiniteenergyidentity}, with the finite-energy exterior integral converging monotonically, gives
\begin{equation}
 I_i+sI_s+I_e=0,\qquad
 I_i=\int_{\Omega_i}|\nabla U|^2,\quad
 I_s=\int_{\Omega_s}|\nabla U|^2,\quad
 I_e=\int_{\mathbb R^2\setminus\overline{\Omega_e}}|\nabla U|^2.
 \label{proofglobalenergy}
\end{equation}
All $I_j$ are real and nonnegative. Since $s=-1+c\ell^p+i\ell$, its imaginary and real parts give, in that order,
\begin{align}
 0&=\operatorname{Im}(I_i+sI_s+I_e)=\ell I_s
       \quad\Longrightarrow\quad I_s=0,\label{proofimaginary}\\
 0&=\operatorname{Re}(I_i+sI_s+I_e)
       =I_i+(-1+c\ell^p)I_s+I_e=I_i+I_e
       \quad\Longrightarrow\quad I_i=I_e=0.\label{proofreal}
\end{align}
The vanishing gradients, continuity across both interfaces and decay~\eqref{prooflayerdecay} imply $U=0$ everywhere.

Finally the other-interface layer has no normal derivative jump across $\Gamma_j$, while the self-layer jump is
\begin{equation}
 \frac{\partial}{\partial \nu_j}U|_+-\frac{\partial}{\partial \nu_j}U|_-
 =\left[(\tfrac12I+\mathcal K_j^*)-
             (-\tfrac12I+\mathcal K_j^*)\right]\psi_j
 =\psi_j,\qquad j=i,e.
 \label{proofdensityjump}
\end{equation}
Since $U=0$, the left side vanishes, so $\psi_i=\psi_e=0$. This proves injectivity. By the index-zero Fredholm property \eqref{prooffredholm}, the operator is also surjective, its inverse is bounded by the bounded inverse theorem.
\end{proof}

\subsection{Proof of Lemma~\ref{layerconstruction}}\label{proof:layerconstruction}
\begin{proof}
Write $s=s_\ell$ and $\lambda=\lambda_\ell$.
Take the density pair in the statement and set $u_\ell=F+\mathcal S_i[\varphi_i]+\mathcal S_e[\varphi_e]$.
The layers are harmonic off their own interfaces. Since $\Delta F=f$ and $f=0$ in a neighborhood of $\overline{\Omega_e}$, this field satisfies the volume equations in \eqref{transmission}. Lemma~\ref{layerflux} proves the two continuity and two flux conditions. The local single-layer mapping properties give $u_\ell\in H^1_{\rm loc}(\mathbb R^2)$. Equation \eqref{prooflayerdecay} and the analogous estimates for $F$ give
\begin{align*}
 u_\ell(x)=O(|x|^{-1})\qquad \text{ and } \qquad \nabla u_\ell(x)=O(|x|^{-2}).
\end{align*}
For $F$ these follow from the same kernel differences integrated against $f$, using $\int f=0$ and compact support. Hence, $\nabla u_\ell\in L^2(\mathbb R^2)$ and the prescribed decay \eqref{pde} holds. To check the global equation directly, the volume equations and continuity give the possible interface measures in $\nabla\cdot(a_\ell\nabla u_\ell)$ with coefficients
\begin{align*}
 s\,\frac{\partial}{\partial \nu_i}u_\ell|_+-\frac{\partial}{\partial \nu_i}u_\ell|_-=0
 \qquad \text{ and } \qquad
 \frac{\partial}{\partial \nu_e}u_\ell|_+-s\,\frac{\partial}{\partial \nu_e}u_\ell|_-=0.
\end{align*}
Thus, there are no interface source terms and \eqref{pde} holds distributionally.
\end{proof}

\subsection{Proof of Lemma~\ref{layerconverse}}\label{proof:layerconverse}
\begin{proof}
Write $s=s_\ell$ and $\lambda=\lambda_\ell$. Distributional pairings below are bilinear, without complex conjugation.
Let $v$ be any solution of \eqref{transmission} in the stated energy class, with decay \eqref{pde}. For $j\in\{i,e\}$, set
\begin{align*}
 q_j^\pm=\frac{\partial}{\partial \nu_j}v|_\pm\in H^{-1/2}(\Gamma_j)
 \qquad \text{ and } \qquad \varphi_j=q_j^+-q_j^-.
\end{align*}
These weak normal traces exist because $v$ is harmonic on the two sides of each interface in a neighborhood of that interface. The transmission conditions read
\begin{equation}
 q_i^-=s q_i^+\qquad \text{ and } \qquad q_e^+=s q_e^-.
 \label{proofconverseflux}
\end{equation}
By the core divergence theorem,
\begin{align*}
 \langle q_i^-,1\rangle_{\Gamma_i}
 =\int_{\Omega_i}\Delta v\,dx=0.
\end{align*}
Since $s\ne0$, \eqref{proofconverseflux} yields
\begin{equation}
 \langle q_i^+,1\rangle_{\Gamma_i}
 =s^{-1}\langle q_i^-,1\rangle_{\Gamma_i}=0 \qquad \text{ and } \qquad
 \langle\varphi_i,1\rangle_{\Gamma_i}
 =\langle q_i^+-q_i^-,1\rangle_{\Gamma_i}=0.
 \label{proofconverseinnercharge}
\end{equation}
Using the outward shell normals gives
\begin{align*}
 0=\int_{\Omega_s}\Delta v\,dx
   =\langle q_e^-,1\rangle_{\Gamma_e}
       -\langle q_i^+,1\rangle_{\Gamma_i}
   =\langle q_e^-,1\rangle_{\Gamma_e}.
\end{align*}
Therefore,
\begin{equation}
 \langle\varphi_e,1\rangle_{\Gamma_e}
 =\langle q_e^+-q_e^-,1\rangle_{\Gamma_e}
 =(s-1)\langle q_e^-,1\rangle_{\Gamma_e}=0.
 \label{proofconverseoutercharge}
\end{equation}
Thus, $\boldsymbol\varphi\in\mathcal X$.

We now calculate the distributional Laplacian, including its sign. Let $\vartheta\in C_c^\infty(\mathbb R^2)$ and choose $R$ so large that $\overline{\Omega_e}$ and $\operatorname{supp}\vartheta$ lie in $B_R$, with $\vartheta=0$ near $\partial B_R$. On each of $\Omega_i,\Omega_s,D_R$, the second Green identity is
\begin{align*}
 \int_D v\Delta\vartheta
 =\int_D(\Delta v)\vartheta
    +\int_{\partial D}\big(v\,\frac{\partial}{\partial \nu_D}\vartheta
                              -\vartheta\,\frac{\partial}{\partial \nu_D}v\big)\,ds.
\end{align*}
Upon summation, the outer-circle terms vanish. The sum of the two boundary contributions on each $\Gamma_j$ equals
\begin{align*}
 \int_{\Gamma_j}
   \left[v|_-\,\frac{\partial}{\partial \nu_j}\vartheta-\vartheta q_j^-
         -v|_+\,\frac{\partial}{\partial \nu_j}\vartheta+\vartheta q_j^+\right]\,ds &=\int_{\Gamma_j}(v|_--v|_+)\frac{\partial}{\partial \nu_j}\vartheta\,ds
          +\langle q_j^+-q_j^-,\vartheta|_{\Gamma_j}\rangle_{\Gamma_j}\\
 &=\langle\varphi_j,\vartheta|_{\Gamma_j}\rangle_{\Gamma_j}.
\end{align*}
Here continuity eliminates the first integral. Only the exterior volume term is nonzero, since $\Delta v=f$ there.
Define the interface distribution by
$\langle\varphi_j\delta_{\Gamma_j},\vartheta\rangle
=\langle\varphi_j,\vartheta|_{\Gamma_j}\rangle_{\Gamma_j}$.
Since $\Delta F=f$, the preceding calculation gives
\begin{equation}
 \Delta(v-F)=\varphi_i\delta_{\Gamma_i}+\varphi_e\delta_{\Gamma_e}.
 \label{proofdistributionjump}
\end{equation}

The fundamental-solution identity also gives, for each $j$,
\begin{align}
 \langle\Delta\mathcal S_j[\varphi_j],\vartheta\rangle
 &=\left\langle\varphi_j,\,
       \left.\int_{\mathbb R^2}\Phi(x-\cdot)\Delta\vartheta(x)\,dx
                 \right|_{\Gamma_j}\right\rangle_{\Gamma_j}\\
                 &=\langle\varphi_j,\vartheta|_{\Gamma_j}\rangle_{\Gamma_j}
   =\langle\varphi_j\delta_{\Gamma_j},\vartheta\rangle.
 \label{prooflayerdistribution}
\end{align}
This identity follows first for smooth densities and subsequently by the single-layer mapping properties and density in $H^{-1/2}$. Subtracting \eqref{prooflayerdistribution} from \eqref{proofdistributionjump} shows that
\begin{align*}
 H=v-F-\mathcal S_i[\varphi_i]-\mathcal S_e[\varphi_e] \qquad \text{ and }
 \qquad \Delta H=0\quad\text{in }\mathcal D'(\mathbb R^2).
\end{align*}
The function $H\in H^1_{\rm loc}$ has a smooth entire harmonic representative by Weyl's lemma. Equations \eqref{pde}, \eqref{prooflayerdecay}, and the decay of $F$ imply $H(x)=O(|x|^{-1})$. For any fixed $x$, the mean-value property on arbitrarily large circles, therefore, gives
\begin{align*}
 H(x)=\lim_{R\to+\infty}\frac1{2\pi}
 \int_0^{2\pi}H(x+R(\cos\theta,\sin\theta))\,d\theta=0.
\end{align*}
Thus, $H=0$, proving the representation for $v$. Applying Lemma~\ref{layerflux} to this representation shows that its densities solve \eqref{proofblocksystem}. Invertibility from Lemmas~\ref{layercharges} and~\ref{layerkernel} identifies them uniquely with
$(\lambda I-\mathbb K^*)^{-1}\mathbf g_F$, and hence, $v=u_\ell$. Finally, the derivative jump \eqref{proofdensityjump} also shows directly that two density pairs representing the same field must coincide. This proves both uniqueness assertions.
\end{proof}

\section{Proofs of Section~\ref{conformal}}\label{app:coefficients}
\subsection{Proof of Lemma~\ref{geometry}}\label{proof:geometry}
\begin{proof}
For $w\ne v$,
\begin{align*}
 \frac{\Psi(w)-\Psi(v)}{w-v}
 =1-\varepsilon\frac{w+v}{w^2v^2} \qquad \text{ and } \qquad
 \left|\varepsilon\frac{w+v}{w^2v^2}\right|\le2|\varepsilon|/\rho^3<1.
\end{align*}
The quotient is nonzero, so $\Psi(w)\ne\Psi(v)$. Moreover
\begin{align*}
 \Psi'(w)=1-2\varepsilon w^{-3}\qquad \text{ and } \qquad
 |\Psi'(w)|\ge1-2|\varepsilon|\rho^{-3}>0.
\end{align*}
Thus, the restrictions to the two circles are analytic embeddings. The image of the annulus is their intervening region: the injective map takes its two boundary components to these curves, and its extension to infinity satisfies $\lim_{|w|\to+\infty}\Psi(w)/w=1$. Consequently the outer image bounds the inner image.

For $|w|,|v|\ge\rho$ the logarithmic series converges absolutely, and
\begin{align*}
 \log\left(1-\varepsilon\frac{w+v}{w^2v^2}\right)
 &=-\sum_{k\ge1}\frac{\varepsilon^k}{k}
       \frac{(w+v)^k}{w^{2k}v^{2k}} =-\sum_{k\ge1}\sum_{j=0}^k
       \frac{\varepsilon^k}{k}\binom{k}{j}
       w^{-(2k-j)}v^{-(k+j)}.
\end{align*}
Setting $n=2k-j$, $m=k+j$ gives $n+m=3k$ and
$\binom{k}{j}=\binom{k}{n-k}$, proving~\eqref{grunsky}.
Interchanging $n,m$ leaves this coefficient unchanged. Since
$\sqrt{nm}\le(n+m)/2=3k/2$ and $\binom{k}{j}\le2^k$,
\begin{equation}
 |G_{nm}|\le\frac32(2|\varepsilon|\rho^{-3})^k
             =\frac32 b_*^{n+m}.
 \label{grunskyentrybound}
\end{equation}
Consequently,
\begin{align*}
 \sum_{n,m}|G_{nm}|^2
 \le\frac94\left(\frac{b_*^2}{1-b_*^2}\right)^2<+\infty.
\end{align*}
hence, $G$ is Hilbert--Schmidt, and its real symmetric matrix defines a compact self-adjoint operator. Its row and column sums satisfy
\begin{align*}
 \sup_n\sum_m|G_{nm}|
 \le\sup_n\frac{3b_*^{n+1}}{2(1-b_*)}
 =\frac{3b_*^2}{2(1-b_*)}<\frac34<1.
\end{align*}
The Schur estimate, therefore, proves the norm bound. Finally,
$k\le n,m\le2k$ implies $n/2\le m\le2n$. Consequently
\begin{align*}
 \|Ge_n\|^2\le\frac94 b_*^{2n}
 \sum_{m\ge\lceil n/2\rceil}b_*^{2m}
 \le\frac9{4(1-b_*^2)}b_*^{3n}.
\end{align*}
\end{proof}

\subsection{Proof of Lemma~\ref{coretraces}}\label{proof:coretraces}
\begin{proof}
Lemma~\ref{geometry} makes $\gamma_i$ a smooth parametrization with $h_i>0$. Pullback by $\gamma_i$ and multiplication by $h_i$ or $h_i^{-1}$ preserve the relevant Sobolev spaces. For even, constant-free traces, Fourier orthogonality gives the coefficients in~\eqref{corenormalization}. The Sobolev norms of their partial sums satisfy
\begin{align*}
 \left\|\sum_{n=1}^N\frac{\mathbf v_n}{\sqrt n}\cos(n\theta)
       \right\|_{H^{1/2}(\mathbb T)}^2
 &\asymp\sum_{n=1}^N\frac{(1+n^2)^{1/2}}{n}|\mathbf v_n|^2
 \asymp\sum_{n=1}^N|\mathbf v_n|^2,\\
 \left\|\sum_{n=1}^N\sqrt n\,\mathbf k_n\cos(n\theta)
       \right\|_{H^{-1/2}(\mathbb T)}^2
 &\asymp\sum_{n=1}^N\frac{n}{(1+n^2)^{1/2}}|\mathbf k_n|^2
 \asymp\sum_{n=1}^N|\mathbf k_n|^2.
\end{align*}
The constants are independent of $N$. Applying the same estimates to tails proves convergence of~\eqref{coreinverses} for every pair of square-summable vectors. Conversely, finite Sobolev norm implies square summability of the normalized coefficients. Fourier uniqueness and the absence of a constant or sine mode show that these series are the inverse maps, proving the isomorphism and norm assertions.

The harmonic extension in~\eqref{coreDN} satisfies
$\|U_v\|_{H^1(\Omega_i)}\le C\|v\|_{H^{1/2}(\Gamma_i)}$.
For a bounded extension operator $\mathcal E:H^{1/2}(\Gamma_i)\to H^1(\Omega_i)$, its weak normal derivative is characterized by
\begin{align*}
 \langle\Lambda_i v,\psi\rangle_{\Gamma_i}
 =\int_{\Omega_i}\nabla U_v\cdot\nabla(\mathcal E\psi)\,dx, \qquad \text{ and } \qquad
 |\langle\Lambda_i v,\psi\rangle_{\Gamma_i}|
 \le C\|v\|_{H^{1/2}(\Gamma_i)}\|\psi\|_{H^{1/2}(\Gamma_i)}.
\end{align*}
If two extensions have the same trace, their difference lies in $H^1_0(\Omega_i)$, and the weak harmonic equation makes the integral of its gradient against $\nabla U_v$ zero. Thus, the definition is independent of the extension and gives a bounded normal trace. Taking the constant extension of $\psi=1$ gives
\begin{align*}
 \langle\widetilde{\Lambda_i v},1\rangle_{\mathbb T}
 =\langle\Lambda_i v,1\rangle_{\Gamma_i}=0.
\end{align*}
Finally $\Psi(\overline w)=\overline{\Psi(w)}$ and $h_i(-\theta)=h_i(\theta)$. For $v\in\mathcal V_i$, reflection preserves both the core Dirichlet problem and its boundary data. Uniqueness, therefore, gives $U_v(\overline z)=U_v(z)$, and its metric normal derivative is even. Together with the zero-charge identity this proves $\Lambda_i v\in\mathcal Q_i$ and the asserted boundedness.
\end{proof}

\subsection{Proof of Lemma~\ref{core}}\label{proof:core}
\begin{proof}
Recall the Faber polynomials $\mathcal F_n$ defined in~\eqref{faberdefinition}. Substitution $z=\Psi(v)$ and expansion for $|w|>|v|$ give
\begin{align*}
 \mathcal F_n(\Psi(v))=v^n+\sum_{m\ge1}n\beta_{nm}v^{-m}.
\end{align*}
To verify this generating calculation, write
\begin{align*}
 \log\frac{\Psi(w)-\Psi(v)}w
 =\log(1-v/w)-\sum_{n,m\ge1}\beta_{nm}w^{-n}v^{-m}
\end{align*}
and equate coefficients of $w^{-n}$. For a finite sequence $a$ define the complex-valued harmonic polynomial
\begin{align*}
 U_a(z)=\sum_n\frac{a_n}{\sqrt n\rho^n}
                      \operatorname{Re}\mathcal F_n(z).
\end{align*}
The real part here acts on the polynomial, not on $a_n$. On $v=re^{i\theta}$, reality of $\beta_{nm}$ gives
\begin{align*}
 U_a(\Psi(re^{i\theta}))
 =\sum_n\frac{a_n}{\sqrt n\rho^n}
     \left(r^n\cos n\theta+\sum_m n\beta_{nm}r^{-m}\cos m\theta\right).
\end{align*}
The trace $v_a=U_a|_{\Gamma_i}$ has no constant or sine coefficient, so $v_a\in\mathcal V_i$. Its $m$th coefficient $(\mathscr D_i v_a)_m$ is
\begin{align*}
 \sqrt m\left(\frac{a_m}{\sqrt m}
       +\sum_n\frac{a_n n\beta_{nm}\rho^{-m}}{\sqrt n\rho^n}\right)
 =a_m+\sum_n\sqrt{mn}\beta_{nm}\rho^{-m-n}a_n
 =[(I+G)a]_m.
\end{align*}
With $h_i=|\rho e^{i\theta}\Psi'(\rho e^{i\theta})|$, the metric flux is
$h_i\frac{\partial}{\partial \nu_i}U_a=\rho\partial_r(U_a\circ\Psi)|_{r=\rho}$.
Hence, the $m$th normalized flux coefficient $(\mathscr N_i\Lambda_i v_a)_m$ is
\begin{align*}
 a_m-\sum_n\frac{mn\beta_{nm}\rho^{-m-n}}{\sqrt m\sqrt n}a_n
 =[(I-G)a]_m.
\end{align*}
Thus, the coordinate identities are
\begin{align*}
 \mathscr D_i v_a=(I+G)a,\qquad
 \mathscr N_i\Lambda_i v_a=(I-G)a \qquad \text{ and } 
 \qquad L(I+G)a=(I-G)a.
\end{align*}

To extend this calculation beyond finite sequences, fix $\mathbf v\in\ell^2$ and choose finitely supported vectors $a^{(N)}$ satisfying $\lim_{N\to+\infty}a^{(N)}=(I+G)^{-1}\mathbf v$ in $\ell^2$. Then,
\begin{align*}
 \lim_{N\to+\infty}v_{a^{(N)}}
 =\lim_{N\to+\infty}\mathscr D_i^{-1}(I+G)a^{(N)}
 =\mathscr D_i^{-1}\mathbf v
       \quad\text{in }H^{1/2}(\Gamma_i).
\end{align*}
Boundedness of $\Lambda_i$ and $\mathscr N_i$ gives
\begin{align*}
 L\mathbf v
 =\lim_{N\to+\infty}\mathscr N_i\Lambda_i v_{a^{(N)}}
 =\lim_{N\to+\infty}(I-G)a^{(N)}
 =(I-G)(I+G)^{-1}\mathbf v.
\end{align*}
This proves the operator identity on all of $\ell^2$.

Finally $\operatorname{Sp}(G)\subset[-\gamma,\gamma]$ and
\begin{align*}
 \frac{1-\gamma}{1+\gamma}
 \le\frac{1-t}{1+t}\le\frac{1+\gamma}{1-\gamma}
 \quad \text{for } -\gamma\le t\le\gamma.
\end{align*}
Functional calculus gives $L=L^*$, these same lower and upper quadratic-form bounds, and $L^{-1}=(I+G)(I-G)^{-1}$. Since $\Psi(\overline w)=\overline{\Psi(w)}$, reflection preserves the domain and the Dirichlet problem. Hence, even boundary data have an even harmonic extension, which proves invariance of the cosine sector.
\end{proof}

\subsection{Proof of Lemma~\ref{smoothdipole}}\label{proof:smoothdipole}
\begin{proof}
Write $\chi(y)=\chi_0(|y|)$. For $|z-z_0|>b$, angular averaging gives
\begin{align*}
 (\Phi*\chi)(z-z_0) =\int_0^b\chi_0(r)r\int_0^{2\pi}
                  \Phi(z-z_0-re^{i\theta})\,d\theta\,dr =2\pi\int_0^b\chi_0(r)r\,dr\,\Phi(z-z_0)
 =\Phi(z-z_0).
\end{align*}
Here the inner integral uses harmonicity in the averaging disk, and the last equality uses $\int\chi=1$. Since $R$ and $\varepsilon$ are real, $z_0=\Psi(R)$ is real. Differentiation of convolution, or integration by parts against the compactly supported bump, now gives
\begin{align*}
 F=\Phi*\partial_{x_1}\chi(\cdot-z_0)
   =\partial_{x_1}(\Phi*\chi)(\cdot-z_0)
   =\frac1{2\pi}\frac{x_1-z_0}{|z-z_0|^2}
\end{align*}
outside that disk, which is~\eqref{dipole}. Globally, convolution with the smooth bump makes $F$ smooth. Outside a large disk $F=O(|z|^{-1})$ and $\nabla F=O(|z|^{-2})$. Hence,
$\int_{|z|>B}|\nabla F|^2\le C\int_B^{+\infty} r^{-3}\,dr<+\infty$, and its energy on the remaining compact set is finite by smoothness.
\end{proof}

\subsection{Proof of Lemma~\ref{source}}\label{proof:source}
\begin{proof}
The source-disk separation and Lemma~\ref{smoothdipole} give~\eqref{dipole} near the closed shell. Differentiate the Grunsky expansion with respect to its second variable at $v=R$. The coefficient bound~\eqref{grunskyentrybound} ensures absolute, locally uniform convergence of the differentiated series for $\rho\le|w|<R$. Thus,
\begin{align*}
 &-\frac{\Psi'(R)}{\Psi(w)-\Psi(R)}+\frac1{w-R}
 =\sum_{n,m\ge1}m\beta_{nm}w^{-n}R^{-m-1},\\
 &\frac1{\Psi(w)-\Psi(R)}
 =\frac1{\Psi'(R)}
 \left[\frac1{w-R}-\sum_{n,m\ge1}m\beta_{nm}w^{-n}R^{-m-1}\right].
\end{align*}
Here $\Psi'(R)=1-2\varepsilon/R^3>0$. Expanding
$1/(w-R)=-\sum_{n\ge0}w^nR^{-n-1}$ and taking real parts in~\eqref{dipole} identifies the coefficients in~\eqref{incident}:
\begin{align*}
 A_n=-\frac{\sqrt n}{2\pi\Psi'(R)R^{n+1}} \quad \text{ and } \quad
 (D_{\rm inc})_n
 =-\frac{\sqrt n}{2\pi\Psi'(R)}
       \sum_{m\ge1}m\beta_{nm}R^{-m-1}
 =\sum_{m\ge1}t_nG_{nm}t_mA_m=(TGTA)_n.
\end{align*}
The last equality uses~\eqref{grunsky}. Since $t_n^\alpha=R^{-2n}$,
\begin{align*}
 |A_n|^2=\kappa nt_n^\alpha \qquad \text{ and } \qquad
 \sum_{n\ge1}|A_n|^2=\frac{\kappa R^{-2}}{(1-R^{-2})^2}<+\infty.
\end{align*}
Boundedness of $T$ and $G$ now gives $D_{\rm inc}=T(GTA)\in\ell^2$, and also $T^{-1}D_{\rm inc}=GTA\in\ell^2$. This proves~\eqref{A} and~\eqref{compatible}.
\end{proof}

\subsection{Proof of Lemma~\ref{coeffinterfaces}}\label{proof:coeffinterfaces}
\begin{proof}
Work first with finite Fourier sums. The trace estimates below justify passage to square-summable vectors. Continuity at the outer boundary gives $b+Tq=A+D_{\rm inc}+h$. In normalized flux coordinates, its exterior derivative is $A-D_{\rm inc}-h=2A-b-Tq$. Thus,
\begin{align*}
 2A-b-Tq=s_\ell(b-Tq), \qquad 
 (1+s_\ell)b+(1-s_\ell)Tq=2A \qquad \text{ and } \qquad 
 \frac{1+s_\ell}{1-s_\ell}b+Tq=\frac2{1-s_\ell}A.
\end{align*}
Since $(1+s_\ell)/(1-s_\ell)=\zeta$ and $2/(1-s_\ell)=1+\zeta$, this gives
\begin{equation}
 \zeta b+Tq=(1+\zeta)A.
 \label{outer}
\end{equation}
At the inner boundary, $L(Tb+q)=s_\ell(Tb-q)$. Using~\eqref{impedance}, multiplication by $I+G$ gives
\begin{align*}
 &(I-G)(Tb+q)=s_\ell(I+G)(Tb-q),\\
 &[(1-s_\ell)I-(1+s_\ell)G]Tb
  +[(1+s_\ell)I-(1-s_\ell)G]q=0.
\end{align*}
Division by $1-s_\ell$ yields
\begin{equation}
 (I-\zeta G)Tb+(\zeta I-G)q=0.
 \label{inner}
\end{equation}
Substituting $q=GTb+Cz$ into the inner equation, and using $CG=GC$, gives
\begin{align*}
 0&=(I-\zeta G)Tb+(\zeta I-G)(GTb+Cz) =(I-G^2)Tb+(\zeta I-G)Cz\\
  &=C[CTb+(\zeta I-G)z].
\end{align*}
Since $C\ge\sqrt{1-\gamma^2}I$, this is equivalent to
$CTb+(\zeta I-G)z=0$. The outer equation becomes
$(\zeta I+TGT)b+TCz=(1+\zeta)A$. These are exactly the two rows of
$(\zeta-J)(b,z)=(1+\zeta)g$.
\end{proof}

\subsection{Proof of Lemma~\ref{coeffmetric}}\label{proof:coeffmetric}
\begin{proof}
Use the factorization
\begin{align*}
 J=QUQ,\qquad Q=\operatorname{diag}(T,I) \qquad \text{ and } \qquad
 U=\begin{pmatrix}-G&-C\\-C&G\end{pmatrix}.
\end{align*}
Since $G$ and $C$ are commuting self-adjoint operators with $G^2+C^2=I$, one has $U^*=U$ and $U^2=I$. Consequently $J^*=J$. Its blocks are compact because $T$ and $G$ are compact. If $JX=0$, injectivity of $Q$ gives $UQX=0$, and multiplication by $U$ gives $QX=0$, hence, $X=0$.

For $X=(b,z)$, put $q=GTb+Cz$ and $v=-CTb+Gz$. The norm-preserving property of $U$ yields
\begin{align*}
 UQX=(-q,v),\qquad JX=(-Tq,v) \qquad \text{ and } \qquad
 \|q\|^2+\|v\|^2=\|Tb\|^2+\|z\|^2.
\end{align*}
Using $W^2=I-T^2$ and $MX=(Wb,Wq)$, we obtain
\begin{align*}
 \|X\|^2-\|JX\|^2
 &=\|b\|^2+\|z\|^2-\|Tq\|^2-\|v\|^2\\
 &=\|b\|^2-\|Tb\|^2+\|q\|^2-\|Tq\|^2
 =\|Wb\|^2+\|Wq\|^2=\|MX\|^2.
\end{align*}
Equality of these quadratic forms proves $I-J^2=M^*M$. Moreover,
\begin{align*}
 M^{-1}=\begin{pmatrix}
 W^{-1}&0\\-C^{-1}GTW^{-1}&C^{-1}W^{-1}
 \end{pmatrix},\qquad
 W\ge\sqrt{1-\rho^2}I,\qquad C\ge\sqrt{1-\gamma^2}I.
\end{align*}
Thus, $M^{-1}$ is bounded, and
\begin{align*}
 \|JX\|^2=\|X\|^2-\|MX\|^2
 \le(1-\|M^{-1}\|^{-2})\|X\|^2,
\end{align*}
which gives $\|J\|<1$. Finally, for $y=\operatorname{Im}\zeta>0$,
\begin{align*}
 y\|X\|^2=\operatorname{Im}\langle X,(\zeta-J)X\rangle
 \le\|X\|\|(\zeta-J)X\|.
\end{align*}
Hence, $\zeta-J$ is injective with closed range. Its adjoint $\overline\zeta-J$ has zero kernel by the same estimate, so the range is also dense. Therefore, $\zeta-J$ is onto and $\|(\zeta-J)^{-1}\|\le1/y$.
\end{proof}

\subsection{Proof of Lemma~\ref{coeffenergy}}\label{proof:coeffenergy}
\begin{proof}
The Jacobian of a conformal map is $|\Psi'|^2$, while squared gradient length is multiplied by $|\Psi'|^{-2}$, so these factors cancel in the Dirichlet integral. For the $n$th shell mode, put $d_n=t_nq_n$ and
$v_n=(b_nr^n+d_nr^{-n})\cos(n\theta)/\sqrt n$. Angular orthogonality removes products with different indices, and
\begin{align*}
 \int_0^{2\pi}(|\partial_rv_n|^2+r^{-2}|\partial_\theta v_n|^2)\,d\theta
 &=\pi n\bigl(|b_nr^{n-1}-d_nr^{-n-1}|^2
                   +|b_nr^{n-1}+d_nr^{-n-1}|^2\bigr)\\
 &=2\pi n(|b_n|^2r^{2n-2}+|d_n|^2r^{-2n-2}).
\end{align*}
Thus,
\begin{align*}
 \int_\rho^1\int_0^{2\pi}|\nabla v_n|^2r\,d\theta\,dr
 =\pi[|b_n|^2(1-\rho^{2n})+|d_n|^2(\rho^{-2n}-1)] =\pi(1-t_n^2)(|b_n|^2+|q_n|^2).
\end{align*}
Summing proves the first identity in~\eqref{power}, by convergence of the partial sums in energy. Also
\begin{align*}
 \binom{Wb}{Wq}=M\binom bz=(1+\zeta)MV,
 \qquad \text{ and } \qquad E=\pi\ell|1+\zeta|^2\|MV\|^2
        =2\pi y\langle V,(I-J^2)V\rangle.
\end{align*}
The inequalities in~\eqref{coercive} follow from
$0\le\|JV\|^2\le\|J\|^2\|V\|^2$.
Subtracting the outer equation~\eqref{outer} from continuity gives
\begin{align*}
 h=(1-\zeta)b+\zeta A-D_{\rm inc},
\end{align*}
which proves the coefficient formula in~\eqref{observe}. The exterior expansion and its energy are
\begin{align*}
 (u_\ell-F)(\Psi(re^{i\theta}))=
 \sum_{n\ge1}\frac{h_nr^{-n}}{\sqrt n}\cos(n\theta) \qquad \text{ and } \qquad
 S_r=\pi\sum_{n\ge1}r^{-2n}|h_n|^2.
\end{align*}

In the last energy calculation one can integrate each exterior mode explicitly:
\begin{align*}
 2\pi n|h_n|^2\int_r^{+\infty} s^{-2n-1}\,ds
 =\pi r^{-2n}|h_n|^2.
\end{align*}
\end{proof}

\subsection{Proof of Lemma~\ref{coeffreconstruction}}\label{proof:coeffreconstruction}
\begin{proof}
Use the shell series specified before Lemma~\ref{coeffinterfaces} and the exterior coefficient vector of Lemma~\ref{coeffenergy}.
The vectors $b,z,q,h$ lie in $\ell^2$ by their bounded-operator formulas. The shell energy calculation and the trace vectors $Tb+q$, $b+Tq$, therefore, give an $H^1$ shell field and $H^{1/2}$ boundary values. Extend the inner trace harmonically into $\Omega_i$ using Lemma~\ref{core}, and define the exterior as $F$ plus the displayed scattered series. The latter has finite energy and, by Cauchy--Schwarz for $r\ge2$,
\begin{align*}
 \left|\sum_{n\ge1}\frac{h_nr^{-n}}{\sqrt n}\cos n\theta\right|
 \le\|h\|\left(\sum_{n\ge1}r^{-2n}/n\right)^{1/2}
 \le C\|h\|r^{-1}.
\end{align*}
For the gradient, the corresponding sum is
\begin{align*}
 \sum_{n\ge1}nr^{-2n-2}=\frac{r^{-4}}{(1-r^{-2})^2}=O(r^{-4}).
\end{align*}
Taking its square root gives the gradient estimate $O(r^{-2})$. The conformal map satisfies $\Psi(w)=w+O(w^{-2})$, with $\lim_{|w|\to+\infty}\Psi'(w)=1$. These relations give the required decay estimates in physical coordinates.

The zeroth shell mode could initially be written $B_0+Q_0\log r$. Its inner flux integral is $2\pi Q_0$. The core is harmonic, so transmission gives
\begin{align*}
 0=\int_{\Gamma_i}\frac{\partial}{\partial \nu_i}u|_-\,ds
   =s_\ell\int_{\Gamma_i}\frac{\partial}{\partial \nu_i}u|_+\,ds
   =2\pi s_\ell Q_0 \qquad \text{ and } 
 \qquad Q_0=0.
\end{align*}
The scattered exterior has neither a logarithmic nor a constant mode by decay. The incident zeroth trace coefficient is $F_0$, so outer continuity gives $B_0=F_0$. Thus, constants match at both interfaces.
The nonconstant traces match by the equations in Lemma~\ref{coeffinterfaces}, and the flux differences are zero in $H^{-1/2}$. For a compactly supported test function $\vartheta$, regional integration by parts consequently yields
\begin{align*}
 -\int_{\mathbb R^2}a_\ell\nabla u\cdot\nabla\vartheta
 =\int_{\mathbb R^2}f\vartheta
 +\langle s_\ell q_i^+-q_i^-,\vartheta\rangle
 +\langle q_e^+-s_\ell q_e^-,\vartheta\rangle
 =\int f\vartheta,
\end{align*}
where $q_j^\pm=\frac{\partial}{\partial \nu_j}u|_\pm$. This verifies the weak PDE. Theorem~\ref{BIE}, whose uniqueness proof is already in the energy class, identifies this field with its layer solution. No second formal uniqueness argument is needed. Reflection of that unique solution is again a solution with the same source, so it equals the original field and has no sine component.
\end{proof}

\subsection{Proof of Lemma~\ref{densityrecovery}}\label{proof:densityrecovery}
\begin{proof}
Theorem~\ref{BIE} identifies each density with its normal derivative jump. Writing $q_j^\pm=\frac{\partial}{\partial \nu_j}u|_\pm$, the two transmission conditions give
\begin{align*}
 \varphi_i=q_i^+-q_i^-=q_i^+-s_\ell q_i^+=(1-s_\ell)q_i^+
 \qquad \text{ and } \qquad
 \varphi_e=q_e^+-q_e^-=s_\ell q_e^--q_e^-=(s_\ell-1)q_e^-.
\end{align*}
Differentiate the shell series to obtain
\begin{align*}
 r\partial_r(u\circ\Psi)
 =\sum_{n\ge1}\sqrt n(b_nr^n-t_nq_nr^{-n})\cos n\theta.
\end{align*}
At $r=\rho$ and $r=1$ this gives, respectively,
\begin{align*}
 h_iq_i^+=\sum_n\sqrt n(t_nb_n-q_n)\cos n\theta
 \qquad \text{ and } \qquad h_eq_e^-=\sum_n\sqrt n(b_n-t_nq_n)\cos n\theta,
\end{align*}
which proves~\eqref{densityinner}--\eqref{densityouter}. For each finite partial sum,
\begin{align*}
 \int_{\Gamma_i}\varphi_i\,ds
 =(1-s_\ell)\sum_n\sqrt n(t_nb_n-q_n)
                        \int_0^{2\pi}\cos n\theta\,d\theta=0.
\end{align*}
The outer charge is zero by the identical calculation with
$(s_\ell-1)(b_n-t_nq_n)$. To pass to the infinite sums, observe that
\begin{align*}
 \sum_n(1+n^2)^{-1/2}
        |\sqrt n(t_nb_n-q_n)|^2
 \le\sum_n|t_nb_n-q_n|^2<+\infty.
\end{align*}
The outer series satisfies the same estimate. Smooth multiplication by $h_j^{-1}$ and smooth change of boundary coordinates preserve $H^{-1/2}$. Pairing with $1\in H^{1/2}$ is continuous. The two charge identities, therefore, hold for the actual densities.
\end{proof}

\subsection{Proof of Lemma~\ref{tracecoordinates}}\label{proof:tracecoordinates}
\begin{proof}
For each $n$, the two boundary values satisfy
\begin{align*}
 (v_e)_n=b_n+t_nq_n \qquad \text{ and } \qquad (v_i)_n=t_nb_n+q_n.
\end{align*}
Subtracting $t_n$ times the second equation from the first, and subsequently reversing their roles, gives
\begin{align*}
 (v_e)_n-t_n(v_i)_n=(1-t_n^2)b_n \qquad \text{ and }
 \qquad (v_i)_n-t_n(v_e)_n=(1-t_n^2)q_n.
\end{align*}
This is~\eqref{traceinverse}. Since $1-t_n^2\ge1-\rho^2$,
\begin{align*}
 \|b\|\le\frac{\|v_e\|+\rho\|v_i\|}{1-\rho^2} \qquad \text{ and } \qquad
 \|q\|\le\frac{\|v_i\|+\rho\|v_e\|}{1-\rho^2}.
\end{align*}
Thus, the inverse maps are bounded. Lemma~\ref{coeffreconstruction} fixes the constant mode at $F_0$ and excludes logarithmic and sine modes, so the recovered vectors describe the full shell field. Applying Lemma~\ref{coeffinterfaces} to its transmission conditions gives
\begin{align*}
 (\zeta-J)\binom{b}{C^{-1}(q-GTb)}=(1+\zeta)g.
\end{align*}
Since $C^{-1}$ is bounded and $1+\zeta=2/(1-s_\ell)\ne0$, division by $1+\zeta$ gives the asserted vector $V$ and equation $(\zeta-J)V=g$.
\end{proof}

\section{Proofs of Section~\ref{sec:uniform}}\label{app:resolvents}
\subsection{Proof of Lemma~\ref{relative}}\label{proof:relative}
\begin{proof}
Initially define $B=GT^{-1}$ on finite sequences. The column bound gives
\begin{align*}
 \sum_n\|Be_n\|^2\le C_G^2\sum_n(\eta/\rho)^{2n}<+\infty,
\end{align*}
so $B$ has a Hilbert--Schmidt extension with $BT=G$ and $TB^*=G$. Put $D=I-C=G(I+C)^{-1}G$ and
\begin{equation}
 K=\begin{pmatrix}-B^*G(I+C)^{-1}T&-B^*\\GT&-D\end{pmatrix}.
 \label{K}
\end{equation}
The identities $BT=G$ and $TB^*=G$ hold first on coordinate vectors and hence, on the whole space by continuity. Each block of $K$ contains a compact factor $B^*$ or $G$, or is $-D$ with $D=G(I+C)^{-1}G$. Thus, $K$ is compact. The four product blocks are
\begin{align*}
 \begin{aligned}
 (J_0K)_{11}&=-TGT,& (J_0K)_{12}&=TD,\\
 (J_0K)_{21}&=TB^*G(I+C)^{-1}T=DT,
 & (J_0K)_{22}&=TB^*=G.
 \end{aligned}
\end{align*}
Since $D=I-C$, these are the blocks of $J-J_0$. Hence, $J=J_0(I+K)$ as an identity of bounded operators. If $(I+K)v=0$, then, $Jv=0$, and Theorem~\ref{physical} gives $v=0$. Since identity plus compact is Fredholm of index zero, zero kernel implies zero cokernel and surjectivity. The bounded inverse theorem gives $(I+K)^{-1}\in\mathcal B(\HH)$.

The $n$th row of $B^*=T^{-1}G$ has norm at most $C_G\eta^n/\rho^n$. Right multiplication by the bounded factors in~\eqref{K} preserves this bound up to constants. The remaining blocks satisfy
\begin{align*}
 \|e_n^*GT\|\le\rho C_G\eta^n \qquad \text{ and } \qquad
 \|e_n^*D\|\le\frac{\gamma C_G}{1+\sqrt{1-\gamma^2}}\eta^n.
\end{align*}
More explicitly, write $K_{ij}$ for its four blocks. There is a constant $C$ independent of $n$ with
\begin{align*}
 \|e_n^*K_{1j}\|\le C(\eta/\rho)^n
 \qquad \text{ and } \qquad \|e_n^*K_{2j}\|\le C\eta^n\quad(j=1,2).
\end{align*}
Define the matrix of $L_s$ by multiplying the $n$th row in each component by $\rho^{-sn}$. Then,
\begin{align*}
 \|L_s\|_{\rm HS}^2
 \le2C^2\sum_{n\ge1}\left[
       (\eta\rho^{-1-s})^{2n}+(\eta\rho^{-s})^{2n}\right]<+\infty.
\end{align*}
This constructs $L_s$ as a bounded Hilbert--Schmidt operator directly from its matrix. Equality $K=\TT^sL_s$ holds entrywise and, therefore, on finite sequences. Continuity extends it to $\HH$. No unbounded inverse of $J_0$ or $\TT$ has been applied to an arbitrary vector.
\end{proof}

\subsection{Proof of Lemma~\ref{multiplier}}\label{proof:multiplier}
\begin{proof}
At index $n$, the orthonormal vectors $(e_n,-e_n)/\sqrt2$ and $(e_n,e_n)/\sqrt2$ diagonalize $J_0$ with eigenvalues $t_n$ and $-t_n$. They also diagonalize $\TT^s$ with eigenvalue $t_n^s$. The required operator norm is consequently
\begin{align*}
 \sup_{n\ge1,\,\epsilon\in\{-1,1\}}
          \frac{a t_n^s}{|\zeta-\epsilon t_n|}.
\end{align*}
For $t=t_n\le a/2$, the reverse triangle inequality gives
$|\zeta-\epsilon t|\ge a-t\ge a/2$, so the quotient is at most $2t^s\le2^{1-s}a^s$.
For $a/2<t<2a$, the imaginary part gives $|\zeta-\epsilon t|\ge y$, so it is at most $2^sa^{1+s}/y$.
For $t\ge2a$, the denominator is at least $t-a\ge t/2$, giving $2at^{s-1}$. If $0<s\le1$ this is at most $2a(2a)^{s-1}=2^sa^s$, whereas if $s\ge1$ it is at most $2a\rho^{s-1}$. Taking the supremum over the three regions proves the estimate.
\end{proof}

\subsection{Proof of Lemma~\ref{uniformresolvent}}\label{proof:uniformresolvent}
\begin{proof}
Along the exact path, $a\asymp\ell^\tau$ and $y\asymp\ell$. The strict coupling condition allows
$s>\max(0,1/\tau-1)$ with $\eta<\rho^{1+s}$. Lemmas~\ref{relative} and~\ref{multiplier} show
\begin{align*}
 \|\zeta(\zeta-J_0)^{-1}K\|
 =O(\ell^{\tau\min(s,1)}+\ell^{\tau(1+s)-1})=o(1).
\end{align*}
To derive the comparison identity without an unbounded inverse, put $R_0=(\zeta-J_0)^{-1}$ and multiply
$(\zeta-J_0-J_0K)V=g$ by $R_0$. Since $R_0J_0=\zeta R_0-I$, this gives
\begin{align*}
 V-(\zeta R_0-I)KV=R_0g=V_0.
\end{align*}
Thus,
\begin{align*}
 F_\zeta V=V_0 \qquad \text{ and } \qquad
 \lim_{\ell\downarrow0}F_{\zeta(\ell)}
 =I+K-\lim_{\ell\downarrow0}\zeta(\ell)(\zeta(\ell)-J_0)^{-1}K=I+K
\end{align*}
in operator norm. Write $H=I+K$ and $R_\zeta=F_\zeta-H$. For sufficiently small loss, $\|H^{-1}R_\zeta\|\le1/2$, so
\begin{align*}
 F_\zeta^{-1}=(I+H^{-1}R_\zeta)^{-1}H^{-1} \qquad \text{ and } \qquad
 \|F_\zeta^{-1}\|\le2\|H^{-1}\|,
\end{align*}
which implies that
\begin{align*}
 F_\zeta^{-1}-H^{-1}=-F_\zeta^{-1}R_\zeta H^{-1}
 \qquad \text{ and } \qquad
 \lim_{\ell\downarrow0}\|F_{\zeta(\ell)}^{-1}-H^{-1}\|=0.
\end{align*}
Also $\|F_\zeta\|\le\|H\|+o(1)$, and consequently
\begin{align*}
 \|F_\zeta\|^{-1}\|V_0\|\le\|V\|
 \le\|F_\zeta^{-1}\|\|V_0\|.
\end{align*}
The two energy forms satisfy
\begin{align*}
 2\pi(1-\|J\|^2)y\|V\|^2\le E\le2\pi y\|V\|^2
 \qquad \text{ and } \qquad
 2\pi(1-\rho^2)y\|V_0\|^2\le E_0\le2\pi y\|V_0\|^2.
\end{align*}
Combining these inequalities proves both uniform comparisons, with constants independent of $g$.
\end{proof}

\subsection{Proof of Lemma~\ref{leadingresolvent}}\label{proof:leadingresolvent}
\begin{proof}
Use $F_\zeta$ from Lemma~\ref{uniformresolvent} and set $H=I+K$.
Fix $g$ and a positive sequence $(\ell_j)$ with $\lim_{j\to+\infty}\ell_j=0$ and $\lim_{j\to+\infty}\|V_0(\ell_j)\|=+\infty$. Throughout this proof, evaluate all loss-dependent quantities at $\ell=\ell_j$ and set $\zeta=\zeta(\ell_j)$ and $N_\zeta=\|V_0(\ell_j)\|$.
At every fixed modal index $n$,
\begin{align*}
 \lim_{j\to+\infty}
 \begin{pmatrix}\zeta&t_n\\t_n&\zeta\end{pmatrix}^{-1}
 =\begin{pmatrix}0&t_n^{-1}\\t_n^{-1}&0\end{pmatrix}.
\end{align*}
Thus, every fixed coordinate of $V_0$ remains bounded, while $\lim_{j\to+\infty}N_\zeta=+\infty$. Writing $w_\zeta=V_0/N_\zeta$, we obtain $\lim_{j\to+\infty}(w_\zeta)_n=0$ at every fixed index. Since $\|w_\zeta\|=1$, finite-coordinate approximation gives, for $v\in\HH$, 
\begin{align*}
 \lim_{j\to+\infty}\langle w_\zeta,v\rangle=0.
\end{align*}
The operator $H^{-1}-I=-H^{-1}K$ is compact, so
$\lim_{j\to+\infty}\|(H^{-1}-I)w_\zeta\|=0$. Therefore,
\begin{align*}
 \frac{V-V_0}{N_\zeta}
 =(F_\zeta^{-1}-H^{-1})w_\zeta+(H^{-1}-I)w_\zeta \qquad \text{ and } \qquad
 \lim_{j\to+\infty}\frac{\|V-V_0\|}{N_\zeta}=0.
\end{align*}
In particular, $\lim_{j\to+\infty}\|V\|^2/N_\zeta^2=1$. The resolvent equations give $JV=\zeta V-g$ and $J_0V_0=\zeta V_0-g$, hence,
\begin{align*}
 \frac{\|JV\|+\|J_0V_0\|}{N_\zeta}
 \le|\zeta|\left(\frac{\|V\|}{N_\zeta}+1\right)
       +\frac{2\|g\|}{N_\zeta} \qquad \text{ and } \qquad
 \lim_{j\to+\infty}\frac{\|JV\|+\|J_0V_0\|}{N_\zeta}=0.
\end{align*}
Finally,
\begin{align*}
 \lim_{j\to+\infty}\frac E{2\pi yN_\zeta^2}
 &=\lim_{j\to+\infty}
   \frac{\|V\|^2-\|JV\|^2}{N_\zeta^2}=1 \qquad \text{ and } \qquad
 \lim_{j\to+\infty}\frac{E_0}{2\pi yN_\zeta^2}
 &=\lim_{j\to+\infty}
   \left(1-\frac{\|J_0V_0\|^2}{N_\zeta^2}\right)=1.
\end{align*}
Consequently, $\lim_{j\to+\infty}E(\ell_j)/E_0(\ell_j)=1$, completing~\eqref{leading}.
\end{proof}

\section{Proofs of Section~\ref{sec:power}}\label{app:sums}
\subsection{Proof of Lemma~\ref{circularsum}}\label{proof:circularsum}
\begin{proof}
The block form of $J_0$ in~\eqref{comparator} and the forcing $g=(A,0)$ give the block equation below. Since $\operatorname{Im}\zeta>0$ and $t_n>0$, its determinant $\zeta^2-t_n^2$ does not vanish.
Thus,
\begin{align*}
 \begin{pmatrix}\zeta&t_n\\t_n&\zeta\end{pmatrix}(V_0)_n
 =\binom{A_n}0
 \qquad \text{ and } \qquad
 \begin{pmatrix}\zeta&t_n\\t_n&\zeta\end{pmatrix}^{-1}
 =\frac1{\zeta^2-t_n^2}
       \begin{pmatrix}\zeta&-t_n\\-t_n&\zeta\end{pmatrix}.
\end{align*}
Multiplying by $(A_n,0)^{\mathsf T}$ and taking the squared Euclidean norm gives both identities in~\eqref{V0}. Moreover, for real $t$,
\begin{align*}
 \frac1{|\zeta-t|^2}+\frac1{|\zeta+t|^2}
 =\frac{|\zeta+t|^2+|\zeta-t|^2}{|\zeta^2-t^2|^2}
 =\frac{2(a^2+t^2)}{|\zeta^2-t^2|^2}.
\end{align*}
Since $I-J_0^2=\operatorname{diag}(I-T^2,I-T^2)$, the energy in~\eqref{comparator} becomes
\begin{align*}
 E_0 =2\pi y\sum_{n\ge1}(1-t_n^2)\|(V_0)_n\|^2 =2\pi y\sum_{n\ge1}|A_n|^2(1-t_n^2)
       \frac{a^2+t_n^2}{|\zeta^2-t_n^2|^2}.
\end{align*}
The preceding identity with $t=t_n$ now gives~\eqref{positivesum}. All summands are nonnegative because $0<t_n<1$. The sum is finite for each $y>0$: the resolvent of the self-adjoint operator $J_0$ satisfies $\|V_0\|\le y^{-1}\|g\|$, and $0<I-J_0^2\le I$.
\end{proof}

\subsection{Proof of Lemma~\ref{geometricsums}}\label{proof:geometricsums}
\begin{proof}
Throughout, $1-\rho^2\le1-t_n^2\le1$ and $|A_n|^2=\kappa n t_n^\alpha$. For $0<z<1$, differentiating the geometric series gives
\begin{equation}
 \sum_{n>N}nz^n
 =z^{N+1}\left(\frac{N+1}{1-z}+\frac z{(1-z)^2}\right).
 \label{geometrictail}
\end{equation}
For $z>1$, on the other hand,
\begin{align*}
 Nz^N\le\sum_{n=1}^N nz^n
 =z^N\sum_{j=0}^{N-1}(N-j)z^{-j}
 \le\frac{Nz^N}{1-z^{-1}}.
\end{align*}
For $z=1$ the sum is $N(N+1)/2$, while for $0<z<1$ its partial sums are bounded above by $z/(1-z)^2$ and below by $z$.
Given small $q>0$, choose $N$ with $t_{N+1}<q\le t_N$. Then,
\begin{align*}
 q\le t_N<q/\rho \qquad \text{ and } \qquad N\asymp\log(1/q),
\end{align*}
and splitting at $N$ yields
\begin{align*}
 \sum_n\frac{nt_n^\alpha}{(q+t_n)^2} \asymp\sum_{n\le N}nt_n^{\alpha-2}
             +q^{-2}\sum_{n>N}nt_n^\alpha \asymp H_\alpha(q).
\end{align*}
Indeed the tail is $\asymp Nq^{\alpha-2}$ by~\eqref{geometrictail}, while the head is respectively $\asymp Nq^{\alpha-2}$, $\asymp N^2$, or $\asymp1$ as $\alpha<2$, $\alpha=2$, or $\alpha>2$.
\end{proof}

\subsection{Proof of Lemma~\ref{sectorialsums}}\label{proof:sectorialsums}
\begin{proof}
If $|x|\le C y$, then,
\begin{align*}
 (x\mp t)^2+y^2\le2t^2+(2C^2+1)y^2
 \qquad \text{ and } \qquad
 t^2+y^2\le2(t\mp x)^2+(2C^2+1)y^2.
\end{align*}
Thus, $|\zeta\mp t|^2\asymp t^2+y^2\asymp(t+y)^2$, with constants depending only on $C$. Inserting this into~\eqref{positivesum} and applying Lemma~\ref{geometricsums} with $q=y$ proves $E_0\asymp yH_\alpha(y)$.
\end{proof}

\subsection{Proof of Lemma~\ref{tangentialsums}}\label{proof:tangentialsums}
\begin{proof}
Since $q=|x|$, the two denominators in~\eqref{positivesum} are $(q-t_n)^2+y^2$ and $(q+t_n)^2+y^2$, possibly in reversed order. Separating the matching-branch term at the chosen index $m$ gives~\eqref{nearestpoledecomposition} exactly. By symmetry, we may, therefore, assume $x=q>0$. Eventually $0<y\le q$. Throughout this range, the opposite branch satisfies
\begin{align*}
 (q+t_n)^2\le|\zeta+t_n|^2\le2(q+t_n)^2,
\end{align*}
so its full contribution is $\asymp yH_\alpha(q)$. Let $m$ minimize $|q-t_m|$. It is one of the two indices bracketing $q$, and hence,
$t_m\asymp q$, $m\asymp\log(1/q)$. The matching-branch term at $m$ is, therefore,
\begin{align*}
 \frac{y m t_m^\alpha(1-t_m^2)}{(q-t_m)^2+y^2}
 \asymp\frac{yq^\alpha\log(1/q)}{\mathscr{D}(q)^2+y^2}.
\end{align*}
Here and below a tied nearest index can be selected arbitrarily.
To bound all other terms uniformly, note that for $n\ne m$,
\begin{align*}
 |q-t_n|\ge\tfrac12|t_n-t_m|
\end{align*}
by the triangle inequality and the nearest-index property. If $q/2\le t_n\le2q$, geometric spacing gives
$|t_n-t_m|\ge(1-\rho)\max(t_n,t_m)\ge(1-\rho)t_n$.
If $t_n<q/2$ or $t_n>2q$, directly $|q-t_n|\ge\tfrac12\max(q,t_n)$. Combining these cases gives, for some $c_\rho>0$ and for $n\ne m$
\begin{align*}
 |q-t_n|\ge c_\rho\max(q,t_n).
\end{align*}
This also covers ties. Therefore, the remaining matching-branch terms have sum at most
\begin{align*}
 Cy\sum_{n\ne m}\frac{nt_n^\alpha}{(q+t_n)^2}
 \le CyH_\alpha(q).
\end{align*}
The opposite branch supplies a lower bound for $E_0^{\mathrm{rem}}$, while the last estimate and the opposite-branch upper bound give
\begin{align*}
 yH_\alpha(q)\lesssim E_0^{\mathrm{rem}}
 \lesssim yH_\alpha(q).
\end{align*}
Together with the estimate for the matching term in~\eqref{nearestpoleenergy}, this proves~\eqref{nearestpolebounds}. Adding the two nonnegative contributions in~\eqref{nearestpoledecomposition} yields~\eqref{envelopes}.
\end{proof}

\subsection{Proof of Lemma~\ref{anglelower}}\label{proof:anglelower}
\begin{proof}
Choose $n$ with $t_{n+1}<a\le t_n$. Then, $a\le t_n<a/\rho$, $n\asymp\log(1/a)$, and
$|\zeta\mp t_n|\le a+t_n\le(1+\rho^{-1})a$.
Keeping just one positive term in~\eqref{positivesum} gives
\begin{align*}
 E_0\ge C y\frac{nt_n^\alpha}{a^2}
 \ge C'ya^{\alpha-2}\log(1/a),
\end{align*}
as claimed for the specified range $0<\alpha<2$.
\end{proof}

\subsection{Proof of Lemma~\ref{powersector}}\label{proof:powersector}
\begin{proof}
If $p\ge1$, relations~\eqref{crossings} and~\eqref{coordinate} give $|x|\lesssim y\asymp\ell$. Proposition~\ref{comparison} and Lemma~\ref{sums}, therefore, yield
\begin{align*}
 E\asymp\ell H_\alpha(\ell)=
 \begin{cases}
 \ell^{\alpha-1}\log(1/\ell),&0<\alpha<2,\\
 \ell\log(1/\ell)^2,&\alpha=2,\\
 \ell,&\alpha>2.
 \end{cases}
\end{align*}
For $\alpha<1$ the negative power diverges, whereas at $\alpha=1$ the logarithm diverges. For $1<\alpha<2$ the positive power dominates the logarithm and the expression tends to zero. The last two lines also tend to zero. This proves the first classification.
\end{proof}

\subsection{Proof of Lemma~\ref{powercrossings}}\label{proof:powercrossings}
\begin{proof}
For $0<p<1$, the exact coordinate~\eqref{crossings} gives
$x(\ell)\sim(c/2)\ell^p$ and $y(\ell)\sim\ell/2$. In particular,
$q(\ell)\sim(|c|/2)\ell^p$ and $\lim_{\ell\downarrow0}y(\ell)/q(\ell)=0$.
Differentiating that coordinate directly yields
\begin{align*}
 \zeta'(\ell)=\frac{2(cp\ell^{p-1}+i)}
                    {(2-c\ell^p-i\ell)^2} \qquad \text{ and } 
 \qquad
 \lim_{\ell\downarrow0}\ell^{1-p}\zeta'(\ell)=\frac{cp}{2}.
\end{align*}
Taking real parts gives $x'(\ell)\sim(cp/2)\ell^{p-1}$. Since
$\operatorname{sgn}x(\ell)=\operatorname{sgn}c$ for sufficiently small $\ell$,
\begin{align*}
 q'(\ell)=\operatorname{sgn}(c)x'(\ell)
 \sim\frac{|c|p}{2}\ell^{p-1}>0.
\end{align*}
Thus, $q$ is strictly increasing on some interval $(0,\ell_*)$, with
$\lim_{\ell\downarrow0}q(\ell)=0$. For all sufficiently large $n$, both
$t_n$ and $2\rho t_n/(1+\rho)$ lie in its image, giving unique crossings
$\ell_n$ and gaps $\widehat\ell_n$ in this interval. Finally,
$q(\ell)\asymp\ell^p$ gives
$\ell_n\asymp\widehat\ell_n\asymp t_n^{1/p}$.
\end{proof}

\subsection{Proof of Lemma~\ref{powerenvelopes}}\label{proof:powerenvelopes}
\begin{proof}
Use $q=|x|\asymp\ell^p$, $y\asymp\ell$, and the quantities $B_\ell,P_\ell$ in the statement.
The uniform comparison and~\eqref{envelopes}, together with
$\mathscr{D}^2+y^2\ge y^2$, give
\begin{align*}
 B_\ell\lesssim E\lesssim B_\ell+P_\ell.
\end{align*}
At a crossing $\mathscr{D}=0$, so $E\asymp B_\ell+P_\ell$.
At a gap $q=2\rho t_n/(1+\rho)$, the adjacent distances are
\begin{align*}
 t_n-q=\frac{1-\rho}{1+\rho}t_n \qquad \text{ and }\qquad
 q-t_{n+1}=\frac{\rho(1-\rho)}{1+\rho}t_n.
\end{align*}
All other poles are farther away, so $\mathscr{D}\asymp q$. Its extra term is
$O(yq^{\alpha-2}\log(1/q))=O(yH_\alpha(q))$ for every $\alpha>0$, and thus, $E\asymp B_\ell$ on gaps.
\end{proof}

\section{Proofs of Section~\ref{sec:weighted}}\label{app:weighted}
\subsection{Proof of Lemma~\ref{twosided}}\label{proof:twosided}
\begin{proof}
Symmetry and~\eqref{gdecay} imply
$|G_{nm}|\le C_G\min(\eta^n,\eta^m)\le C_G\eta^{(n+m)/2}$. Consequently
\begin{align*}
 \sum_{n,m}|\rho^{-\sigma n}G_{nm}\rho^{-\sigma m}|^2
 \le C_G^2\sum_{n,m}(\eta\rho^{-2\sigma})^{n+m}<+\infty.
\end{align*}
This constructs $G_\sigma$ as a Hilbert--Schmidt matrix. Equality on finite sequences and continuity give $G=T^\sigma G_\sigma T^\sigma$. The one-sided weighted matrices are also Hilbert--Schmidt, since
\begin{align*}
 \sum_n\|Ge_n\|^2t_n^{-2\sigma}
 \le C_G^2\sum_n(\eta\rho^{-\sigma})^{2n}<+\infty.
\end{align*}
Taking adjoints gives the left-weighted version. The hypothesis implies
$\eta<\rho^{2\sigma}<\rho^\sigma$. Define
\begin{align*}
 D_\sigma=(T^{-\sigma}G)(I+C)^{-1}(GT^{-\sigma}).
\end{align*}
This is a bounded compact self-adjoint operator, and
$T^\sigma D_\sigma T^\sigma=G(I+C)^{-1}G=D$.
All negative powers here denote the bounded extensions just constructed, not an inverse applied to arbitrary vectors. The four blocks of $\TT^\sigma P_\sigma\TT^\sigma$ are
\begin{align*}
 \begin{pmatrix}
 -T(T^\sigma G_\sigma T^\sigma)T&T(T^\sigma D_\sigma T^\sigma)\\
 (T^\sigma D_\sigma T^\sigma)T&T^\sigma G_\sigma T^\sigma
 \end{pmatrix}
 =\begin{pmatrix}-TGT&TD\\DT&G\end{pmatrix}=J-J_0.
\end{align*}
Each block of $P_\sigma$ contains a compact factor. Its diagonal blocks are self-adjoint and $(TD_\sigma)^*=D_\sigma T$. This proves all the assertions.
\end{proof}

\subsection{Proof of Lemma~\ref{weightedlimit}}\label{proof:weightedlimit}
\begin{proof}
Let $Z_{\sigma,\zeta}=\TT^\sigma(\zeta-J_0)^{-1}\TT^\sigma$. Its eigenmultipliers are $t^{2\sigma}/(\zeta\mp t)$. For $t\le2a$ their modulus is at most $(2a)^{2\sigma}/y$, and for $t>2a$ at most $2t^{2\sigma-1}$, uniformly bounded because $\sigma\ge1/2$. Coordinatewise convergence and uniform boundedness yield the strong limit
\begin{align*}
 Z_{\sigma,0}=\begin{pmatrix}0&T^{2\sigma-1}\\T^{2\sigma-1}&0\end{pmatrix}.
\end{align*}
Here is the convergence argument explicitly. Let $\Pi_N$ project onto the first $N$ indices in both components. For any fixed $v\in\HH$,
\begin{align*}
 \|(Z_{\sigma,\zeta}-Z_{\sigma,0})v\|
 \le\|(Z_{\sigma,\zeta}-Z_{\sigma,0})\Pi_Nv\|
       +C\|(I-\Pi_N)v\|.
\end{align*}
For fixed $N$, the scalar multiplier limits give
\begin{align*}
 \lim_{\zeta\to0}\|(Z_{\sigma,\zeta}-Z_{\sigma,0})\Pi_Nv\|=0
 \qquad \text{ and } \qquad
 \lim_{N\to+\infty}\|(I-\Pi_N)v\|=0.
\end{align*}
Taking these limits in this order proves
$\lim_{\zeta\to0}\|(Z_{\sigma,\zeta}-Z_{\sigma,0})v\|=0$.
Here and below, $\zeta$ approaches zero within the upper half-plane under the bound in Lemma~\ref{weightedlimit}.
Choose finite-rank $P_{\sigma,N}$ with $\lim_{N\to+\infty}\|P_\sigma-P_{\sigma,N}\|=0$. Then,
\begin{align*}
 \|(Z_{\sigma,\zeta}-Z_{\sigma,0})P_\sigma\|
 \le C\|P_\sigma-P_{\sigma,N}\|
       +\|(Z_{\sigma,\zeta}-Z_{\sigma,0})P_{\sigma,N}\|.
\end{align*}
Strong convergence is uniform on the unit ball of a fixed finite-dimensional range, so the last term tends to zero for fixed $N$. It follows that
$\lim_{\zeta\to0}\|F_{\sigma,\zeta}-F_{\sigma,0}\|=0$.
\end{proof}

\subsection{Proof of Lemma~\ref{weightedkernel}}\label{proof:weightedkernel}
\begin{proof}
We must check the limiting kernel in this weighted space. Let $X=(X_1,X_2)\in\HH$ satisfy $F_{\sigma,0}X=0$. Its two component equations are
\begin{align*}
 \begin{aligned}
 X_1&=T^{2\sigma-1}(D_\sigma TX_1+G_\sigma X_2) \qquad \text{ and } \qquad X_2&=T^{2\sigma}(-G_\sigma TX_1+D_\sigma X_2).
 \end{aligned}
\end{align*}
If $\sigma\ge1$, the first equation factors as $T^\sigma$ times the bounded vector $T^{\sigma-1}(D_\sigma TX_1+G_\sigma X_2)$, and the second likewise factors as $T^\sigma$ times a bounded vector. Thus, $X\in\operatorname{Ran}\TT^\sigma$.
For $1/2\le\sigma<1$, use the second equation
\begin{align*}
 X_2=T^{2\sigma}(-G_\sigma TX_1+D_\sigma X_2),
\end{align*}
to set $u_2=T^\sigma(-G_\sigma TX_1+D_\sigma X_2)\in\ell^2$, so $X_2=T^\sigma u_2$. In the first equation, the two weighted identities give
\begin{align*}
 T^{2\sigma-1}D_\sigma T=T^\sigma(T^{-1}D)T^{1-\sigma} \qquad \text{ and } 
 \qquad T^{2\sigma-1}G_\sigma T^\sigma=T^\sigma(T^{-1}G).
\end{align*}
It, therefore, becomes
\begin{equation}
 X_1=T^\sigma[(T^{-1}D)T^{1-\sigma}X_1+(T^{-1}G)u_2].
 \label{kernelrange}
\end{equation}
The bracket is a bounded expression: $\eta<\rho^{2\sigma}\le\rho$ implies bounded $T^{-1}G$ and $T^{-1}D=(T^{-1}G)(I+C)^{-1}G$. Identities with negative powers can first be checked on finite sequences and subsequently extended by boundedness. Thus, in both cases $X=\TT^\sigma u$ for $u\in\HH$.

Multiply the kernel identity by $J_0$. Direct multiplication gives
$J_0Z_{\sigma,0}=-\TT^{2\sigma}$, so with $X=\TT^\sigma u$,
\begin{align*}
 J_0\TT^\sigma u=-\TT^{2\sigma}P_\sigma\TT^\sigma u
                  =-\TT^\sigma Pu \qquad \text{ and } \qquad \TT^\sigma(J_0+P)u=\TT^\sigma Ju=0.
\end{align*}
Injectivity of $\TT^\sigma$ and $J$ gives $u=X=0$.
\end{proof}

\subsection{Proof of Lemma~\ref{weightedinverse}}\label{proof:weightedinverse}
\begin{proof}
Lemma~\ref{weightedkernel} gives $\ker F_{\sigma,0}=\{0\}$, while Lemma~\ref{weightedlimit} gives norm convergence.
Since $Z_{\sigma,0}P_\sigma$ is compact, $F_{\sigma,0}$ is Fredholm of index zero and hence, bijective with bounded inverse. Norm convergence and the Neumann-series argument give, for small $\zeta$,
\begin{align*}
 \|F_{\sigma,\zeta}^{-1}\|\le2\|F_{\sigma,0}^{-1}\| \qquad \text{ and } 
 \qquad
 F_{\sigma,\zeta}^{-1}-F_{\sigma,0}^{-1}
 =F_{\sigma,\zeta}^{-1}(F_{\sigma,0}-F_{\sigma,\zeta})F_{\sigma,0}^{-1}.
\end{align*}
Hence, $\lim_{\zeta\to0}\|F_{\sigma,\zeta}^{-1}-F_{\sigma,0}^{-1}\|=0$ under the hypotheses of Lemma~\ref{weightedlimit}.
To obtain the claimed identity, start from $(\zeta-J_0)V=g+PV$, multiply by $(\zeta-J_0)^{-1}$, and subsequently by $\TT^\sigma$. Factorization~\eqref{weightedP} yields
\begin{align*}
 \TT^\sigma V=\TT^\sigma V_0+
       Z_{\sigma,\zeta}P_\sigma\TT^\sigma V.
\end{align*}
Rearrangement and inversion prove~\eqref{weightedcomparisonidentity}. Finally,
$F_{\sigma,0}^{-1}-I=F_{\sigma,0}^{-1}Z_{\sigma,0}P_\sigma$ is bounded times compact, hence, compact. At $\sigma=1/2$, $Z_{\sigma,0}$ itself is not compact, but only its product with $P_\sigma$ is needed.
\end{proof}

\section{Proofs of Section~\ref{sec:localisation}}\label{app:localisation}
\subsection{Proof of Lemma~\ref{localsufficiency}}\label{proof:localsufficiency}
\begin{proof}
Set $\sigma_c=1/(2\tau)$. Given $r>r_c$, choose $\sigma>\sigma_c$ sufficiently close that $\eta<\rho^{2\sigma}$ and $\rho^{-\sigma}<r$. This is possible because the inequalities in the hypotheses are strict. Then,
\begin{align*}
 \frac{a^{2\sigma}}{y}=O(\ell^{2\sigma\tau-1}) \qquad \text{ and }
 \qquad
 \lim_{\ell\downarrow0}\frac{a^{2\sigma}}{y}=0.
\end{align*}
Fix a positive sequence $(\ell_j)$ with $\lim_{j\to+\infty}\ell_j=0$ and $\lim_{j\to+\infty}E(\ell_j)=+\infty$. In all subsequent limits in this proof, loss-dependent quantities are evaluated at $\ell=\ell_j$. Write
\begin{align*}
 \|\TT^\sigma V_0\|^2
 =\sum_{t_n\le2a}t_n^{2\sigma}\|(V_0)_n\|^2
    +\sum_{t_n>2a}t_n^{2\sigma}\|(V_0)_n\|^2=:I_<+I_>.
\end{align*}
Since $E\asymp y\|V_0\|^2$,
\begin{align*}
 \frac{I_<}{E}\le\frac{(2a)^{2\sigma}\|V_0\|^2}{E}
 \le C\frac{a^{2\sigma}}y \qquad \text{ and }
 \qquad \lim_{j\to+\infty}\frac{I_<}{E}=0.
\end{align*}
For $t_n>2a$, both $|\zeta\mp t_n|\ge t_n-a\ge t_n/2$ and
$a^2+t_n^2\le5t_n^2/4$. Thus,~\eqref{V0} gives
$\|(V_0)_n\|^2\le20|A_n|^2/t_n^2$. The remaining weighted sum is at most
\begin{align*}
 C\sum_{t_n>2a}n t_n^e \qquad \text{ with } \qquad e=\alpha+2\sigma-2.
\end{align*}
The head-sum calculations in the proof of Lemma~\ref{geometricsums} give the three cases explicitly:
\begin{align*}
 I_>\le C
 \begin{cases}
 1,&e>0,\\
 \log(1/a)^2,&e=0,\\
 a^e\log(1/a),&e<0.
 \end{cases}
\end{align*}
When $e>0$, $\lim_{j\to+\infty}E(\ell_j)=+\infty$ gives $\lim_{j\to+\infty}I_>/E=0$. When $e<0$, necessarily $\alpha<2$, so~\eqref{lowera} and $E\asymp E_0$ give
\begin{align*}
 \frac{I_>}{E}\le C\frac{a^{\alpha+2\sigma-2}\log(1/a)}
                      {ya^{\alpha-2}\log(1/a)}
 =C\frac{a^{2\sigma}}y \qquad \text{ and }
 \qquad \lim_{j\to+\infty}\frac{I_>}{E}=0.
\end{align*}
When $e=0$, $\alpha=2-2\sigma\in(0,2)$ and the same lower bound gives
\begin{align*}
 \frac{I_>}{E}\le C\frac{a^{2\sigma}}y\log(1/a)
 \le C'\ell_j^{2\sigma\tau-1}\log(1/\ell_j) \qquad \text{ and }
 \qquad \lim_{j\to+\infty}\frac{I_>}{E}=0,
\end{align*}
because $2\sigma\tau-1>0$. Thus, $\lim_{j\to+\infty}\|\TT^\sigma V_0\|/\sqrt E=0$.
Lemma~\ref{weightedcomparison} now gives
\begin{align*}
 \frac{\|\TT^\sigma V\|}{\sqrt E}
 \le\|F_{\sigma,\zeta}^{-1}\|
          \frac{\|\TT^\sigma V_0\|}{\sqrt E} \qquad \text{ and }
 \qquad
 \lim_{j\to+\infty}\frac{\|\TT^\sigma V\|}{\sqrt E}=0.
\end{align*}
Since $D_{\rho^{-\sigma}}=T^\sigma$, the observation formula implies
\begin{align*}
 \frac{\sqrt{S_{\rho^{-\sigma}}}}{\sqrt{\pi E}}
 \le |1-\zeta^2|\frac{\|\TT^\sigma V\|}{\sqrt E}
       +\frac{|\zeta|\|A\|+\|D_{\rm inc}\|}{\sqrt E}.
\end{align*}
Consequently,
\begin{align*}
 \lim_{j\to+\infty}\frac{S_{\rho^{-\sigma}}(\ell_j)}{E(\ell_j)}=0.
\end{align*}
For $r>\rho^{-\sigma}$, each $r^{-2n}\le\rho^{2\sigma n}$, so
$S_r\le S_{\rho^{-\sigma}}$. This proves the asserted gradient convergence.

For local convergence of the potential itself, let $K$ be any compact subset of $\OO_{r_c}$. Choose $r_c<r'<r_-<r_+$ with
$\Psi^{-1}(K)\subset\{r_-<|w|<r_+\}$. Let $v=(u_\ell-F)\circ\Psi$. Orthogonality gives
\begin{align*}
 \int_{r_-}^{r_+}\int_0^{2\pi}|v|^2r\,d\theta\,dr
 &=\pi\sum_n\frac{|h_n|^2}{n}\int_{r_-}^{r_+}r^{1-2n}\,dr\\
 &\le\frac{\pi(r_+^2-r_-^2)}2
                \sum_n|h_n|^2(r')^{-2n}
 =\frac{r_+^2-r_-^2}{2}S_{r'}.
\end{align*}
The physical Jacobian is bounded on this compact annulus, and its physical gradient energy is at most $S_{r'}$ by conformal invariance. Therefore,
\begin{align*}
 \frac{\|u_\ell-F\|_{H^1(K)}}{\sqrt E}
 \le C_K\sqrt{S_{r'}/E},
 \end{align*}
 which gives
 \begin{align*}
 \lim_{j\to+\infty}\frac{\|u_{\ell_j}-F\|_{H^1(K)}}{\sqrt{E(\ell_j)}}=0
 \qquad\text{and}\qquad
 \lim_{j\to+\infty}\frac{\|F\|_{H^1(K)}}{\sqrt{E(\ell_j)}}=0.
\end{align*}
The last estimate is valid even when $K$ meets the smooth source support. The triangle inequality proves 
\begin{align*}
    \lim_{j\to+\infty}u_{\ell_j}/\sqrt{E(\ell_j)}=0 \quad \text{ in } H^1_{\rm loc}(\OO_{r_c}).
\end{align*}
\end{proof}

\subsection{Proof of Lemma~\ref{localsharpness}}\label{proof:localsharpness}
\begin{proof}
Throughout this proof, the limit is $n\to+\infty$ along the loss sequence specified in each case. Quantities written without an explicit loss argument are evaluated at $\ell=\ell_n$.
For $p<1$ on~\eqref{crossings}, write $t=t_n$, $y_n=\operatorname{Im}\zeta(\ell_n)$, and $\epsilon=\operatorname{sgn}c$. Then, $\zeta(\ell_n)=\epsilon t+iy_n$, $y_n\asymp t^{1/p}$, and $\lim_{n\to+\infty}y_n/t_n=0$. The matching pole contribution to~\eqref{positivesum} is exactly
\begin{align*}
 Q_n=\frac{\pi|A_n|^2(1-t^2)}{y_n}
 \asymp\frac{nt^\alpha}{y_n}.
\end{align*}
Lemma~\ref{tangentialsums} bounds every remaining contribution, including the opposite branch at index $n$, by $Cy_nH_\alpha(t)$. Dividing by $Q_n$ gives
\begin{align*}
 \frac{E_0-Q_n}{Q_n}
 \le C\frac{y_n^2H_\alpha(t)}{nt^\alpha}.
\end{align*}
Since $\log(1/t)=n\log(1/\rho)$, this is bounded by
\begin{equation}
 C\begin{cases}
 (y_n/t)^2,&0<\alpha<2,\\
 n(y_n/t)^2,&\alpha=2,\\
 t^{2/p-\alpha}/n,&\alpha>2.
 \end{cases}
 \label{dominance}
\end{equation}
For $\alpha<2$ the first bound tends to zero because $\lim_{n\to+\infty}y_n/t_n=0$. For $\alpha=2$, it is $O(nt_n^{2/p-2})$, with $\lim_{n\to+\infty}nt_n^{2/p-2}=0$. For $\alpha>2$ and $\alpha\le1/p$, the exponent $2/p-\alpha\ge1/p>0$, so the last bound tends to zero as well. Also
$Q_n\asymp nt_n^{\alpha-1/p}$ and $\lim_{n\to+\infty}Q_n=+\infty$ throughout the stated source range. Thus, $E_0\sim\pi|A_n|^2/y_n$, $\lim_{n\to+\infty}E_0(\ell_n)=+\infty$, and $\lim_{n\to+\infty}E(\ell_n)/E_0(\ell_n)=1$ by~\eqref{leading}.
The corresponding unweighted squared norm is
\begin{align*}
 \|V_0\|^2=\frac12\sum_m|A_m|^2
    \left(\frac1{|\zeta-t_m|^2}+\frac1{|\zeta+t_m|^2}\right).
\end{align*}
Since $1-\rho^2\le1-t_m^2\le1$, the preceding remainder bound also proves
\begin{align*}
 \|V_0\|^2\sim\frac{|A_n|^2}{2y_n^2}.
\end{align*}
For its first component at the matching index,~\eqref{V0} gives the exact expression
\begin{align*}
 |(V_0)_{1,n}|^2
 =\frac{|A_n|^2(t^2+y_n^2)}{y_n^2(4t^2+y_n^2)}
 \sim\frac{|A_n|^2}{4y_n^2}\sim\tfrac12\|V_0\|^2.
\end{align*}
Put $k_\zeta=\zeta A-D_{\rm inc}$, so $\|k_\zeta\|=O(1)$. Then,
\begin{align*}
 |h_n-(V_0)_{1,n}|
 &\le |1-\zeta^2|\|V-V_0\|
           +|\zeta|^2|(V_0)_{1,n}|+\|k_\zeta\|\\ 
           &=o(\|V_0\|)=o(|(V_0)_{1,n}|).
\end{align*}
Here $\lim_{n\to+\infty}\|V_0(\ell_n)\|=+\infty$, so the bounded source correction is negligible. Consequently $|h_n|^2\sim |A_n|^2/(4y_n^2)$. With $\beta=2\log r/\log(1/\rho)$ one has $r^{-2n}=t_n^\beta$, and the single exterior mode gives
\begin{equation}
 \frac{S_r}{E}\gtrsim\frac{t_n^\beta}{y_n}
 \asymp t_n^{\beta-1/p}.
 \label{tangentiallower}
\end{equation}

For $p\ge1$, instead set $\ell_n=t_n$. Then, $a\asymp y\asymp t_n$. Both factors $|\zeta\mp t_n|$ are bounded below by $y\asymp t_n$ and above by $a+t_n\asymp t_n$, so $|\zeta^2-t_n^2|\asymp t_n^2$. Formula~\eqref{V0}, therefore, gives
\begin{align*}
 |(V_0)_{1,n}|^2
 =\frac{|A_n|^2a^2}{|\zeta^2-t_n^2|^2}
 \asymp nt_n^{\alpha-2}.
\end{align*}
The sectorial sum estimate and the energy bounds give, for $\alpha\le1$,
\begin{align*}
 \|V_0\|^2\asymp E_0/y\asymp H_\alpha(t_n)
 \asymp nt_n^{\alpha-2} \qquad \text{ and } \qquad
 E(\ell_n)\asymp nt_n^{\alpha-1},
 \qquad\lim_{n\to+\infty}E(\ell_n)=+\infty.
\end{align*}
Thus, $|(V_0)_{1,n}|\ge c\|V_0\|$ with $\lim_{n\to+\infty}|(V_0(\ell_n))_{1,n}|=+\infty$, and the previous error estimate for $h_n-(V_0)_{1,n}$ again yields $|h_n|\ge c'\|V_0\|$ for large $n$. Keeping this single mode in $S_r$ gives
\begin{equation}
 S_r/E\gtrsim t_n^{\beta-1}.
 \label{sectorlower}
\end{equation}
Equations~\eqref{tangentiallower}--\eqref{sectorlower} prove divergence below $r_c$ and a positive lower limit at $r_c$.

For a ring $r<|w|<r_1$, direct integration of the same mode gives
\begin{align*}
 \int_{\Psi(\{r<|w|<r_1\})}|\nabla(u_\ell-F)|^2
 =\pi\sum_m|h_m|^2(r^{-2m}-r_1^{-2m})
 \ge\pi|h_n|^2r^{-2n}[1-(r/r_1)^{2n}].
\end{align*}
For fixed $r_1>r$ the last bracket tends to one. Given any proposed inner radius $r_0<r_c$ with $r_0>1$, choose $r_0<r<r_c$ and $r_1>r$. Let $K$ be the closure of the physical ring $\Psi(\{r<|w|<r_1\})$. Then, $K$ lies compactly in $\OO_{r_0}$, while the normalized gradient energy diverges on it along the same sequence. The fixed field satisfies $\lim_{n\to+\infty}\|\nabla F\|_{L^2(K)}/\sqrt{E(\ell_n)}=0$, so the reverse triangle inequality gives the same failure for $u_\ell/\sqrt E$. This proves sharpness of the open exterior region. A ring beginning at $r_c$ is not relatively compact in $\OO_{r_c}$, so the endpoint seminorm lower bound does not contradict local convergence there.
\end{proof}

\subsection{Proof of Lemma~\ref{localendpoint}}\label{proof:localendpoint}
\begin{proof}
All limits in this proof are taken as $n\to+\infty$ along the crossings~\eqref{crossings}. Loss-dependent quantities without an explicit argument are evaluated at $\ell=\ell_n$.
Now let $p<1$, $0<\alpha\le1/p$, and $\sigma=1/(2p)>1/2$. Write $y_n=\operatorname{Im}\zeta(\ell_n)$ on the exact crossing sequence. On crossings $a^{2\sigma}/y=O(1)$. We first prove
\begin{equation}
 \|\TT^\sigma V_0\|^2/E_0=O(1),
 \label{weightedbounded}
\end{equation}
and that all indices other than $n$ contribute $o(E_0)$ to the weighted squared norm. Write $N_n=\|V_0\|^2$ along this crossing sequence. The concentration calculation in the proof of Lemma~\ref{localsharpness} gives
\begin{align*}
 E_0\sim2\pi y_nN_n,\qquad
 \sum_{m\ne n}\|(V_0)_m\|^2=o(N_n) \qquad \text{ and }  \qquad
 \frac{t_n^{2\sigma}}{y_n}=O(1).
\end{align*}
For $t_m\le2a$ the full weighted contribution satisfies
\begin{align*}
 \frac1{E_0}\sum_{t_m\le2a}t_m^{2\sigma}\|(V_0)_m\|^2
 \le\frac{(2a)^{2\sigma}N_n}{E_0}=O(1).
\end{align*}
If $m=n$ is excluded, the same right-hand side contains $o(N_n)$ in place of $N_n$ and is $o(1)$. For $t_m>2a$, the denominator estimate in the proof of Lemma~\ref{localsufficiency} bounds the sum by $C\sum_{t_m>2a}m t_m^e$, with $e=\alpha+2\sigma-2$. Since $a\asymp t_n$, $\log(1/a)\asymp n$, and $E_0\asymp n t_n^{\alpha-2\sigma}$, division gives respectively
\begin{align*}
 O(E_0^{-1})\ (e>0),\qquad
 O(n t_n^{4\sigma-2})\ (e=0) \qquad \text{ and } \qquad
 O(t_n^{4\sigma-2})\ (e<0).
\end{align*}
For $e=0$, the algebra is
$n^2/(nt_n^{\alpha-2\sigma})=nt_n^{2\sigma-\alpha}=nt_n^{4\sigma-2}$, since $\alpha=2-2\sigma$. For $e<0$ it is
$nt_n^e/(nt_n^{\alpha-2\sigma})=t_n^{4\sigma-2}$.
All three cases tend to zero because $\lim_{n\to+\infty}E_0(\ell_n)=+\infty$ and $4\sigma-2>0$. This proves both assertions.

Put $w_n=\TT^\sigma V_0/\sqrt E$. By $\lim_{n\to+\infty}E(\ell_n)/E_0(\ell_n)=1$ and~\eqref{weightedbounded}, this sequence is bounded. Every fixed coordinate of $V_0$ stays bounded, while $\lim_{n\to+\infty}E(\ell_n)=+\infty$. Hence, $\lim_{n\to+\infty}(w_n)_m=0$ for each fixed index $m$. Finite-coordinate approximation proves $\lim_{n\to+\infty}\langle w_n,v\rangle=0$ for every $v\in\HH$. Set $C_\sigma=F_{\sigma,0}^{-1}-I$, which is compact by Lemma~\ref{weightedcomparison}. That lemma gives
\begin{align*}
 \frac{\TT^\sigma(V-V_0)}{\sqrt E}
 &=(F_{\sigma,\zeta}^{-1}-F_{\sigma,0}^{-1})w_n+C_\sigma w_n
.
\end{align*}
The first term has norm limit zero as $n\to+\infty$ by operator-norm convergence and boundedness, and the second by compactness and weak convergence. Thus,
\begin{align*}
 \lim_{n\to+\infty}\frac{\|\TT^\sigma(V-V_0)\|}{\sqrt E}=0.
\end{align*}
Since $D_{r_c}=T^\sigma$,~\eqref{observe} implies
\begin{align*}
 \frac{\|T^\sigma h-T^\sigma P_1V_0\|}{\sqrt E}
 &\le |1-\zeta^2|\frac{\|\TT^\sigma(V-V_0)\|}{\sqrt E}
 +|\zeta|^2\frac{\|\TT^\sigma V_0\|}{\sqrt E}
 +\frac{|\zeta|\|A\|+\|D_{\rm inc}\|}{\sqrt E}.
\end{align*}
Each term on the right has limit zero as $n\to+\infty$, so
\begin{align*}
 \lim_{n\to+\infty}\frac{\|T^\sigma h-T^\sigma P_1V_0\|}{\sqrt E}=0.
\end{align*}
The two vectors being compared are bounded after division by $\sqrt E$. The identity
$|\|v\|^2-\|w\|^2|\le\|v-w\|(\|v\|+\|w\|)$ and $\lim_{n\to+\infty}E(\ell_n)/E_0(\ell_n)=1$, therefore, yield
\begin{align*}
 \lim_{n\to+\infty}\left(
 \frac{S_{r_c}}E-\frac{\pi\|T^\sigma P_1V_0\|^2}{E_0}
 \right)=0.
\end{align*}
All nonmatching indices in the second quotient have just been shown to contribute $o(1)$. The matching index satisfies
\begin{align*}
 |(V_0)_{1,n}|^2
 =\frac{|A_n|^2(t_n^2+y_n^2)}{y_n^2(4t_n^2+y_n^2)}
 \sim\frac{|A_n|^2}{4y_n^2} \qquad \text{ and } \qquad
 E_0\sim\frac{\pi|A_n|^2}{y_n}.
\end{align*}
Consequently its contribution gives
\begin{align*}
 \frac{\pi\|T^\sigma P_1V_0\|^2}{E_0}
 \sim\frac{t_n^{2\sigma}}{4y_n}.
\end{align*}
Finally the exact coordinate implies
\begin{align*}
 \ell_n\sim(2t_n/|c|)^{1/p} \qquad \text{ and } \qquad
 y_n\sim\tfrac12(2t_n/|c|)^{1/p}.
\end{align*}
Since $2\sigma=1/p$, substitution yields
\begin{align*}
 \lim_{n\to+\infty}\frac{t_n^{2\sigma}}{4y_n}
 =
 \frac1{2(2/|c|)^{1/p}}
 =\frac12\left(\frac{|c|}{2}\right)^{1/p}.
\end{align*}
This proves the asserted endpoint limit. The source amplitude factor $\Psi'(R)^{-2}$ cancels through the common factor $|A_n|^2$.
\end{proof}

\section{Numerical details for Section~\ref{numerics}}\label{numericalappendix}
\paragraph*{Discretization and source.}
The figures use the first $N=80$ modes in each component of~\eqref{JM}, solving $(\zeta I-J_N)V_N=(A_N,0)^T$ in double precision. We evaluate $C_N=(I-G_N^2)^{1/2}$ with 24 nonconstant terms of its convergent binomial series and reconstruct $b_N,q_N,h_N$ through~\eqref{reconstruct}--\eqref{observe}. The energies are
\begin{align*}
 E_N=\pi\ell\sum_{m=1}^N(1-\rho^{2m})(|b_m|^2+|q_m|^2)
 \qquad \text{ and } \qquad S_{r,N}=\pi\sum_{m=1}^N r^{-2m}|h_m|^2.
\end{align*}
Thus, $S_{r,N}$ concerns the full unbounded exterior, not the displayed window. For each $\alpha$, set $R=\rho^{-\alpha/2}$ and use the source in Section~\ref{conformal} with a radial unit-mass bump of radius $b=0.04$, proportional to $\exp[-1/(1-|x|^2/b^2)]$ on $|x|<b$ and zero elsewhere. Its support is separated from the shell. Radial averaging preserves the exterior dipole potential, while the source itself remains smooth and fixed along each loss sequence.

\paragraph*{Sampling and field maps.}
For $p=0.75$, we use $n=6,\ldots,20$ and solve the exact material relation for
\begin{align*}
 \operatorname{Re}\zeta(\ell_n)=\rho^n \qquad \text{ and } \qquad
 \operatorname{Re}\zeta(\widehat\ell_n)=\frac{2\rho}{1+\rho}\rho^n.
\end{align*}
The inverse map $d+i\ell=2\zeta/(1+\zeta)$ is constrained by $d=\ell^p$. These are circular reference crossings and gaps, not computed noncircular eigenvalues. For $p=1.2$ we use $\ell_n=\rho^n$. Physical gradients are reconstructed analytically with the conformal derivative factor $1/\Psi'(w)$. The maps use 220 radial and 512 angular cells, with $\rho<|w|<1$ in the shell and $1<|w|<2.3$ outside. Each heat-map figure has one logarithmic inferno scale, floored at $10^{-7}$ of its maximum.

\paragraph*{Accuracy checks.}
Across all 120 regime samples, increasing the truncation from 80 to 104 modes changes the power by at most $8.64\times10^{-10}$ relatively, while an independent finite trace-system solve differs by at most $3.52\times10^{-10}$. At the intermittent crossing $n=14$, physical-coordinate quadrature agrees with the finite-mode shell and exterior energies to relative error below $10^{-13}$ on the tested meshes (40--160 radial and 128--512 angular nodes). The exterior comparator is $S_1-S_{2.3}$, matching the integration annulus. An independent boundary-integral reconstruction at $d+i\ell=0.001+0.003i$, $R=\sqrt2$, gives relative gradient errors $1.48\times10^{-2}$, $1.37\times10^{-4}$ and $1.17\times10^{-8}$ with 96, 192 and 384 nodes per interface at five fixed probes. 

\bibliographystyle{abbrv}
\bibliography{references}{}

\end{document}